\documentclass[11pt]{amsart}
\usepackage{mathrsfs,amssymb}
\newtheorem{theorem}{Theorem}[section]
\newtheorem{proposition}[theorem]{Proposition}
\newtheorem{lemma}[theorem]{Lemma}
\newtheorem{corollary}[theorem]{Corollary}
\theoremstyle{definition}
\newtheorem{definition}[theorem]{Definition}
\newtheorem{remark}[theorem]{Remark}

\begin{document}

\title{Ancient solutions to the Ricci flow in dimension $3$}
\author{Simon Brendle}
\address{Department of Mathematics, Columbia University, 2990 Broadway, New York, NY 10027, USA}
\begin{abstract}
It follows from work of Perelman that any finite-time singularity of the Ricci flow on a compact three-manifold is modeled on an ancient $\kappa$-solution.

We prove that the every noncompact ancient $\kappa$-solution in dimension $3$ is isometric to a family of shrinking cylinders (or a quotient thereof), or to the Bryant soliton. This confirms a conjecture of Perelman.
\end{abstract}
\maketitle
\tableofcontents

\section{Introduction}

A central problem in the study of geometric flows is to understand singularity formation. It turns out that singularities can often be modeled on ancient solutions; these are solutions which are defined on $(-\infty,T]$. The notion of an ancient solution was first introduced in Hamilton's work \cite{Hamilton4}. Perelman \cite{Perelman1} studied ancient solutions to the Ricci flow in dimension $3$ which are complete; non-flat; $\kappa$-noncollapsed; and have bounded and nonnegative curvature. These solutions are referred to as ancient $\kappa$-solutions. Perelman \cite{Perelman1} proved that every finite time singularity of the Ricci flow in dimension $3$ is modeled on an ancient $\kappa$-solution. Moreover, he proved a deep structure theorem for ancient $\kappa$-solutions. Roughly speaking, this theorem asserts that every noncompact ancient $\kappa$-solution with positive sectional curvature consists of a tube with a positively curved cap attached on one side. 

The purpose of this paper is to give a classification of all noncompact ancient $\kappa$-solutions in dimension $3$. In the first part of this paper, we classify all noncompact ancient $\kappa$-solutions with rotational symmetry:

\begin{theorem}
\label{main.thm.a}
Assume that $(M,g(t))$ is a three-dimensional ancient $\kappa$-solution which is noncompact and has positive sectional curvature. If $(M,g(t))$ is rotationally symmetric, then $(M,g(t))$ is isometric to the Bryant soliton up to scaling.
\end{theorem}

In the second part of this paper, we reduce the classification of noncompact ancient $\kappa$-solutions to the rotationally symmetric case:

\begin{theorem}
\label{main.thm.b}
Assume that $(M,g(t))$ is a three-dimensional ancient $\kappa$-solution which is noncompact and has positive sectional curvature. Then $(M,g(t))$ is rotationally symmetric. 
\end{theorem}

Theorem \ref{main.thm.b} extends our earlier work \cite{Brendle}, where we proved that the Bryant soliton is the only noncollapsed steady gradient Ricci soliton in dimension $3$. Note that, by work of Chen, every complete ancient solution to the Ricci flow in dimension $3$ has nonnegative sectional curvature (see \cite{Chen}, Corollary 2.4). 

Combining Theorem \ref{main.thm.a} and Theorem \ref{main.thm.b}, we obtain the following result: 

\begin{theorem}
\label{main.thm.c}
Assume that $(M,g(t))$ is a three-dimensional ancient $\kappa$-solution which is noncompact. Then $(M,g(t))$ is isometric to either a family of shrinking cylinders (or a quotient thereof), or to the Bryant soliton. 
\end{theorem}

Combining Theorem \ref{main.thm.c} with work of Perelman \cite{Perelman1}, we can draw the following conclusion:

\begin{corollary}
\label{first.singular.time}
Consider a solution to the Ricci flow on a compact three-manifold which forms a singularity in a finite time. Then, at the first singular time, the only possible blow-up limits are quotients of the round sphere $S^3$, quotients of the standard cylinder $S^2 \times \mathbb{R}$, and the Bryant soliton.
\end{corollary}

Let us sketch how Corollary \ref{first.singular.time} follows from Theorem \ref{main.thm.c}. Consider a smooth solution of the Ricci flow on a compact three-manifold which is defined on a finite time interval $[0,T)$ and becomes singular as $t \to T$. By work of Perelman \cite{Perelman1}, every blow-up limit as $t \to T$ is an ancient $\kappa$-solution. If a blow-up limit is compact with strictly positive sectional curvature, then the original flow will have positive sectional curvature for $t$ sufficiently close to $T$. A classical theorem of  Hamilton \cite{Hamilton1} then implies that the original flow becomes round as $t \to T$. If a blow-up limit is noncompact with strictly positive sectional curvature, then it must be the Bryant soliton by Theorem \ref{main.thm.c}. Finally, if a blow-up limit does not have strictly positive sectional curvature, then it must be a quotient of the cylinder by standard isometries.

Let us mention some related work. In \cite{Daskalopoulos-Hamilton-Sesum}, Daskalopoulos, Hamilton, and \v Se\v sum obtained a classification of all compact ancient solutions to the Ricci flow in dimension $2$. In \cite{Brendle-Choi}, it was shown that the bowl soliton is the only noncompact ancient solution to the mean curvature flow in $\mathbb{R}^3$ which is noncollapsed and strictly convex. Angenent, Daskalopoulos, and \v Se\v sum \cite{Angenent-Daskalopoulos-Sesum2} later obtained a classification of all compact ancient solutions to mean curvature flow in $\mathbb{R}^3$ which are noncollapsed and strictly convex. 

We now give an overview of the main ideas involved in the proof of Theorem \ref{main.thm.a} and Theorem \ref{main.thm.b}. 

In the first part of this paper, we classify noncompact ancient $\kappa$-solutions with rotational symmetry. In Section \ref{barrier.construction}, we set up a barrier argument for solutions to the Ricci flow with rotational symmetry. One important ingredient in our barrier construction are the steady gradient Ricci solitons with singularity at the tip which were found by Robert Bryant \cite{Bryant}. In Section \ref{asymptotics}, we study the asymptotic behavior of a noncompact ancient $\kappa$-solution with rotational symmetry. To that end, we focus on the cylindrical region, and carry out a spectral decomposition in Hermite polynomials. As in \cite{Angenent-Daskalopoulos-Sesum1} and \cite{Brendle-Choi}, a subtle point here is that we need to control certain error terms arising from the cutoff functions. In our work, this is done using barrier arguments. Using the spectral analysis, we obtain precise asymptotics for the solution in the cylindrical region. Combining these estimates with the barrier arguments, we conclude that $\liminf_{t \to -\infty} R_{\text{\rm max}}(t) > 0$ (see Proposition \ref{lower.bounds.for.d(t).and.R_max(t)}). In Section \ref{uniqueness.rot.sym}, we complete the proof of Theorem \ref{main.thm.a}. The idea is to consider a quantity which is constant on the Bryant soliton, and then show that this quantity must be constant on any noncompact ancient $\kappa$-solution with rotational symmetry.

In the second part of this paper, we show that every noncompact ancient $\kappa$-solution must be rotationally symmetric. In Section \ref{lie.derivatives}, we derive a crucial evolution equation for the Lie derivative of the metric along a vector field. In Sections \ref{model.problem} -- \ref{NIT}, we establish a Neck Improvement Theorem for the Ricci flow, which tells us that a neck-like region becomes more symmetric under the evolution. The proof of the Neck Improvement Theorem is based on the vector field method developed in \cite{Brendle}, and requires a careful analysis of the parabolic Lichnerowicz equation on the cylinder. Finally, in Section \ref{rot.sym}, we complete the proof of Theorem \ref{main.thm.b}. The idea is as follows. Since our solution is of Type II, we can find a sequence of points $(\hat{p}_k,\hat{t}_k)$ in spacetime such that, if we rescale the flow around $(\hat{p}_k,\hat{t}_k)$, then the rescaled flows converge to a steady gradient Ricci soliton as $k \to \infty$. By \cite{Brendle}, this limiting soliton must be the Bryant soliton. In particular, we can find a sequence $\hat{\varepsilon}_k \to 0$ such that the flow is $\hat{\varepsilon}_k$-symmetric at time $\hat{t}_k$ (see Definition \ref{symmetry.of.cap} for a precise definition). We now move forward in time, starting from time $\hat{t}_k$. As long as the solution is nearly rotationally symmetric, it will remain close to the Bryant soliton by Theorem \ref{main.thm.a}. On the other hand, as long as the cap is close to the Bryant soliton, we are able to show that the symmetry improves under the evolution (see Proposition \ref{cap.improvement.theorem}). Using a continuity argument, we are able to show that there exists a sequence $\varepsilon_k \geq 2\hat{\varepsilon}_k$ such that $\varepsilon_k \to 0$ and the flow is $\varepsilon_k$-symmetric at time $t$ for all $t \in [\hat{t}_k,0]$. Passing to the limit as $k \to \infty$, it follows that $(M,g(t))$ is rotationally symmetric for all $t$. 

\begin{remark} 
The proof of Theorem \ref{main.thm.b} can be adapted to the compact case. This will imply that every ancient $\kappa$-solution on $S^3$ must be rotationally symmetric. 
\end{remark}

\textbf{Acknowledgements.} I am grateful to Robert Bryant for sharing with me his insights on singular Ricci solitons (cf. Theorem \ref{singular.soliton}), and to Kyeongsu Choi for pointing out to me his variant of the Anderson-Chow estimate (cf. Theorem \ref{anderson.chow.estimate}). I would like to thank Keaton Naff for comments on an earlier version of this paper. I am grateful to T\"ubingen University, where part of this work was carried out. This project was supported by the National Science Foundation under grant DMS-1806190 and by the Simons Foundation.

\part{Proof of Theorem \ref{main.thm.a}}

\section{A barrier construction} 

\label{barrier.construction}

In this section, we study the Ricci flow in the rotationally symmetric setting. In this case, the Ricci flow reduces to a parabolic equation for a single scalar function (see \cite{Angenent-Isenberg-Knopf}). We first construct a family of functions $\psi_a$ which will serve as barriers. A key ingredient in our construction is the following result due to Robert Bryant \cite{Bryant} (see also \cite{Alexakis-Chen-Fournodavlos}, Proposition 2.1):

\begin{theorem}[R.~Bryant \cite{Bryant}, Section 3.4]
\label{singular.soliton}
There exists a steady gradient Ricci soliton which is rotationally symmetric, singular at the tip, and asymptotic to the Bryant soliton near infinity. This soliton can be written in the form $\varphi(r)^{-1} \, dr \otimes dr + r^2 \, g_{S^2}$, where $\varphi(r)$ is a positive function defined on the interval $(0,\infty)$ satisfying 
\[\varphi(r) \varphi''(r) - \frac{1}{2} \, \varphi'(r)^2 + r^{-2} \, (1-\varphi(r)) \, (r\varphi'(r)+2\varphi(r)) = 0.\] 
The function $\varphi(r)$ satisfies $\varphi(r) \to \infty$ as $r \to 0$. Like the Bryant soliton, $\varphi(r)$ satisfies an asymptotic expansion of the form $\varphi(r) = r^{-2} + 2r^{-4} + O(r^{-6})$ as $r \to \infty$. 
\end{theorem}

\textbf{Proof.} 
We sketch how Theorem \ref{singular.soliton} follows from Robert Bryant's results. In equation (3.26) in \cite{Bryant}, Robert Bryant considers the ODE 
\[\frac{du}{ds} = \frac{u(1-u^2)s^2}{(u+s)(2-s^2)}.\] 
It is shown in \cite{Bryant} that this ODE admits a solution $u(s)$ which is defined for $s \in (-\sqrt{2},0)$, takes values in the interval $(0,1)$, and satisfies $u(s) \to 1$ as $s \searrow -\sqrt{2}$ and $u(s) \to 0$ as $s \nearrow 0$. Moreover, this solution satisfies $u(s)+s < 0$ for all $s \in (-\sqrt{2},0)$. Given a solution $u(s)$ of this ODE, the metric 
\[g = \frac{1-u^2}{(u+s)^2(2-s^2)^2} \, ds \otimes ds + \frac{1-u^2}{u^2(2-s^2)} \, g_{S^2}\] 
will be a steady gradient Ricci soliton (cf. equation (3.28) in \cite{Bryant}). Using the differential equation for $u$, we compute 
\[\frac{1}{2} \, \frac{d}{ds} \Big ( \frac{1-u^2}{u^2(2-s^2)} \Big ) = \frac{(1-u^2)s}{u(u+s)(2-s^2)^2} > 0\] 
for all $s \in (-\sqrt{2},0)$. Consequently, the function $s \mapsto \frac{1-u^2}{u^2(2-s^2)}$ is strictly monotone increasing. Moreover, $\frac{1-u^2}{u^2(2-s^2)} \to 0$ as $s \searrow -\sqrt{2}$, and $\frac{1-u^2}{u^2(2-s^2)} \to \infty$ as $s \nearrow 0$. Hence, the metric $g$ can be rewritten as 
\[g = \varphi(r)^{-1} \, dr \otimes dr + r^2 \, g_{S^2},\] 
where $\varphi$ is defined by $\varphi \big ( \sqrt{\frac{1-u^2}{u^2(2-s^2)}} \big ) = \frac{s^2}{2-s^2}$. The function $\varphi(r)$ is defined for all $r \in (0,\infty)$, and satisfies the ODE 
\[\varphi(r) \varphi''(r) - \frac{1}{2} \, \varphi'(r)^2 + r^{-2} \, (1-\varphi(r)) \, (r\varphi'(r)+2\varphi(r)) = 0.\] 
Moreover, $\varphi(r) \to \infty$ as $r \to 0$. Finally, after replacing $\varphi(r)$ by $\varphi(cr)$ for a suitable constant $c>0$, the function $\varphi(r)$ will have the desired asymptotic expansion as $r \to \infty$. From this, Theorem \ref{singular.soliton} follows. \\


\begin{remark} 
Robert Bryant proved that there is a one-parameter family of singular steady gradient Ricci solitons, which satisfy $\varphi(r) \sim r^{-2(\sqrt{2}-1)}$ as $r \to 0$. However, for the purposes of this paper, one example is sufficient. 
\end{remark}

In the following, we fix a function $\varphi$ as in Theorem \ref{singular.soliton}. Moreover, we fix a positive number $r_*$ such that $\varphi(r_*)=2$. 

Let us choose a smooth function $\zeta$ such that 
\[\frac{d}{ds} \big [ (s^{-2}-1)^{-1} \zeta(s) \big ] = (s^{-2}-1)^{-2} \, \Big [ 2s^{-3}-5s^{-6} - \frac{1}{2} \, s^{27} \Big ].\] 
Note that 
\[(s^{-2}-1)^{-2} \, \Big [ 2s^{-3}-5s^{-6} - \frac{1}{2} \, s^{27} \Big ] = -5s^{-2} + O(1)\] 
as $s \to 0$, and 
\[(s^{-2}-1)^{-2} \, \Big [ 2s^{-3}-5s^{-6} - \frac{1}{2} \, s^{27} \Big ] = -\frac{7}{8} \, (1-s)^{-2} + O(1)\] 
as $s \to 1$. The first statement gives $\zeta(s) = 5s^{-3} + O(s^{-2})$ as $s \to 0$. The second statement implies that $\zeta(s)$ is indeed smooth at $s=1$, and $\zeta(1)=-\frac{7}{4}$. By continuity, we can find a small constant $\theta \in (0,\frac{1}{100})$ such that $2s^{-4}+\zeta(s) \geq \frac{1}{8}$ for all $s \in [1-\theta,1+\theta]$. \\


\begin{lemma}
\label{positivity}
We can find a large constant $N$ with the following property. If $a$ is sufficiently large, then 
\[\varphi(as) - a^{-2} + a^{-4} \zeta(s) \geq a^{-2} (s^{-2}-1) + \frac{1}{16} \, a^{-4}\] 
for all $s \in [1-\theta,1+\frac{1}{100} \, a^{-2}]$, and 
\[\varphi(as) - a^{-2} + a^{-4} \zeta(s) \geq \frac{1}{32} \, a^{-4}\] 
for all $s \in [Na^{-1},1+\frac{1}{100} \, a^{-2}]$. 
\end{lemma} 

\textbf{Proof.}
Since $2s^{-4}+\zeta(s) \geq \frac{1}{8}$ for all $s \in [1-\theta,1+\theta]$, we obtain 
\begin{align*} 
\varphi(as) - a^{-2} + a^{-4} \zeta(s) 
&= a^{-2} (s^{-2}-1) + a^{-4}(2s^{-4}+\zeta(s)) + O(a^{-6}) \\ 
&\geq a^{-2}(s^{-2}-1) + \frac{1}{8} \, a^{-4} + O(a^{-6}) 
\end{align*}
for all $s \in [1-\theta,1+\theta]$. Consequently, if $a$ is sufficiently large, then 
\[\varphi(as) - a^{-2} + a^{-4} \zeta(s) \geq a^{-2} (s^{-2}-1) + \frac{1}{16} \, a^{-4}\] 
for all $s \in [1-\theta,1+\theta]$. This proves the first statement. 

In particular, if $a$ is sufficiently large, then 
\[\varphi(as) - a^{-2} + a^{-4} \zeta(s) \geq \frac{1}{32} \, a^{-4}\] 
for all $s \in [1-\theta,1+\frac{1}{100} \, a^{-2}]$. We next observe that 
\[\varphi(as) - a^{-2} + a^{-4} \zeta(s) = a^{-2} s^{-2} - a^{-2} + O(a^{-4} s^{-4})\] 
for all $s \in [r_* a^{-1},1-\theta]$. Hence, if we choose $N$ sufficiently large (depending on $\theta$), then 
\[\varphi(as) - a^{-2} + a^{-4} \zeta(s) \geq (1-\theta) a^{-2}s^{-2} - a^{-2} \geq [(1-\theta)^{-1}-1] \, a^{-2}\] 
for all $s \in [Na^{-1},1-\theta]$. In particular, if $a$ is sufficiently large (depending on $\theta$), then 
\[\varphi(as) - a^{-2} + a^{-4} \zeta(s) \geq \frac{1}{32} \, a^{-4}\] 
for all $s \in [Na^{-1},1-\theta]$. Putting these facts together, the second statement follows. The proof Lemma \ref{positivity} is now complete. \\

\begin{lemma}
\label{barrier.1}
We can find a large constant $N$ with the following property. Suppose that $a$ is sufficiently large, and let 
\[\psi_a(s) := \varphi(as) - a^{-2} + a^{-4} \zeta(s)\] 
for $s \in [Na^{-1},1+\frac{1}{100} \, a^{-2}]$. Then 
\[\psi_a(s) \psi_a''(s) - \frac{1}{2} \, \psi_a'(s)^2 + s^{-2} \, (1-\psi_a(s)) \, (s \psi_a'(s) + 2 \psi_a(s)) - s \psi_a'(s) < 0\] 
for $s \in [Na^{-1},1+\frac{1}{100} \, a^{-2}]$.
\end{lemma} 

\textbf{Proof.} 
The function $\zeta$ satisfies 
\begin{align*} 
s^{-2} (s \zeta'(s) + 2 \zeta(s)) - s \zeta'(s) 
&= s \, (s^{-2}-1)^2 \, \frac{d}{ds} \big [ (s^{-2}-1)^{-1} \zeta(s) \big ] \\ 
&= 2s^{-2} - 5s^{-5} - \frac{1}{2} \, s^{28}. 
\end{align*}
Since $-r\varphi'(r) - 2r^{-2} = 8r^{-4} + O(r^{-6})$ as $r \to \infty$, it follows that 
\begin{align*} 
&s^{-2} (s \psi_a'(s) + 2 \psi_a(s)) - s \psi_a'(s) \\ 
&= s^{-2} (as \varphi'(as) + 2 \varphi(as)) - as \varphi'(as) - 2a^{-2}s^{-2} \\ 
&+ 2a^{-4}s^{-2} - 5a^{-4}s^{-5} - \frac{1}{2} \, a^{-4} s^{28} \\ 
&= s^{-2} (as \varphi'(as) + 2 \varphi(as)) + 8a^{-4}s^{-4} \\ 
&+ 2a^{-4}s^{-2} - 5a^{-4}s^{-5} - \frac{1}{2} \, a^{-4} s^{28} + O(a^{-6}s^{-6})
\end{align*} 
for $s \in [r_*a^{-1},1+\frac{1}{100} \, a^{-2}]$. Moreover, using the identity $\varphi''(r) - r^{-2} \, (r\varphi'(r) + 4\varphi(r)) = 4r^{-4} + O(r^{-6})$ as $r \to \infty$, we obtain 
\begin{align*} 
&\psi_a(s) \psi_a''(s) - \frac{1}{2} \, \psi_a'(s)^2 - s^{-2} \, \psi_a(s) \, (s \psi_a'(s) + 2 \psi_a(s)) \\ 
&= a^2 \varphi(as) \varphi''(as) - \frac{1}{2} \, a^2 \varphi'(as)^2 - s^{-2} \, \varphi(as) \, (as \varphi'(s) + 2 \varphi(as)) \\ 
&- 2a^{-4} s^{-2} - \Big [ \varphi''(as) - a^{-2} s^{-2} \, (as\varphi'(as) + 4\varphi(as)) \Big ] \\ 
&+ a^{-4} \, \Big [ \varphi(as) \, \zeta''(s) + a^2 \, \varphi''(as) \, \zeta(s) - a \, \varphi'(as) \, \zeta'(s) \\ 
&\hspace{15mm} - s^{-1} \, \varphi(as) \, \zeta'(s) - as^{-1} \, \varphi'(as) \, \zeta(s) - 4s^{-2} \, \varphi(as) \, \zeta(s) \Big ] \\  
&- a^{-6} \, \Big [ \zeta''(s) - s^{-2} \, (s\zeta'(s) + 4\zeta(s)) \Big ] \\ 
&+ a^{-8} \, \Big [ \zeta(s) \zeta''(s) - \frac{1}{2} \, \zeta'(s)^2 - s^{-2} \, \zeta(s) \, (s \zeta'(s) + 2 \zeta(s)) \Big ] \\ 
&= a^2 \varphi(as) \varphi''(as) - \frac{1}{2} \, a^2 \varphi'(as)^2 - s^{-2} \, \varphi(as) \, (as \varphi'(s) + 2 \varphi(as)) \\ 
&- 2a^{-4} s^{-2} - 4a^{-4} s^{-4} + O(a^{-6} s^{-7})
\end{align*}
for $s \in [r_*a^{-1},1+\frac{1}{100} \, a^{-2}]$. Adding both identities yields 
\begin{align*} 
&\psi_a(s) \psi_a''(s) - \frac{1}{2} \, \psi_a'(s)^2 + s^{-2} \, (1-\psi_a(s)) \, (s \psi_a'(s) + 2 \psi_a(s)) - s \psi_a'(s) \\ 
&= a^2 \, \Big [ \varphi(as) \varphi''(as) - \frac{1}{2} \, \varphi'(as)^2 + (as)^{-2} \, (1-\varphi(as)) \, (as \varphi'(s) + 2 \varphi(as)) \Big ] \\ 
&+ 4a^{-4}s^{-4}-5a^{-4}s^{-5} - \frac{1}{2} \, a^{-4}s^{28} + O(a^{-6} s^{-7}) 
\end{align*} 
for $s \in [r_*a^{-1},1+\frac{1}{100} \, a^{-2}]$. Using the ODE for $\varphi(r)$, we conclude that 
\begin{align*} 
&\psi_a(s) \psi_a''(s) - \frac{1}{2} \, \psi_a'(s)^2 + s^{-2} \, (1-\psi_a(s)) \, (s \psi_a'(s) + 2 \psi_a(s)) - s \psi_a'(s) \\ 
&= 4a^{-4}s^{-4}-5a^{-4}s^{-5} - \frac{1}{2} \, a^{-4}s^{28} + O(a^{-6} s^{-7}) 
\end{align*}
for $s \in [r_*a^{-1},1+\frac{1}{100} \, a^{-2}]$. Clearly, the expression on the right hand side is negative if $s \in [Na^{-1},1+\frac{1}{100} \, a^{-2}]$ and $N$ is sufficiently large. This completes the proof of Lemma \ref{barrier.1}. \\

From now on, we will fix a large number $N$ so that the conclusions of Lemma \ref{positivity} and Lemma \ref{barrier.1} hold. For $a$ sufficiently large, we can find a smooth function $\beta_a(r)$ such that $\beta_a(N)=a^{-3} \zeta(Na^{-1})-a^{-1}$, $\beta_a'(N)=a^{-4} \zeta'(Na^{-1})$, and 
\begin{align*} 
&\varphi(r) \beta_a''(r) + \varphi''(r) \beta_a(r) - \varphi'(r) \beta_a'(r) \\ 
&+ r^{-2} \, (1-\varphi(r)) \, (r\beta_a'(r)+2\beta_a(r)) \\ 
&- r^{-2} \, \beta_a(r) \, (r\varphi'(r)+2\varphi(r)) \\ 
&= -1 
\end{align*} 
for $r \in [r_*,N]$. Note that $\beta_a(N)$ and $\beta_a'(N)$ are uniformly bounded independent of $a$. Consequently, the function $\beta_a$ and all its derivatives are uniformly bounded on the interval $[r_*,N]$, and the bounds are independent of $a$. \\

\begin{lemma}
\label{barrier.2}
Suppose that $a$ is sufficiently large, and let 
\[\psi_a(s) := \varphi(as) + a^{-1} \beta_a(as)\] 
for $s \in [r_* a^{-1},Na^{-1}]$. Then  
\[\psi_a(s) \psi_a''(s) - \frac{1}{2} \, \psi_a'(s)^2 + s^{-2} \, (1-\psi_a(s)) \, (s \psi_a'(s) + 2 \psi_a(s)) - s \psi_a'(s) < 0\] 
for all $s \in [r_* a^{-1},Na^{-1}]$.
\end{lemma}

\textbf{Proof.} 
Using the ODEs for $\varphi(r)$ and $\beta_a(r)$, we obtain 
\begin{align*} 
&\psi_a(s) \psi_a''(s) - \frac{1}{2} \, \psi_a'(s)^2 + s^{-2} \, (1-\psi_a(s)) \, (s \psi_a'(s) + 2 \psi_a(s)) \\ 
&= a^2 \Big [ \varphi(as) \varphi''(as) - \frac{1}{2} \, \varphi'(as)^2 + (as)^{-2} \, (1-\varphi(as)) \, (as \varphi'(as) + 2 \varphi(as)) \Big ] \\ 
&+ a \Big [ \varphi(as) \beta_a''(as) + \varphi''(as) \beta_a(as) - \varphi'(as) \beta_a'(as) \\ 
&\hspace{10mm} + (as)^{-2} \, (1-\varphi(as)) \, (as \beta_a'(as)+2\beta_a(as)) \\ 
&\hspace{10mm} - (as)^{-2} \, \beta_a(as) \, (as\varphi'(as)+2\varphi(as)) \Big ] \\ 
&+ \Big [ \beta_a(as) \beta_a''(as) - \frac{1}{2} \, \beta_a'(as)^2 - (as)^{-2} \, \beta_a(as) \, (as \beta_a'(as) + 2\beta_a(as)) \Big ] \\ 
&\leq -a + C  
\end{align*} 
for all $s \in [r_* a^{-1},Na^{-1}]$. On the other hand, 
\[s \psi_a'(s) = as \varphi'(as) + s \beta_a'(as) \geq -C\]
for all $s \in [r_* a^{-1},Na^{-1}]$. Hence, if $a$ is sufficiently large, then 
\[\psi_a(s) \psi_a''(s) - \frac{1}{2} \, \psi_a'(s)^2 + s^{-2} \, (1-\psi_a(s)) \, (s \psi_a'(s) + 2 \psi_a(s)) < s \psi_a'(s)\] 
for all $s \in [r_* a^{-1},Na^{-1}]$. This completes the proof of Lemma \ref{barrier.2}. \\

After these preparations, we now give the definition of our barriers:

\begin{definition}
\label{definition.of.psi_a}
Suppose that $a$ is sufficiently large. We define a function $\psi_a: [r_* a^{-1},1+\frac{1}{100} \, a^{-2}] \to \mathbb{R}$ by 
\[\psi_a(s) := \begin{cases} \varphi(as) - a^{-2} + a^{-4} \zeta(s) & \text{\rm for $s \in [Na^{-1},1+\frac{1}{100} \, a^{-2}]$} \\ \varphi(as) + a^{-1} \beta_a(as) & \text{\rm for $s \in [r_* a^{-1},Na^{-1}]$.} \end{cases}\] 
\end{definition}

Using Lemma \ref{positivity}, Lemma \ref{barrier.1}, and Lemma \ref{barrier.2}, we can draw the following conclusion:

\begin{proposition}
\label{properties.of.psi_a} 
Suppose that $a$ is sufficiently large. Then $\psi_a$ is continuously differentiable, and 
\[\psi_a(s) \psi_a''(s) - \frac{1}{2} \, \psi_a'(s)^2 + s^{-2} \, (1-\psi_a(s)) \, (s \psi_a'(s) + 2 \psi_a(s)) - s \psi_a'(s) < 0\] 
for all $s \in [r_* a^{-1},1+\frac{1}{100} \, a^{-2}]$. Moreover, we have $\psi_a(s) \geq \frac{1}{32} \, a^{-4}$ for all $s \in [r_* a^{-1},1+\frac{1}{100} \, a^{-2}]$, and $\psi_a(s) \geq a^{-2} (s^{-2}-1) + \frac{1}{16} \, a^{-4}$ for all $s \in [1-\theta,1+\frac{1}{100} \, a^{-2}]$.
\end{proposition} 

\textbf{Proof.} 
Recall that $\beta_a(N)=a^{-3} \zeta(Na^{-1})-a^{-1}$, $\beta_a'(N)=a^{-4} \zeta'(Na^{-1})$. This implies that $\psi_a$ is continuously differentiable at the point $s=Na^{-1}$. This proves the first statement. The second statement follows from Lemma \ref{barrier.1} and Lemma \ref{barrier.2}. Finally, the third and fourth statement follow directly from Lemma \ref{positivity}. \\

\begin{corollary}
\label{parabolic.barrier}
The function $\Psi_a(r,t) := \psi_a \big ( \frac{r}{\sqrt{-2t}} \big )$ satisfies 
\[\Psi_{a,t} > \Psi_a \Psi_{a,rr} - \frac{1}{2} \, \Psi_{a,r}^2 + r^{-2} \, (1-\Psi_a) \, (r \Psi_{a,r}+2\Psi_a)\] 
for $r \in [r_* a^{-1} \sqrt{-2t},(1+\frac{1}{100} \, a^{-2}) \sqrt{-2t}]$.
\end{corollary}

\textbf{Proof.} 
This follows immediately from Proposition \ref{properties.of.psi_a}. \\

In the remainder of this section, we will set up a barrier argument based on the functions $\psi_a$. We will assume throughout that $(M,g(t))$ is a three-dimensional ancient $\kappa$-solution which is noncompact, has positive sectional curvature, and is rotationally symmetric. After a reparametrization, the metric can be written in the form $\tilde{g}(t) = u(r,t)^{-1} \, dr \otimes dr + r^2 \, g_{S^2}$. For each $t$, the function $r \mapsto u(r,t)$ is defined on an interval $[0,r_{\text{\rm max}}(t))$, where $r_{\text{\rm max}}(t)$ may be finite or infinite. 

The Ricci and scalar curvature of $\tilde{g}$ are given by 
\[\text{\rm Ric}_{\tilde{g}} = -\frac{1}{r} \, u^{-1} u_r \, dr \otimes dr + \Big ( 1-u - \frac{1}{2} \, ru_r \Big ) \, g_{S^2}\] 
and 
\[R_{\tilde{g}} = \frac{2}{r^2} \, (1-u-ru_r)\] 
(cf. \cite{Angenent-Knopf}, p.~497). Since the original metrics $g(t)$ evolve by the Ricci flow, the reparametrized metrics $\tilde{g}(t)$ satisfy an evolution equation of the form 
\[\frac{\partial}{\partial t} \tilde{g} = -2 \, \text{\rm Ric}_{\tilde{g}} + \mathscr{L}_V(\tilde{g})\] 
where $V$ is a radial vector field of the form $V = v(r,t) \, \frac{\partial}{\partial r}$ which may depend on time. 

Clearly, $\frac{\partial}{\partial t} \tilde{g} = -u^{-2} \, u_t \, dr \otimes dr$. Moreover, $\mathscr{L}_V(r) = v$ and $\mathscr{L}_V(dr) = v_r \, dr$. This gives 
\[\mathscr{L}_V(\tilde{g}) = (-u^{-2} u_r v + 2 u^{-1} v_r) \, dr \otimes dr + 2r v \, g_{S^2},\] 
hence 
\begin{align*} 
\text{\rm Ric}_{\tilde{g}} - \frac{1}{2} \, \mathscr{L}_V(\tilde{g}) 
&= \Big ( -\frac{1}{r} \, u^{-1} u_r + \frac{1}{2} \, u^{-2} u_r v - u^{-1} v_r \Big ) \, dr \otimes dr \\ 
&+ \Big ( 1 - u - \frac{1}{2} \, ru_r - rv \Big ) \, g_{S^2}. 
\end{align*} 
Putting these facts together, we conclude that 
\[v = \frac{1}{r} \, (1-u - \frac{1}{2} \, ru_r)\]
and 
\begin{align*} 
u_t 
&= 2 \, \Big ( -\frac{1}{r} \, u u_r + \frac{1}{2} \, u_r v - u v_r \Big ) \\ 
&= uu_{rr} - \frac{1}{2} \, u_r^2 + r^{-2} \, (1-u) \, (ru_r+2u). 
\end{align*}
The function $u$ has a natural geometric interpretation. Namely, we can view the radius $r$ as a scalar function on $M$. Then $u = |dr|_{\tilde{g}(t)}^2$. In particular, $u$ is very small on a neck. 

\begin{lemma}
\label{behavior.at.tip}
We have $u(r,t) \leq 1$, $u_r(r,t) \leq 0$, and $v(r,t) \geq 0$ at each point in space-time. Moreover, $1-u(r,t)=O(r^2)$ and $v(r,t)=O(r)$ as $r \to 0$. 
\end{lemma}

\textbf{Proof.} 
Since the metric is smooth at the tip, we obtain $1-u(r,t)=O(r^2)$ and $v(r,t)=O(r)$ as $r \to 0$. Since $(M,g(t))$ has positive Ricci curvature, we have $u_r(r,t) = -r \, {\text{\rm Ric}_r}^r \leq 0$ at each point in space-time. Integrating over $r$, we obtain $u(r,t) \leq 1$ at each point in space-time. Finally, $v(r,t) = \frac{r}{2} \, (R-{\text{\rm Ric}_r}^r) \geq 0$ at each point in space-time. \\

\begin{lemma} 
\label{diameter.sphere.of.symmetry}
If a sphere of symmetry in $(M,g(t))$ has radius $r$, then its diameter in $(M,g(t))$ is at least $2r$.
\end{lemma}

\textbf{Proof.} 
By Lemma \ref{behavior.at.tip}, we have $u(r,t) \leq 1$. Consequently, the metric satisfies $u(r,t)^{-1} \, dr \otimes dr + r^2 \, g_{S^2} \geq dr \otimes dr + r^2 \, g_{S^2}$. This allows us to compare the distance function in $(M,g(t))$ to the distance function in Euclidean space. In particular, if we consider two antipodal points on a sphere of radius $r$, then their geodesic distance in $(M,g(t))$ is at least $2r$. \\

\begin{lemma}
\label{u.converges.to.0}
Given any $\delta>0$, we have $\liminf_{t \to -\infty} \sup_{r \geq \delta \sqrt{-t}} u(r,t) = 0$. 
\end{lemma}

\textbf{Proof.} 
Let $\varepsilon>0$ be given. For each $t$, we denote by $R_{\text{\rm max}}(t)$ the supremum of the scalar curvature of $(M,g(t))$. By work of Perelman \cite{Perelman1}, the set of all points in $(M,g(t))$ which do not lie on an $\varepsilon$-neck has diameter less than $C(\varepsilon) \, R_{\text{\rm max}}(t)^{-\frac{1}{2}}$ (see Theorem \ref{canonical.neighborhood.theorem} and Corollary \ref{bound.for.R_max}). Hence, if $r > C(\varepsilon) \, R_{\text{\rm max}}(t)^{-\frac{1}{2}}$ at some point in spacetime, then that point lies on an $\varepsilon$-neck, and we have $u \leq 2\varepsilon$. Thus, $\sup_{r > C(\varepsilon) \, R_{\text{\rm max}}(t)^{-\frac{1}{2}}} u(r,t) \leq 2\varepsilon$ for each $t$. On the other hand, we know that our ancient solution is of Type II (cf. \cite{Zhang}), so that $\limsup_{t \to \infty} (-t) \, R_{\text{\rm max}}(t) = \infty$. Putting these facts together, we conclude that $\liminf_{t \to -\infty} \sup_{r \geq \delta \sqrt{-t}} u(r,t) \leq 2\varepsilon$ for each $\delta>0$. Since $\varepsilon>0$ is arbitrary, it follows that $\liminf_{t \to -\infty} \sup_{r \geq \delta \sqrt{-t}} u(r,t) = 0$ for each $\delta>0$. This completes the proof of Lemma \ref{u.converges.to.0}. \\

\begin{proposition}
\label{maximum.principle}
There exists a large number $K$ with the following property. Suppose that $a \geq K$ and $\bar{t} \leq 0$. Moreover, suppose that $\bar{r}(t) \in [0,r_{\text{\rm max}}(t))$ is a function satisfying $\big | \frac{\bar{r}(t)}{\sqrt{-2t}}-1 \big | \leq \frac{1}{100} \, a^{-2}$ and $u(\bar{r}(t),t) < \frac{1}{32} \, a^{-4}$ for all $t \leq \bar{t}$. Then $u(r,t) \leq \psi_a(\frac{r}{\sqrt{-2t}})$ whenever $t \leq \bar{t}$ and $r_* a^{-1} \sqrt{-2t} \leq r \leq \bar{r}(t)$. In particular, $u(r,t) \leq C a^{-2}$ whenever $t \leq \bar{t}$ and $\frac{1}{2} \, \sqrt{-2t} \leq r \leq \bar{r}(t)$.
\end{proposition}

\textbf{Proof.} 
By Proposition \ref{properties.of.psi_a}, we can find a large constant $K$ such that $\psi_a(s) \geq \frac{1}{32} \, a^{-4}$ for all $s \in [r_* a^{-1},1+\frac{1}{100} \, a^{-2}]$ and all $a \geq K$. Moreover, we can arrange that $1+a^{-1} \beta_a(r_*) > 0$ for all $a \geq K$.

We claim that $K$ has the desired property. To see this, we fix an arbitrary number $a \geq K$. Moreover, suppose that $\bar{r}(t) \in [0,r_{\text{\rm max}}(t))$ is a function satisfying $\big | \frac{\bar{r}(t)}{\sqrt{-2t}}-1 \big | \leq \frac{1}{100} \, a^{-2}$ and $u(\bar{r}(t),t) < \frac{1}{32} \, a^{-4}$ for all $t \leq \bar{t}$. Then 
\[\psi_a \Big ( \frac{\bar{r}(t)}{\sqrt{-2t}} \Big ) - u(\bar{r}(t),t) \geq \frac{1}{32} \, a^{-4} - u(\bar{r}(t),t) > 0\] 
for all $t \leq \bar{t}$. Moreover, since $\varphi(r_*)=2$ and $u \leq 1$, we have 
\begin{align*} 
&\psi_a(r_* a^{-1}) - u(r_* a^{-1} \sqrt{-2t},t) \\ 
&= 2+a^{-1} \beta_a(r_*) - u(r_* a^{-1} \sqrt{-2t},t) \geq 1+a^{-1} \beta_a(r_*) > 0 
\end{align*}
for all $t \leq \bar{t}$. On the other hand, Lemma \ref{u.converges.to.0} implies 
\[\limsup_{t \to -\infty} \inf_{r_* a^{-1} \sqrt{-2t} \leq r \leq \bar{r}(t)} \Big [ \psi_a \Big ( \frac{r}{\sqrt{-2t}} \Big ) - u(r,t) \Big ] > 0.\] 
Using Corollary \ref{parabolic.barrier} and the maximum principle, we obtain 
\[\psi_a \Big ( \frac{r}{\sqrt{-2t}} \Big ) - u(r,t) \geq 0\] 
whenever $t \leq \bar{t}$ and $r_* a^{-1} \sqrt{-2t} \leq r \leq \bar{r}(t)$. This completes the proof of Proposition \ref{maximum.principle}. \\

\begin{proposition}
\label{distance.of.reference.point.to.tip}
Suppose that there exists a function $\bar{r}(t) \in [0,r_{\text{\rm max}}(t))$ such that $\bar{r}(t) = \sqrt{-2t} + O(1)$ and $u(\bar{r}(t),t) \leq O(\frac{1}{(-t)})$ as $t \to -\infty$. Then we can find a large constant $K \geq 100$ with the property that 
\[\frac{r}{\sqrt{-2t+Ka^2}} \leq 1+\frac{1}{100} \, a^{-2}\] 
and 
\[u(r,t) \leq \psi_a \Big ( \frac{r}{\sqrt{-2t+Ka^2}} \Big )\] 
whenever $a \geq K$, $t \leq -K^2a^2$, and $r_* a^{-1} \sqrt{-2t+Ka^2} \leq r \leq \bar{r}(t)$. Note that $K$ is independent of $a$ and $t$. Moreover, 
\[\liminf_{t \to -\infty} (-t)^{-1} \int_0^{\bar{r}(t)} u(r,t)^{-\frac{1}{2}} \, dr > 0.\] 
\end{proposition} 

\textbf{Proof.}
We choose $K \geq 100$ sufficiently large so that the following holds: 
\begin{itemize}
\item $\frac{\bar{r}(t)}{\sqrt{-(2+K^{-1})t}} \geq 1-\theta$ for all $t \leq -K^4$. 
\item $\bar{r}(t)^2+2t \leq \frac{\sqrt{K}}{10} \, \bar{r}(t)$ for all $t \leq -K^4$.
\item $u(\bar{r}(t),t) \leq \frac{K}{2\bar{r}(t)^2}$ for all $t \leq -K^4$. 
\item $\psi_a(s) \geq a^{-2} (s^{-2}-1) + \frac{1}{16} \, a^{-4} > 0$ for all $s \in [1-\theta,1+\frac{1}{100} \, a^{-2}]$ and all $a \geq K$.
\item $1+a^{-1} \beta_a(r_*) > 0$ for all $a \geq K$. 
\end{itemize} 
We claim that $K$ has the desired property. To prove this, we fix an arbitrary number $a \geq K$. Clearly, $\frac{\bar{r}(t)}{\sqrt{-2t+Ka^2}} \geq \frac{\bar{r}(t)}{\sqrt{-(2+K^{-1})t}} \geq 1-\theta$ for all $t \leq -K^2a^2$. Moreover, using the inequality $\bar{r}(t)^2+2t \leq \frac{\sqrt{K}}{10} \, \bar{r}(t)$, we obtain 
\begin{align*} 
\frac{-2t+Ka^2}{\bar{r}(t)^2} - 1 + \frac{1}{100} \, a^{-2} 
&= \frac{Ka^2}{\bar{r}(t)^2} - \frac{\bar{r}(t)^2+2t}{\bar{r}(t)^2} + \frac{1}{100} \, a^{-2} \\ 
&\geq \frac{Ka^2}{\bar{r}(t)^2} - \frac{\sqrt{K}}{10\bar{r}(t)} + \frac{1}{100} \, a^{-2} \\ 
&= \frac{3Ka^2}{4\bar{r}(t)^2} + \Big ( \frac{\sqrt{K} a}{2\bar{r}(t)} - \frac{1}{10} \, a^{-1} \Big )^2
\end{align*}
for all $t \leq -K^2a^2$. Since the right hand side is positive, it follows that 
\[\frac{\bar{r}(t)}{\sqrt{-2t+Ka^2}} \leq \Big ( 1-\frac{1}{100} \, a^{-2} \Big )^{-\frac{1}{2}} \leq 1+\frac{1}{100} \, a^{-2}\] 
for all $t \leq -K^2a^2$. This proves the first statement. 

We next observe that 
\begin{align*} 
&\psi_a \Big ( \frac{\bar{r}(t)}{\sqrt{-2t+Ka^2}} \Big ) - u(\bar{r}(t),t) \\ 
&\geq a^{-2} \, \Big ( \frac{-2t+Ka^2}{\bar{r}(t)^2} - 1 \Big ) + \frac{1}{100} \, a^{-4} - u(\bar{r}(t),t) \\ 
&\geq \frac{3K}{4\bar{r}(t)^2} - u(\bar{r}(t),t) \\ 
&> 0
\end{align*}
for all $t \leq -K^2a^2$. Moreover, since $\varphi(r_*)=2$ and $u \leq 1$, we have 
\begin{align*} 
&\psi_a(r_* a^{-1}) - u(r_* a^{-1} \sqrt{-2t+Ka^2},t) \\ 
&= 2+a^{-1} \beta_a(r_*) - u(r_* a^{-1} \sqrt{-2t+Ka^2},t) \geq 1+a^{-1} \beta_a(r_*) > 0 
\end{align*} 
for all $t \leq -K^2a^2$. On the other hand, Lemma \ref{u.converges.to.0} implies 
\[\limsup_{t \to -\infty} \inf_{r_* a^{-1} \sqrt{-2t+Ka^2} \leq r \leq \bar{r}(t)} \Big [ \psi_a \Big ( \frac{r}{\sqrt{-2t+Ka^2}} \Big ) - u(r,t) \Big ] > 0.\] 
Using Corollary \ref{parabolic.barrier} and the maximum principle, we obtain
\[\psi_a \Big ( \frac{r}{\sqrt{-2t+Ka^2}} \Big ) - u(r,t) \geq 0\] 
for all $t \leq -K^2a^2$ and $r_* a^{-1} \sqrt{-2t+Ka^2} \leq r \leq \bar{r}(t)$. This proves the second statement. 

To prove the last statement, we recall that $\frac{\bar{r}(t)}{\sqrt{-2t+Ka^2}} \geq 1-\theta$ for all $t \leq -K^2a^2$. Consequently, 
\begin{align*} 
\int_{r_* a^{-1} \sqrt{-2t+Ka^2}}^{\bar{r}(t)} u(r,t)^{-\frac{1}{2}} \, dr 
&\geq \int_{r_* a^{-1} \sqrt{-2t+Ka^2}}^{(1-\theta) \sqrt{-2t+Ka^2}} \psi_a \Big ( \frac{r}{\sqrt{-2t+Ka^2}} \Big )^{-\frac{1}{2}} \, dr \\ 
&= \sqrt{-2t+Ka^2} \int_{r_* a^{-1}}^{1-\theta} \psi_a(s)^{-\frac{1}{2}} \, ds \\ 
&\geq \frac{a}{C} \, \sqrt{-2t+Ka^2}
\end{align*}
for $t \leq -K^2a^2$. To summarize, we have shown that $\int_0^{\bar{r}(t)} u(r,t)^{-\frac{1}{2}} \, dr \geq \frac{a}{C} \, \sqrt{-t}$ whenever $a \geq K$ and $t \leq -K^2 a^2$. Putting $t = -K^2a^2$, we conclude that $\int_0^{\bar{r}(t)} u(r,t)^{-\frac{1}{2}} \, dr \geq \frac{1}{CK} \, (-t)$ for $t \leq -K^4$. This completes the proof of Proposition \ref{distance.of.reference.point.to.tip}. \\

\section{Asymptotics of ancient $\kappa$-solutions with rotational symmetry} 

\label{asymptotics}

We continue to assume that $(M,g(t))$ is a three-dimensional ancient $\kappa$-solution which is noncompact, has positive sectional curvature, and is rotationally symmetric. Let $q \in M$ be a fixed reference point satisfying $\sup_{t \leq 0} (-t) \, R(q,t) \leq 100$; such a point exists by Theorem \ref{consequence.of.kleiner.lott}.

\begin{proposition}
\label{asymptotic.shrinking.soliton}
If we dilate $(M,g(t))$ around $q$ by the factor $(-t)^{-\frac{1}{2}}$, then the rescaled manifolds converge in the Cheeger-Gromov sense to a cylinder of radius $\sqrt{2}$.
\end{proposition}

\textbf{Proof.} 
Recall that $\sup_{t \leq 0} (-t) \, R(q,t) \leq 100$ by our choice of $q$. Let $\ell$ denote the reduced distance from $(q,0)$. Moreover, let us consider an arbitrary sequence of times $t_k \to -\infty$. Then $\ell(q,t_k) \leq \frac{1}{2\sqrt{-t_k}} \int_{t_k}^0 \sqrt{-t} \, R(q,t) \, dt \leq 1000$ if $k$ is sufficiently large. Let us dilate the flow $(M,g(t))$ around $(q,t_k)$ by the factor $(-t_k)^{-\frac{1}{2}}$. By work of Perelman, the rescaled flows converge in the Cheeger-Gromov sense to a shrinking gradient Ricci soliton (see \cite{Perelman1}, Section 11), and this asymptotic soliton must be a cylinder (cf. \cite{Perelman2}, Section 1). This completes the proof. \\

For each $t$, we denote by $\bar{r}(t) \in [0,r_{\text{\rm max}}(t))$ the radius of the sphere of symmetry passing through the point $q$. By Proposition \ref{asymptotic.shrinking.soliton}, $\frac{\bar{r}(t)}{\sqrt{-2t}} \to 1$ as $t \to -\infty$. Since $q$ is fixed, $\bar{r}(t)$ satisfies the following ODE: 
\[\frac{d}{dt} \bar{r}(t) = -v(\bar{r}(t),t) = -\frac{1}{\bar{r}(t)} \, \big ( 1-u(\bar{r}(t),t)-\frac{1}{2} \, \bar{r}(t) \, u_r(\bar{r}(t),t) \big ).\] 
We define a function $F(z,t)$ by 
\[F \bigg ( \int_{\bar{r}(t)}^\rho u(r,t)^{-\frac{1}{2}} \, dr,t \bigg ) = \rho.\] 
In other words, for each time $t$, the function $F(z,t)$ tells us the radius as a function of the signed distance $z$ from the reference point $q$. For each $t$, the function $z \mapsto F(z,t)$ is defined on the interval $[-d(t),\infty)$, where 
\[d(t) = \int_0^{\bar{r}(t)} u(r,t)^{-\frac{1}{2}} \, dr\] 
denotes the distance of the reference point $q$ from the tip. Note that $F(-d(t),t) = 0$.

\begin{proposition}
\label{evolution.of.F}
The function $F$ satisfies 
\begin{align*} 
0 
&= F_t(z,t) - F_{zz}(z,t) + F(z,t)^{-1} \, (1+F_z(z,t)^2) \\ 
&+ 2 \, F_z(z,t) \, \bigg [ -F(0,t)^{-1} \, F_z(0,t) + \int_{F(0,t)}^{F(z,t)} \frac{1}{r^2} \, u(r,t)^{\frac{1}{2}} \, dr \bigg ]. 
\end{align*} 
\end{proposition}

\textbf{Proof.} 
Differentiating the identity 
\[\rho = F \bigg ( \int_{\bar{r}(t)}^\rho u(r,t)^{-\frac{1}{2}} \, dr,t \bigg )\] 
with respect to $\rho$ gives 
\[1 = F_z \bigg ( \int_{\bar{r}(t)}^\rho u(r,t)^{-\frac{1}{2}} \, dr,t \bigg ) \, u(\rho,t)^{-\frac{1}{2}}.\] 
Taking another derivative with respect to $\rho$ gives 
\begin{align*} 
0 &= F_{zz} \bigg ( \int_{\bar{r}(t)}^\rho u(r,t)^{-\frac{1}{2}} \, dr,t \bigg ) \, u(\rho,t)^{-1} \\ 
&- \frac{1}{2} \, F_z \bigg ( \int_{\bar{r}(t)}^\rho u(r,t)^{-\frac{1}{2}} \, dr,t \bigg ) \, u(\rho,t)^{-\frac{3}{2}} \, u_r(\rho,t). 
\end{align*} 
Therefore, 
\[F_z \bigg ( \int_{\bar{r}(t)}^\rho u(r,t)^{-\frac{1}{2}} \, dr,t \bigg ) = u(\rho,t)^{\frac{1}{2}}\] 
and 
\[F_{zz} \bigg ( \int_{\bar{r}(t)}^\rho u(r,t)^{-\frac{1}{2}} \, dr,t \bigg ) = \frac{1}{2} \, u_r(\rho,t).\] 
Using the identity 
\begin{align*} 
\frac{\partial}{\partial t}  (u^{-\frac{1}{2}}) 
&= -\frac{1}{2} \, u^{-\frac{3}{2}} \, \Big ( uu_{rr} - \frac{1}{2} \, u_r^2 + \frac{1}{r^2} \, (1-u) \, (ru_r+2u) \Big ) \\ 
&= \frac{\partial}{\partial r} \Big ( \frac{1}{r} \, u^{-\frac{1}{2}} \, (1+u-\frac{1}{2} ru_r) \Big ) + \frac{2}{r^2} \, u^{\frac{1}{2}},
\end{align*}
we obtain 
\begin{align*} 
\frac{\partial}{\partial t} \bigg ( \int_{\bar{r}(t)}^\rho u(r,t)^{-\frac{1}{2}} \, dr \bigg ) 
&= \frac{1}{\rho} \, u(\rho,t)^{-\frac{1}{2}} \, (1+u(\rho,t)-\frac{1}{2} \, \rho \, u_r(\rho,t)) \\ 
&- \frac{1}{\bar{r}(t)} \, u(\bar{r}(t),t)^{-\frac{1}{2}} \, (1+u(\bar{r}(t),t)-\frac{1}{2} \, \bar{r}(t) \, u_r(\bar{r}(t),t)) \\ 
&- u(\bar{r}(t),t)^{-\frac{1}{2}} \, \frac{d}{dt} \bar{r}(t) + \int_{\bar{r}(t)}^\rho \frac{2}{r^2} \, u(r,t)^{\frac{1}{2}} \, dr \\ 
&= \frac{1}{\rho} \, u(\rho,t)^{-\frac{1}{2}} \, (1+u(\rho,t)-\frac{1}{2} \, \rho \, u_r(\rho,t)) \\ 
&- \frac{2}{\bar{r}(t)} \, u(\bar{r}(t),t)^{\frac{1}{2}} + \int_{\bar{r}(t)}^\rho \frac{2}{r^2} \, u(r,t)^{\frac{1}{2}} \, dr. 
\end{align*} 
Hence, if we differentiate the identity 
\[\rho = F \bigg ( \int_{\bar{r}(t)}^\rho u(r,t)^{-\frac{1}{2}} \, dr,t \bigg )\] 
with respect to $t$, we find  
\begin{align*} 
0 &= F_t \bigg ( \int_{\bar{r}(t)}^\rho u(r,t)^{-\frac{1}{2}} \, dr,t \bigg ) \\ 
&+ F_z \bigg ( \int_{\bar{r}(t)}^\rho u(r,t)^{-\frac{1}{2}} \, dr,t \bigg ) \, \frac{1}{\rho} \, u(\rho,t)^{-\frac{1}{2}} \, (1+u(\rho,t)-\frac{1}{2} \, \rho \, u_r(\rho,t)) \\ 
&+ 2 \, F_z \bigg ( \int_{\bar{r}(t)}^\rho u(r,t)^{-\frac{1}{2}} \, dr,t \bigg ) \, \bigg [ -\frac{1}{\bar{r}(t)} \, u(\bar{r}(t),t)^{\frac{1}{2}} + \int_{\bar{r}(t)}^\rho \frac{1}{r^2} \, u(r,t)^{\frac{1}{2}} \, dr \bigg ]. 
\end{align*} 
Putting these facts together, we conclude that  
\begin{align*}
0 &= F_t \bigg ( \int_{\bar{r}(t)}^\rho u(r,t)^{-\frac{1}{2}} \, dr,t \bigg ) - F_{zz} \bigg ( \int_{\bar{r}(t)}^\rho u(r,t)^{-\frac{1}{2}} \, dr,t \bigg ) \\ 
&+ F \bigg ( \int_{\bar{r}(t)}^\rho u(r,t)^{-\frac{1}{2}} \, dr,t \bigg )^{-1} \, \bigg [ 1 + F_z \bigg ( \int_{\bar{r}(t)}^\rho u(r,t)^{-\frac{1}{2}} \, dr,t \bigg )^2 \bigg ] \\ 
&+ 2 \, F_z \bigg ( \int_{\bar{r}(t)}^\rho u(r,t)^{-\frac{1}{2}} \, dr,t \bigg ) \, \bigg [ -\frac{1}{\bar{r}(t)} \, u(\bar{r}(t),t)^{\frac{1}{2}} + \int_{\bar{r}(t)}^\rho \frac{1}{r^2} \, u(r,t)^{\frac{1}{2}} \, dr \bigg ]. 
\end{align*} 
Using the relations $F(0,t)=\bar{r}(t)$ and $F_z(0,t)=u(\bar{r}(t),t)^{\frac{1}{2}}$, the assertion follows. \\ 

\begin{corollary}
\label{evolution.of.F.2}
The function $F$ satisfies 
\begin{align*} 
&\Big | F_t(z,t) - F_{zz}(z,t) + F(z,t)^{-1} \, (1+F_z(z,t)^2) \Big | \\ 
&\leq 2 \, F(0,t)^{-1} \, F_z(0,t) \, F_z(z,t) \\ 
&+ 2 \, \Big | \frac{1}{F(z,t)} - \frac{1}{F(0,t)} \Big | \, \max \{F_z(z,t),F_z(0,t)\} \, F_z(z,t). 
\end{align*}
\end{corollary}

\textbf{Proof.} 
By Lemma \ref{behavior.at.tip}, the function $u(r,t)$ is monotone decreasing in $r$. Hence, if $r$ lies in between $F(0,t)$ and $F(z,t)$, then $u(r,t)^{\frac{1}{2}} \leq \max \{F_z(z,t),F_z(0,t)\}$. This implies 
\[\bigg | \int_{F(0,t)}^{F(z,t)} \frac{1}{r^2} \, u(r,t)^{\frac{1}{2}} \, dr \bigg | \leq \Big | \frac{1}{F(z,t)} - \frac{1}{F(0,t)} \Big | \, \max \{F_z(z,t),F_z(0,t)\}.\] 
Therefore, the assertion follows from Proposition \ref{evolution.of.F}. \\

\begin{proposition}
\label{pointwise.derivative.bound.for.F}
We have the pointwise estimate 
\[F(z,t)^m \, |\partial_z^{m+1} F(z,t)| \leq C(m) \, (1+F(z,t) \, |F_{zz}(z,t)|)^m\] 
for each $m \geq 0$.
\end{proposition}

\textbf{Proof.} 
We argue by induction on $m$. Lemma \ref{behavior.at.tip} implies $0 \leq F_z \leq 1$ at each point in space-time. Consequently, the assertion holds for $m=0$. Moreover, the assertion clearly holds for $m=1$.

Suppose now that $m \geq 2$, and the assertion holds for all integers less than $m$. Using the standard formula for the scalar curvature of a warped product, we obtain 
\[R = 2 \, F^{-2} \, (1-F_z^2-2 \, F \, F_{zz}).\] 
Differentiating this identity with respect to $z$ gives 
\begin{align*} 
&\partial_z^{m-1} R + 4 \, F^{-1} \, \partial_z^{m+1} F \\ 
&= \sum_{k=0}^{m+1} \sum_{\substack{i_1 \geq 0,\hdots,i_k \geq 0 \\ i_1+\hdots+i_k \leq m-1}} c_{i_1 \hdots i_k} \, F^{i_1+\hdots+i_k-m-1} \, \partial_z^{i_1+1} F \cdots \partial_z^{i_k+1} F. 
\end{align*}
Using the induction hypothesis, we obtain 
\[|\partial_z^{m-1} R + 4 \, F^{-1} \, \partial_z^{m+1} F| \leq C(m) \, F^{-m-1} \, (1+F \, |F_{zz}|)^{m-1}.\] 
On the other hand, Perelman's pointwise curvature derivative estimate (cf. \cite{Perelman1}) implies 
\[|\partial_z^{m-1} R| \leq C(m) \, R^{\frac{m+1}{2}} \leq C(m) \, F^{-m-1} \, (1+F \, |F_{zz}|)^{\frac{m+1}{2}}.\] 
Putting these facts together, we conclude that 
\[|F^{-1} \, \partial_z^{m+1} F| \leq C(m) \, F^{-m-1} \, (1+F \, |F_{zz}|)^m.\] 
This completes the proof of Proposition \ref{pointwise.derivative.bound.for.F}. \\

We now perform a rescaling. For $\tau \leq 0$, we define 
\[G(\xi,\tau) := e^{\frac{\tau}{2}} \, F(e^{-\frac{\tau}{2}} \xi,-e^{-\tau}) - \sqrt{2}.\] 
Since $u(r,t)>0$ and $u_r(r,t) \leq 0$, it follows that $G_\xi(\xi,\tau)>0$ and $G_{\xi\xi}(\xi,\tau) \leq 0$.

\begin{proposition}
\label{G.converges.to.0.smoothly.on.compact.sets}
As $\tau \to -\infty$, $G(\xi,\tau) \to 0$ in $C_{\text{\rm loc}}^\infty$. 
\end{proposition}

\textbf{Proof.} 
This follows from Proposition \ref{asymptotic.shrinking.soliton}. \\

\begin{proposition}
\label{evolution.of.G}
The function $G$ satisfies 
\begin{align*} 
&\Big | G_\tau(\xi,\tau) - G_{\xi\xi}(\xi,\tau) + \frac{1}{2} \, \xi \, G_\xi(\xi,\tau) \\ 
&- \frac{1}{2} \, (\sqrt{2}+G(\xi,\tau)) + (\sqrt{2}+G(\xi,\tau))^{-1} \, (1+G_\xi(\xi,\tau)^2) \Big | \\ 
&\leq 2 \, (\sqrt{2}+G(0,\tau))^{-1} \, G_\xi(0,\tau) \, G_\xi(\xi,\tau) \\ 
&+ 2 \, \Big | \frac{1}{\sqrt{2}+G(\xi,\tau)} - \frac{1}{\sqrt{2}+G(0,\tau)} \Big | \, \max \{G_\xi(\xi,\tau),G_\xi(0,\tau)\} \, G_\xi(\xi,\tau). 
\end{align*}
\end{proposition}

\textbf{Proof.} This follows immediately from Corollary \ref{evolution.of.F.2}. \\

For each $k$, we define 
\begin{align*} 
\delta_k 
&:= \sup_{\tau \leq -k} |G(0,\tau)|+G_\xi(0,\tau) \\ 
&= \sup_{t \leq -e^k} \Big | \frac{F(0,t)}{\sqrt{-t}}-\sqrt{2} \Big |+F_z(0,t) \\ 
&= \sup_{t \leq -e^k} \Big | \frac{\bar{r}(t)}{\sqrt{-t}}-\sqrt{2}|+u(\bar{r}(t),t)^{\frac{1}{2}}. 
\end{align*}
By definition, $\delta_k$ is a decreasing sequence of positive numbers. Moreover, $\delta_k \to 0$ by Proposition \ref{G.converges.to.0.smoothly.on.compact.sets}. 

\begin{lemma}
\label{C1.bound.for.G}
We have $|G(\xi,\tau)|+|G_\xi(\xi,\tau)| \leq C\delta_k^{\frac{1}{4}}$ for $\tau \leq -k$ and $|\xi| \leq 2\delta_k^{-\frac{1}{100}}$. 
\end{lemma}

\textbf{Proof.} 
By definition of $\delta_k$, we have $\big | \frac{\bar{r}(t)}{\sqrt{-2t}}-1 \big | \leq \delta_k$ and $u(\bar{r}(t),t) \leq \delta_k^2$ for all $t \leq -e^k$. We now apply Proposition \ref{maximum.principle} with $\bar{t}=-e^k$ and $a = \frac{1}{10} \, \delta_k^{-\frac{1}{2}}$. Using Proposition \ref{maximum.principle}, we conclude that $u(r,t) \leq C\delta_k$ for all $t \leq -e^k$ and all $\frac{1}{2} \, \sqrt{-2t} \leq r \leq \bar{r}(t)$. This implies $0 \leq G_\xi(\xi,\tau) \leq C\delta_k^{\frac{1}{2}}$ whenever $\tau \leq -k$ and $-\frac{1}{\sqrt{2}} \leq G(\xi,\tau) \leq G(0,\tau)$. Since $|G(0,\tau)| \leq \delta_k$ for $\tau \leq -k$, we conclude that $|G(\xi,\tau)|+|G_\xi(\xi,\tau)| \leq C \delta_k^{\frac{1}{4}}$ for all $\tau \leq -k$ and $-2\delta_k^{-\frac{1}{100}} \leq \xi \leq 0$. 

On the other hand, using the inequality $G_{\xi\xi}(\xi,\tau) \leq 0$, we obtain $0 \leq G_\xi(\xi,\tau) \leq G_\xi(0,\tau) \leq \delta_k$ for all $\tau \leq -k$ and $0 \leq \xi \leq 2\delta_k^{-\frac{1}{100}}$. Since $|G(0,\tau)| \leq \delta_k$, it follows that $|G(\xi,\tau)|+|G_\xi(\xi,\tau)| \leq C \delta_k^{\frac{1}{4}}$ for $\tau \leq -k$ and $0 \leq \xi \leq 2\delta_k^{-\frac{1}{100}}$. This completes the proof of Lemma \ref{C1.bound.for.G}. \\

\begin{lemma}
\label{C2.bound.for.G}
We have $|G_{\xi\xi}(\xi,\tau)| \leq C \delta_k^{\frac{1}{8}}$ for $\tau \leq -k$ and $|\xi| \leq \delta_k^{-\frac{1}{100}}$. 
\end{lemma} 

\textbf{Proof.} 
Applying Proposition \ref{pointwise.derivative.bound.for.F} with $m=2$, we obtain $|G_{\xi\xi\xi}(\xi,\tau)| \leq C \, (1+|G_{\xi\xi}(\xi,\tau)|)^2$  for $\tau \leq -k$ and $|\xi| \leq 2\delta_k^{-\frac{1}{100}}$. Moreover, Lemma \ref{C1.bound.for.G} implies 
\[\inf_{\xi' \in [\xi-\delta_k^{\frac{1}{8}},\xi+\delta_k^{\frac{1}{8}}]} |G_{\xi\xi}(\xi',\tau)| \leq C \delta_k^{\frac{1}{8}}\] 
for $\tau \leq -k$ and $|\xi| \leq \delta_k^{-\frac{1}{100}}$. Putting these facts together, we conclude that $|G_{\xi\xi}(\xi,\tau)| \leq C \delta_k^{\frac{1}{8}}$ for $\tau \leq -k$ and $|\xi| \leq \delta_k^{-\frac{1}{100}}$. This completes the proof of Lemma \ref{C2.bound.for.G}. \\

\begin{lemma}
\label{higher.derivative.bound.for.G}
We have $|\partial_\xi^{m+1} G(\xi,\tau)| \leq C(m)$ for $\tau \leq -k$ and $|\xi| \leq \delta_k^{-\frac{1}{100}}$. 
\end{lemma}

\textbf{Proof.} 
Using Proposition \ref{pointwise.derivative.bound.for.F}, we obtain $|\partial_\xi^{m+1} G(\xi,\tau)| \leq C(m) \, (1+|G_{\xi\xi}(\xi,\tau)|)^m$ for $\tau \leq -k$ and $|\xi| \leq 2\delta_k^{-\frac{1}{100}}$. Moreover, Lemma \ref{C2.bound.for.G} implies $|G_{\xi\xi}(\xi,\tau)| \leq C \delta_k^{\frac{1}{8}}$ for $\tau \leq -k$ and $|\xi| \leq \delta_k^{-\frac{1}{100}}$. Putting these facts together, the assertion follows. \\

\begin{lemma}
\label{integral.estimates.for.G}
We have 
\[|G_\xi(0,\tau)|^4 \leq C \delta_k^{\frac{1}{100}} \int_{\{|\xi| \leq \delta_k^{-\frac{1}{100}}\}} e^{-\frac{\xi^2}{4}} \, |G(\xi,\tau)|^2 \, d\xi\] 
and  
\begin{align*} 
\int_{\{|\xi| \leq \delta_k^{-\frac{1}{100}}\}} e^{-\frac{\xi^2}{4}} \, |G_\xi(\xi,\tau)|^4 \, d\xi 
&\leq C\delta_k^{\frac{1}{100}} \int_{\{|\xi| \leq \delta_k^{-\frac{1}{100}}\}} e^{-\frac{\xi^2}{4}} \, |G(\xi,\tau)|^2 \, d\xi \\ 
&+ C \, \exp(-\frac{1}{8} \, \delta_k^{-\frac{1}{50}}) 
\end{align*}
for $\tau \leq -k$.
\end{lemma}

\textbf{Proof.} 
Using Lemma \ref{higher.derivative.bound.for.G} and standard interpolation inequalities, we obtain 
\begin{align*} 
|G_\xi(0,\tau)|^4 
&\leq C \, \bigg ( \int_{\{|\xi| \leq 1\}} |G(\xi,\tau)|^2 \, d\xi \bigg )^{\frac{3}{2}} \\ 
&\leq C \delta_k^{\frac{1}{100}} \int_{\{|\xi| \leq \delta_k^{-\frac{1}{100}}\}} e^{-\frac{\xi^2}{4}} \, |G(\xi,\tau)|^2 \, d\xi 
\end{align*} 
for $\tau \leq -k$, where in the last step we have used Lemma \ref{C1.bound.for.G}. This proves the first statement. To prove the second statement, we observe that 
\begin{align*} 
&\int_{\{|\xi| \leq \delta_k^{-\frac{1}{100}}\}} e^{-\frac{\xi^2}{4}} \, G_\xi(\xi,\tau)^4 \, d\xi \\ 
&+ 3 \int_{\{|\xi| \leq \delta_k^{-\frac{1}{100}}\}} e^{-\frac{\xi^2}{4}} \, G_\xi(\xi,\tau)^2 \, G(\xi,\tau) \, G_{\xi\xi}(\xi,\tau) \, d\xi \\ 
&- \frac{1}{2} \int_{\{|\xi| \leq \delta_k^{-\frac{1}{100}}\}} e^{-\frac{\xi^2}{4}} \, \xi \, G_\xi(\xi,\tau)^3 \, G(\xi,\tau) \, d\xi \\ 
&= \int_{\{|\xi| \leq \delta_k^{-\frac{1}{100}}\}} \frac{\partial}{\partial \xi} (e^{-\frac{\xi^2}{4}} \, G_\xi(\xi,\tau)^3 \, G(\xi,\tau)) \, d\xi \\ 
&\leq C \, \exp(-\frac{1}{8} \, \delta_k^{-\frac{1}{50}}) 
\end{align*} 
for $\tau \leq -k$, where in the last step we have used Lemma \ref{C1.bound.for.G}. Using Lemma \ref{C2.bound.for.G}, we obtain   
\begin{align*} 
&-3 \int_{\{|\xi| \leq \delta_k^{-\frac{1}{100}}\}} e^{-\frac{\xi^2}{4}} \, G_\xi(\xi,\tau)^2 \, G(\xi,\tau) \, G_{\xi\xi}(\xi,\tau) \, d\xi \\ 
&\leq C \delta_k^{\frac{1}{100}} \int_{\{|\xi| \leq \delta_k^{-\frac{1}{100}}\}} e^{-\frac{\xi^2}{4}} \, G_\xi(\xi,\tau)^2 \, |G(\xi,\tau)| \, d\xi 
\end{align*}
for $\tau \leq -k$. Moreover, Lemma \ref{C1.bound.for.G} implies 
\begin{align*} 
&\frac{1}{2} \int_{\{|\xi| \leq \delta_k^{-\frac{1}{100}}\}} e^{-\frac{\xi^2}{4}} \, \xi \, G_\xi(\xi,\tau)^3 \, G(\xi,\tau) \, d\xi \\ 
&\leq C \delta_k^{\frac{1}{100}} \int_{\{|\xi| \leq \delta_k^{-\frac{1}{100}}\}} e^{-\frac{\xi^2}{4}} \, G_\xi(\xi,\tau)^2 \, |G(\xi,\tau)| \, d\xi 
\end{align*} 
for $\tau \leq -k$. Adding these inequalities gives 
\begin{align*} 
&\int_{\{|\xi| \leq \delta_k^{-\frac{1}{100}}\}} e^{-\frac{\xi^2}{4}} \, G_\xi(\xi,\tau)^4 \, d\xi \\ 
&\leq C\delta_k^{\frac{1}{100}} \int_{\{|\xi| \leq \delta_k^{-\frac{1}{100}}\}} e^{-\frac{\xi^2}{4}} \, G_\xi(\xi,\tau)^2 \, |G(\xi,\tau)| \, d\xi \\ 
&+ C \, \exp(-\frac{1}{8} \, \delta_k^{-\frac{1}{50}}) \\ 
&\leq C\delta_k^{\frac{1}{100}} \int_{\{|\xi| \leq \delta_k^{-\frac{1}{100}}\}} e^{-\frac{\xi^2}{4}} \, (G_\xi(\xi,\tau)^4 + G(\xi,\tau)^2) \, d\xi \\ 
&+ C \, \exp(-\frac{1}{8} \, \delta_k^{-\frac{1}{50}}) 
\end{align*}
for $\tau \leq -k$. Rearranging terms, the assertion follows. \\

\begin{lemma}
\label{evolution.of.G.2}
We have 
\begin{align*} 
&\int_{\{|\xi| \leq \delta_k^{-\frac{1}{100}}\}} e^{-\frac{\xi^2}{4}} \, \Big | G_\tau(\xi,\tau) - G_{\xi\xi}(\xi,\tau) + \frac{1}{2} \, \xi \, G_\xi(\xi,\tau) - G(\xi,\tau) \Big |^2 \, d\xi \\ 
&\leq C\delta_k^{\frac{1}{100}} \int_{\{|\xi| \leq \delta_k^{-\frac{1}{100}}\}} e^{-\frac{\xi^2}{4}} \, |G(\xi,\tau)|^2 \, d\xi + C \, \exp(-\frac{1}{8} \, \delta_k^{-\frac{1}{50}}) 
\end{align*} 
for $\tau \leq -k$. 
\end{lemma}

\textbf{Proof.} 
Note that 
\[\Big | G(\xi,\tau) - \frac{1}{2} \, (\sqrt{2}+G(\xi,\tau)) + (\sqrt{2}+G(\xi,\tau))^{-1} \Big | \leq C \, G(\xi,\tau)^2\] 
for $\tau \leq -k$ and $|\xi| \leq \delta_k^{-\frac{1}{100}}$. Using Proposition \ref{evolution.of.G}, we obtain the pointwise estimate 
\begin{align*} 
&\Big | G_\tau(\xi,\tau) - G_{\xi\xi}(\xi,\tau) + \frac{1}{2} \, \xi \, G_\xi(\xi,\tau) - G(\xi,\tau) \Big | \\ 
&\leq C \, G(\xi,\tau)^2 + C \, G_\xi(\xi,\tau)^2 + C \, G_\xi(0,\tau) \, G_\xi(\xi,\tau) \\ 
&\leq C \, \delta_k^{\frac{1}{100}} \, G(\xi,\tau) + C \, G_\xi(\xi,\tau)^2 + C \, G_\xi(0,\tau)^2 
\end{align*}
for $\tau \leq -k$ and $|\xi| \leq \delta_k^{-\frac{1}{100}}$. Hence, the assertion follows from Lemma \ref{integral.estimates.for.G}. \\ 

We now perform a spectral decomposition for the operator $G_{\xi\xi} - \frac{1}{2} \, \xi \, G_\xi + G$. This operator is symmetric with respect to the inner product $\|G\|_{\mathcal{H}}^2 = \int_{\mathbb{R}} e^{-\frac{|\xi|^2}{4}} \, G^2 \, d\xi$. The eigenvalues of this operator are given by $1-\frac{n}{2}$, where $n \geq 0$. Moreover, the associated eigenfunctions are given by $H_n(\frac{\xi}{2})$, where $H_n$ is the $n$-th Hermite polynomial. Let us write $\mathcal{H} = \mathcal{H}_+ \oplus \mathcal{H}_0 \oplus \mathcal{H}_-$, where the subspace $\mathcal{H}_+$ is defined as the span of $H_0(\frac{\xi}{2})$ and $H_1(\frac{\xi}{2})$,  the subspace $\mathcal{H}_0$ is defined as the span of $H_2(\frac{\xi}{2})$, and $\mathcal{H}_-$ is the orthogonal complement of $\mathcal{H}_+ \oplus \mathcal{H}_0$. Moreover, let $P_+$, $P_0$, and $P_-$ denote the orthogonal projections associated with the direct sum $\mathcal{H} = \mathcal{H}_+ \oplus \mathcal{H}_0 \oplus \mathcal{H}_-$. The eigenvalues of the operator $-G_{\xi\xi} + \frac{1}{2} \, \xi \, G_\xi - G$ on $\mathcal{H}_+$ are bounded from above by $-\frac{1}{2}$. Similarly, the eigenvalues of the operator $-G_{\xi\xi} + \frac{1}{2} \, \xi \, G_\xi - G$ on $\mathcal{H}_-$ are bounded from below by $\frac{1}{2}$.

Let $\chi$ denote a smooth cutoff function satisfying $\chi(s)=1$ for $s \in [-\frac{1}{2},\frac{1}{2}]$, $\chi(s)=0$ for $s \in \mathbb{R} \setminus [-1,1]$, and $s \chi'(s) \leq 0$ for all $s \in \mathbb{R}$. We define
\begin{align*} 
\gamma_j &:= \sup_{\tau \in [-j-1,-j]} \int_{\mathbb{R}} e^{-\frac{\xi^2}{4}} \, |G(\xi,\tau) \, \chi(\delta_j^{\frac{1}{100}} \xi)|^2 \, d\xi, \\ 
\gamma_j^+ &:= \sup_{\tau \in [-j-1,-j]} \int_{\mathbb{R}} e^{-\frac{\xi^2}{4}} \, |P_+ (G(\xi,\tau) \, \chi(\delta_j^{\frac{1}{100}} \xi))|^2 \, d\xi, \\ 
\gamma_j^0 &:= \sup_{\tau \in [-j-1,-j]} \int_{\mathbb{R}} e^{-\frac{\xi^2}{4}} \, |P_0 (G(\xi,\tau) \, \chi(\delta_j^{\frac{1}{100}} \xi))|^2 \, d\xi, \\ 
\gamma_j^- &:= \sup_{\tau \in [-j-1,-j]} \int_{\mathbb{R}} e^{-\frac{\xi^2}{4}} \, |P_- (G(\xi,\tau) \, \chi(\delta_j^{\frac{1}{100}} \xi))|^2 \, d\xi.
\end{align*}
Clearly, $\frac{1}{C} \, \gamma_j \leq \gamma_j^+ + \gamma_j^0 + \gamma_j^- \leq C \, \gamma_j$. Using Lemma \ref{C1.bound.for.G}, we obtain 
\[\gamma_j \leq C \sup_{\tau \in [-j-1,-j]} \sup_{|\xi| \leq \delta_j^{-\frac{1}{100}}} |G(\xi,\tau)|^2 \leq C \delta_j^{\frac{1}{4}}.\] 
In particular, $\gamma_j \to 0$. \\

\begin{lemma}
\label{discrete.dynamical.system}
We have
\begin{align*} 
&\gamma_{j+1}^+ \leq e^{-1} \, \gamma_j^+ + C\delta_j^{\frac{1}{200}} \, (\gamma_j+\gamma_{j+1}) + C \, \exp(-\frac{1}{64} \, \delta_j^{-\frac{1}{50}}), \\ 
&|\gamma_{j+1}^0 - \gamma_j^0| \leq C\delta_j^{\frac{1}{200}} \, (\gamma_j+\gamma_{j+1}) + C \, \exp(-\frac{1}{64} \, \delta_j^{-\frac{1}{50}}), \\ 
&\gamma_{j+1}^- \geq e \, \gamma_j^- - C\delta_j^{\frac{1}{200}} \, (\gamma_j+\gamma_{j+1}) - C \, \exp(-\frac{1}{64} \, \delta_j^{-\frac{1}{50}}).
\end{align*} 
\end{lemma} 

\textbf{Proof.} 
Fix $j$, and define $\hat{G}(\xi,\tau) := G(\xi,\tau) \, \chi(\delta_j^{\frac{1}{100}} \xi)$. Note that 
\[\int_{\mathbb{R}} e^{-\frac{\xi^2}{4}} \, |\hat{G}(\xi,\tau)|^2 \, d\xi \leq \gamma_j+\gamma_{j+1}\] 
for $\tau \in [-j-2,-j]$. Using Lemma \ref{evolution.of.G.2} and Lemma \ref{C1.bound.for.G}, we obtain 
\begin{align*} 
&\int_{\mathbb{R}} e^{-\frac{\xi^2}{4}} \, \Big | \hat{G}_\tau(\xi,\tau) - \hat{G}_{\xi\xi}(\xi,\tau) + \frac{1}{2} \, \xi \, \hat{G}_\xi(\xi,\tau) - \hat{G}(\xi,\tau) \Big |^2 \, d\xi \\ 
&\leq C\delta_j^{\frac{1}{100}} \, (\gamma_j+\gamma_{j+1}) + C \, \exp(-\frac{1}{32} \, \delta_j^{-\frac{1}{50}}) 
\end{align*} 
for $\tau \in [-j-2,-j]$. Consequently, 
\begin{align*} 
\frac{d}{d\tau} \bigg ( \int_{\mathbb{R}} e^{-\frac{\xi^2}{4}} \, |P_+  \hat{G}(\xi,\tau)|^2 \, d\xi \bigg ) 
&\geq \int_{\mathbb{R}} e^{-\frac{\xi^2}{4}} \, |P_+ \hat{G}(\xi,\tau)|^2 \, d\xi \\ 
&- C\delta_j^{\frac{1}{200}} \, (\gamma_j+\gamma_{j+1}) - C \, \exp(-\frac{1}{64} \, \delta_j^{-\frac{1}{50}}),  
\end{align*}
\begin{align*} 
&\bigg | \frac{d}{d\tau} \bigg ( \int_{\mathbb{R}} e^{-\frac{\xi^2}{4}} \, |P_0  \hat{G}(\xi,\tau)|^2 \, d\xi \bigg ) \bigg | \leq C\delta_j^{\frac{1}{200}} \, (\gamma_j+\gamma_{j+1}) + C \, \exp(-\frac{1}{64} \, \delta_j^{-\frac{1}{50}}),  
\end{align*}
\begin{align*} 
\frac{d}{d\tau} \bigg ( \int_{\mathbb{R}} e^{-\frac{\xi^2}{4}} \, |P_-  \hat{G}(\xi,\tau)|^2 \, d\xi \bigg ) 
&\leq -\int_{\mathbb{R}} e^{-\frac{\xi^2}{4}} \, |P_- \hat{G}(\xi,\tau)|^2 \, d\xi \\ 
&+ C\delta_j^{\frac{1}{200}} \, (\gamma_j+\gamma_{j+1}) + C \, \exp(-\frac{1}{64} \, \delta_j^{-\frac{1}{50}})  
\end{align*}
for $\tau \in [-j-2,-j]$. Integrating these inequalities over the interval $[\tau-1,\tau]$ gives 
\begin{align*} 
\int_{\mathbb{R}} e^{-\frac{\xi^2}{4}} \, |P_+ \hat{G}(\xi,\tau-1)|^2 \, d\xi 
&\leq e^{-1} \int_{\mathbb{R}} e^{-\frac{\xi^2}{4}} \, |P_+ \hat{G}(\xi,\tau)|^2 \, d\xi \\ 
&+ C\delta_j^{\frac{1}{200}} \, (\gamma_j+\gamma_{j+1}) + C \, \exp(-\frac{1}{64} \, \delta_j^{-\frac{1}{50}}), 
\end{align*} 
\begin{align*} 
&\bigg | \int_{\mathbb{R}} e^{-\frac{\xi^2}{4}} \, |P_0 \hat{G}(\xi,\tau-1)|^2 \, d\xi - \int_{\mathbb{R}} e^{-\frac{\xi^2}{4}} \, |P_0 \hat{G}(\xi,\tau)|^2 \, d\xi \bigg | \\ 
&\leq C\delta_j^{\frac{1}{200}} \, (\gamma_j+\gamma_{j+1}) + C \, \exp(-\frac{1}{64} \, \delta_j^{-\frac{1}{50}}), 
\end{align*} 
\begin{align*} 
\int_{\mathbb{R}} e^{-\frac{\xi^2}{4}} \, |P_- \hat{G}(\xi,\tau-1)|^2 \, d\xi 
&\geq e \int_{\mathbb{R}} e^{-\frac{\xi^2}{4}} \, |P_- \hat{G}(\xi,\tau)|^2 \, d\xi \\ 
&- C\delta_j^{\frac{1}{200}} \, (\gamma_j+\gamma_{j+1}) - C \, \exp(-\frac{1}{64} \, \delta_j^{-\frac{1}{50}}) 
\end{align*} 
for $\tau \in [-j-1,-j]$. We now define $\tilde{G}(\xi,\tau) := G(\xi,\tau) \, \chi(\delta_{j+1}^{\frac{1}{100}} \xi)$. Using Lemma \ref{C1.bound.for.G}, we obtain 
\begin{align*} 
&\int_{\mathbb{R}} e^{-\frac{\xi^2}{4}} \, |\tilde{G}(\xi,\tau-1) - \hat{G}(\xi,\tau-1)|^2 \, d\xi \\ 
&\leq \int_{\{\frac{1}{2} \delta_j^{-\frac{1}{100}} \leq |\xi| \leq \delta_{j+1}^{-\frac{1}{100}}\}} e^{-\frac{\xi^2}{4}} \, |G(\xi,\tau-1)|^2 \, d\xi \leq C \, \exp(-\frac{1}{32} \, \delta_j^{-\frac{1}{50}}) 
\end{align*} 
for $\tau \in [-j-1,-j]$. Putting these facts together, we conclude that 
\begin{align*} 
\int_{\mathbb{R}} e^{-\frac{\xi^2}{4}} \, |P_+ \tilde{G}(\xi,\tau-1)|^2 \, d\xi 
&\leq e^{-1} \int_{\mathbb{R}} e^{-\frac{\xi^2}{4}} \, |P_+ \hat{G}(\xi,\tau)|^2 \, d\xi \\ 
&+ C\delta_j^{\frac{1}{200}} \, (\gamma_j+\gamma_{j+1}) + C \, \exp(-\frac{1}{64} \, \delta_j^{-\frac{1}{50}}), 
\end{align*} 
\begin{align*} 
&\bigg | \int_{\mathbb{R}} e^{-\frac{\xi^2}{4}} \, |P_0 \tilde{G}(\xi,\tau-1)|^2 \, d\xi - \int_{\mathbb{R}} e^{-\frac{\xi^2}{4}} \, |P_0 \hat{G}(\xi,\tau)|^2 \, d\xi \bigg | \\ 
&\leq C\delta_j^{\frac{1}{200}} \, (\gamma_j+\gamma_{j+1}) + C \, \exp(-\frac{1}{64} \, \delta_j^{-\frac{1}{50}}), 
\end{align*} 
\begin{align*} 
\int_{\mathbb{R}} e^{-\frac{\xi^2}{4}} \, |P_- \tilde{G}(\xi,\tau-1)|^2 \, d\xi 
&\geq e \int_{\mathbb{R}} e^{-\frac{\xi^2}{4}} \, |P_- \hat{G}(\xi,\tau)|^2 \, d\xi \\ 
&- C\delta_j^{\frac{1}{200}} \, (\gamma_j+\gamma_{j+1}) - C \, \exp(-\frac{1}{64} \, \delta_j^{-\frac{1}{50}}) 
\end{align*} 
for $\tau \in [-j-1,-j]$. Taking the supremum over $\tau \in [-j-1,-j]$, the assertion follows. This completes the proof of Lemma \ref{discrete.dynamical.system}. \\ 

We next define 
\begin{align*} 
&\Gamma_k := \sup_{j \geq k} \gamma_j, \\ 
&\Gamma_k^+ := \sup_{j \geq k} \gamma_j^+, \\ 
&\Gamma_k^0 := \sup_{j \geq k} \gamma_j^0, \\ 
&\Gamma_k^- := \sup_{j \geq k} \gamma_j^-. 
\end{align*} 
Clearly, $\frac{1}{C} \, \Gamma_k \leq \Gamma_k^+ + \Gamma_k^0 + \Gamma_k^- \leq C \, \Gamma_k$. The inequality $\gamma_j \leq C \delta_j^{\frac{1}{4}}$ gives $\Gamma_k \leq C \delta_k^{\frac{1}{4}}$. In particular, $\Gamma_k \to 0$. Using Lemma \ref{discrete.dynamical.system}, we obtain 
\begin{align*} 
&\Gamma_{k+1}^+ \leq e^{-1} \, \Gamma_k^+ + C\delta_k^{\frac{1}{200}} \, \Gamma_k + C \, \exp(-\frac{1}{64} \, \delta_k^{-\frac{1}{50}}), \\ 
&|\Gamma_{k+1}^0 - \Gamma_k^0| \leq C\delta_k^{\frac{1}{200}} \, \Gamma_k + C \, \exp(-\frac{1}{64} \, \delta_k^{-\frac{1}{50}}), \\ 
&\Gamma_{k+1}^- \geq e \, \Gamma_k^- - C\delta_k^{\frac{1}{200}} \, \Gamma_k - C \, \exp(-\frac{1}{64} \, \delta_k^{-\frac{1}{50}}).
\end{align*} 
On the other hand, it follows from Lemma \ref{higher.derivative.bound.for.G} and standard interpolation inequalities that 
\[\sup_{\tau \in [-j-1,-j]}  |G(0,\tau)|+G_\xi(0,\tau) \leq C \, \gamma_j^{\frac{1}{4}},\] 
hence 
\[\delta_k = \sup_{\tau \leq -k} |G(0,\tau)|+G_\xi(0,\tau) \leq C \, \Gamma_k^{\frac{1}{4}}.\] 
Consequently, $\exp(-\frac{1}{64} \, \delta_k^{-\frac{1}{50}}) \leq C\delta_k^5 \leq C\delta_k \Gamma_k$. Putting these facts together, we conclude that 
\begin{align*} 
&\Gamma_{k+1}^+ \leq e^{-1} \, \Gamma_k^+ + C\delta_k^{\frac{1}{200}} \, \Gamma_k, \\ 
&|\Gamma_{k+1}^0 - \Gamma_k^0| \leq C\delta_k^{\frac{1}{200}} \, \Gamma_k, \\ 
&\Gamma_{k+1}^- \geq e \, \Gamma_k^- - C\delta_k^{\frac{1}{200}} \, \Gamma_k.
\end{align*} 
The following lemma is inspired by a lemma of Merle and Zaag (cf. \cite{Merle-Zaag}, Lemma A.1):

\begin{lemma}
\label{discrete.merle.zaag.lemma}
We either have $\Gamma_k^0 + \Gamma_k^- \leq o(1) \, \Gamma_k^+$, or $\Gamma_k^++\Gamma_k^- \leq o(1) \, \Gamma_k^0$. 
\end{lemma}

\textbf{Proof.} 
By definition, the sequence $\Gamma_k^-$ is monotone decreasing. This implies $\Gamma_k^- \geq \Gamma_{k+1}^- \geq e \, \Gamma_k^- - o(1) \, \Gamma_k$. Thus, $\Gamma_k^- \leq o(1) \, \Gamma_k$. This gives $\Gamma_k^- \leq o(1) \, (\Gamma_k^++\Gamma_k^0)$.

Let $I$ denote the set of all positive real numbers $\alpha$ with the property that the set $\{k: \Gamma_k^0 < \alpha \, \Gamma_k^+\}$ is finite. Moreover, let $J$ denote the set of all positive real numbers $\alpha$ with the property that the set $\{k: \Gamma_k^0 \geq \alpha \, \Gamma_k^+\}$ is infinite. Clearly, $I \subset J$. 

We claim that $e^{\frac{1}{2}} \alpha \in I$ whenever $\alpha \in J$. To see this, suppose that $\alpha \in J$. We can find a large integer $k_0$ (depending on $\alpha$) such that 
\[\Gamma_{k+1}^+ \leq e^{-1} \, \Gamma_k^+ + \frac{1}{2(1+\alpha)} \, (e^{-\frac{1}{2}}-e^{-1}) \, (\Gamma_k^+ + \Gamma_k^0)\] 
and 
\[|\Gamma_{k+1}^0 - \Gamma_k^0| \leq \frac{\alpha}{2(1+\alpha)} \, (1-e^{-\frac{1}{2}}) \, (\Gamma_k^+ + \Gamma_k^0)\] 
for all $k \geq k_0$. This implies 
\begin{align*} 
\Gamma_{k+1}^0 - e^{\frac{1}{2}} \alpha \, \Gamma_{k+1}^+ 
&\geq \Gamma_k^0 - e^{-\frac{1}{2}} \alpha \, \Gamma_k^+ - \frac{\alpha}{1+\alpha} \, (1-e^{-\frac{1}{2}}) \, (\Gamma_k^++\Gamma_k^0) \\ 
&= \Big ( 1 - \frac{\alpha}{1+\alpha} \, (1-e^{-\frac{1}{2}}) \Big ) \, (\Gamma_k^0 - \alpha \, \Gamma_k^+) 
\end{align*} 
for all $k \geq k_0$. Since $\alpha \in J$, the set $\{k: \Gamma_k^0 \geq \alpha \, \Gamma_k^+\}$ is infinite. Hence, we can find an integer $k_1 \geq k_0$ such that $\Gamma_k^0 - \alpha \, \Gamma_k^+ \geq 0$ for $k=k_1$. Proceeding inductively, we obtain $\Gamma_{k+1}^0 - e^{\frac{1}{2}} \alpha \, \Gamma_{k+1}^+ \geq 0$ for all $k \geq k_1$. Consequently, the set $\{k: \Gamma_k^0 < e^{\frac{1}{2}} \alpha \, \Gamma_k^+\}$ is finite. Thus, $e^{\frac{1}{2}} \alpha \in I$. This proves the claim. 

Thus, we conclude that either $J = \emptyset$ or $I = (0,\infty)$. If $I = (0,\infty)$, we obtain $\Gamma_k^+ \leq o(1) \, \Gamma_k^0$, hence $\Gamma_k^++\Gamma_k^- \leq o(1) \, \Gamma_k^0$. On the other hand, if $J = \emptyset$, then $\Gamma_k^0 \leq o(1) \, \Gamma_k^+$, hence $\Gamma_k^0+\Gamma_k^- \leq o(1) \, \Gamma_k^+$. This completes the proof of Lemma \ref{discrete.merle.zaag.lemma}. \\

In the next step, we show that the second possibility in Lemma \ref{discrete.merle.zaag.lemma} cannot occur. \\

\begin{lemma}
\label{growing.mode.dominates}
We have $\Gamma_k^0 + \Gamma_k^- \leq o(1) \, \Gamma_k^+$. 
\end{lemma} 

\textbf{Proof.} 
Suppose that the assertion is false. In view of Lemma \ref{discrete.merle.zaag.lemma}, we have $\Gamma_k^++\Gamma_k^- \leq o(1) \, \Gamma_k^0$. For each $k$, we can find an integer $j_k \geq k$ and a time $\tau_k \in [-j_k-1,-j_k]$ such that 
\[\Gamma_k = \gamma_{j_k} = \int_{\mathbb{R}} e^{-\frac{\xi^2}{4}} \, |G(\xi,\tau_k) \, \chi(\delta_{j_k}^{\frac{1}{100}} \xi)|^2 \, d\xi.\] 
Note that 
\[\int_{\mathbb{R}} e^{-\frac{\xi^2}{4}} \, |P_+(G(\xi,\tau_k) \, \chi(\delta_{j_k}^{\frac{1}{100}} \xi))|^2 \, d\xi \leq \gamma_{j_k}^+ \leq \Gamma_k^+ \leq o(1) \, \Gamma_k\] 
and 
\[\int_{\mathbb{R}} e^{-\frac{\xi^2}{4}} \, |P_-(G(\xi,\tau_k) \, \chi(\delta_{j_k}^{\frac{1}{100}} \xi))|^2 \, d\xi \leq \gamma_{j_k}^- \leq \Gamma_k^- \leq o(1) \, \Gamma_k.\] 
After passing to a subsequence, the functions $\xi \mapsto \Gamma_k^{-\frac{1}{2}} \, G(\xi,\tau_k) \, \chi(\delta_{j_k}^{\frac{1}{100}} \xi)$ converge, in $\mathcal{H}$, to a non-zero multiple of the function $\xi^2-2$. Since the function $\xi \mapsto G(\xi,\tau_k)$ is monotone increasing, we have 
\[\int_{-3}^{-1} G(\xi,\tau_k) \, d\xi \leq \int_{-1}^1 G(\xi,\tau_k) \, d\xi \leq \int_1^3 G(\xi,\tau_k) \, d\xi\] 
for each $k$. Passing to the limit as $k \to \infty$, we obtain either 
\[\int_{-3}^{-1} (\xi^2-2) \, d\xi \leq \int_{-1}^1 (\xi^2-2) \, d\xi \leq \int_1^3 (\xi^2-2) \, d\xi\] 
or 
\[\int_{-3}^{-1} (2-\xi^2) \, d\xi \leq \int_{-1}^1 (2-\xi^2) \, d\xi \leq \int_1^3 (2-\xi^2) \, d\xi.\] 
In either case, we arrive at a contradiction. This completes the proof of Lemma \ref{growing.mode.dominates}. \\

\begin{lemma} 
\label{asymptotics.of.Gamma_k}
We have $\Gamma_k \leq O(e^{-k})$.
\end{lemma}

\textbf{Proof.} 
Note that $\Gamma_k^0 + \Gamma_k^- \leq o(1) \, \Gamma_k^+$ by Lemma \ref{growing.mode.dominates}. This implies 
\[\Gamma_{k+1}^+ \leq e^{-1} \, \Gamma_k^+ + C\delta_k^{\frac{1}{200}} \, \Gamma_k^+ \leq e^{-\frac{1}{2}} \, \Gamma_k^+\] 
if $k$ is sufficiently large. Iterating this estimate gives $\Gamma_k^+ \leq O(e^{-\frac{k}{2}})$, hence $\Gamma_k \leq O(e^{-\frac{k}{2}})$. Using the estimate $\delta_k \leq C \, \Gamma_k^{\frac{1}{4}}$, we obtain $\delta_k \leq O(e^{-\frac{k}{8}})$. This gives 
\[\Gamma_{k+1}^+ \leq e^{-1} \, \Gamma_k^+ + C\delta_k^{\frac{1}{200}} \, \Gamma_k^+ \leq e^{-1} \, \Gamma_k^+ + e^{-\frac{k}{2000}} \, \Gamma_k^+\] 
if $k$ is sufficiently large. Iterating this estimate, we conclude that $\Gamma_k^+ \leq O(e^{-k})$, hence $\Gamma_k \leq O(e^{-k})$. This completes the proof of Lemma \ref{asymptotics.of.Gamma_k}. \\

\begin{lemma}
\label{asymptotics.of.G}
We have $|G(0,\tau)| \leq O(e^{-\frac{k}{2}})$ and $|G_\xi(0,\tau)| \leq O(e^{-\frac{k}{2}})$ for all $\tau \leq -k$.
\end{lemma}

\textbf{Proof.} 
Lemma \ref{asymptotics.of.Gamma_k} gives 
\[\int_{\{|\xi| \leq 2\}} |G(\xi,\tau)|^2 \, d\xi \leq O(e^{-k})\] 
for all $\tau \leq -k$. Using Lemma \ref{higher.derivative.bound.for.G} and standard interpolation inequalities, we obtain 
\[\sup_{|\xi| \leq 1} |G(\xi,\tau)|+|G_\xi(\xi,\tau)| \leq O(e^{-\frac{k}{3}})\] 
for all $\tau \leq -k$. Hence, Proposition \ref{evolution.of.G} implies 
\[\sup_{|\xi| \leq 1} \Big | G_\tau(\xi,\tau) - G_{\xi\xi}(\xi,\tau) + \frac{1}{2} \, \xi \, G_\xi(\xi,\tau) - G(\xi,\tau) \Big | \leq O(e^{-\frac{2k}{3}})\] 
for all $\tau \leq -k$. Using standard interior estimates for linear parabolic equations, we conclude that $|G(0,\tau)| \leq O(e^{-\frac{k}{2}})$ and $|G_\xi(0,\tau)| \leq O(e^{-\frac{k}{2}})$ for all $\tau \leq -k$. This completes the proof of Lemma \ref{asymptotics.of.G}. \\

After these preparations, we now prove the main result of this section:

\begin{proposition}
\label{lower.bounds.for.d(t).and.R_max(t)}
The function $d(t)$ satisfies $\liminf_{t \to -\infty} (-t)^{-1} \, d(t) > 0$. Moreover, $\liminf_{t \to -\infty} R_{\text{\rm max}}(t) > 0$, where $R_{\text{\rm max}}(t)$ denotes the supremum of the scalar curvature of $(M,g(t))$.
\end{proposition}

\textbf{Proof.} 
By Lemma \ref{asymptotics.of.G}, we have $|G(0,\tau)| \leq O(e^{\frac{\tau}{2}})$ and $G_\xi(0,\tau) \leq O(e^{\frac{\tau}{2}})$. Changing variables gives $|F(0,t)-\sqrt{-2t}| \leq O(1)$ and $F_z(0,t) \leq O(\frac{1}{\sqrt{-t}})$. Since $F(0,t)=\bar{r}(t)$ and $F_z(0,t)=u(\bar{r}(t),t)^{\frac{1}{2}}$, we obtain $|\bar{r}(t)-\sqrt{-2t}| \leq O(1)$ and $u(\bar{r}(t),t) \leq O(\frac{1}{(-t)})$. Applying Proposition \ref{distance.of.reference.point.to.tip}, we conclude that 
\[\liminf_{t \to -\infty} (-t)^{-1} \, d(t) = \liminf_{t \to -\infty} (-t)^{-1} \int_0^{\bar{r}(t)} u(r,t)^{-\frac{1}{2}} \, dr > 0.\] 
We next observe that $d(t) = d_{g(t)}(p,q)$, where $p$ denotes the tip and $q$ is a fixed reference point on the manifold. Using Lemma 8.3(b) in \cite{Perelman1}, we can control how fast the geodesic distance of $p$ and $q$ can grow as we go backwards in time:
\[-\frac{d}{dt} d_{g(t)}(p,q) \leq C \, R_{\text{\rm max}}(t)^{\frac{1}{2}}.\] 
Since $\liminf_{t \to -\infty} (-t)^{-1} \, d_{g(t)}(p,q) > 0$, it follows that $\limsup_{t \to -\infty} R_{\text{\rm max}}(t) > 0$. Since the function $t \mapsto R_{\text{\rm max}}(t)$ is monotone increasing by Hamilton's Harnack inequality \cite{Hamilton2}, we conclude that $\liminf_{t \to -\infty} R_{\text{\rm max}}(t) > 0$. \\

\section{Uniqueness of ancient $\kappa$-solutions with rotational symmetry}

\label{uniqueness.rot.sym}

We continue to assume that $(M,g(t))$ is a three-dimensional ancient $\kappa$-solution which is noncompact, has positive sectional curvature, and is rotationally symmetric.

\begin{proposition}
\label{limit.of.scalar.curvature.at.tip}
Let $p$ denote the tip. Then $\liminf_{t \to -\infty} R(p,t) > 0$.
\end{proposition}

\textbf{Proof.} 
Since the traceless Ricci tensor vanishes at the tip, the tip cannot lie on a neck. Hence, it follows from work of Perelman \cite{Perelman1} that $R_{\text{\rm max}}(t) \leq C \, R(p,t)$ for some uniform constant $C$ (see Corollary \ref{bound.for.R_max} below). Using Proposition \ref{lower.bounds.for.d(t).and.R_max(t)}, we obtain $\liminf_{t \to -\infty} R(p,t) > 0$. This completes the proof of Proposition \ref{limit.of.scalar.curvature.at.tip}. \\

Let $p$ denote the tip. By Hamilton's trace Harnack inequality \cite{Hamilton2}, the function $t \mapsto R(p,t)$ is monotone increasing. Hence, the limit 
\[\mathcal{R} := \lim_{t \to -\infty} R(p,t)\] 
exists. Moreover, $\mathcal{R}>0$ by Proposition \ref{limit.of.scalar.curvature.at.tip}.

\begin{proposition}
\label{convergence.to.Bryant.soliton}
If we dilate $(M,g(t))$ around the tip by the factor $\mathcal{R}^{\frac{1}{2}}$, then the rescaled manifolds converge to the Bryant soliton in the Cheeger-Gromov sense.
\end{proposition}

\textbf{Proof.} 
Let $p$ denote the tip, and let $t_k$ be a sequence of times such that $t_k \to -\infty$. Let us dilate the flow around the point $(p,t_k)$ by the factor $\mathcal{R}^{\frac{1}{2}}$. The rescaled flows have uniformly bounded curvature. Hence, the rescaled flows converge in the Cheeger-Gromov sense to an eternal solution which is rotationally symmetric. Moreover, on the limiting eternal solution, the scalar curvature at the tip is equal to $1$ at all times. Therefore, the limiting solution attains equality in Hamilton's Harnack inequality \cite{Hamilton2}. Consequently, the limit must be a steady gradient Ricci soliton \cite{Hamilton3}. Therefore, the limit must be the Bryant soliton. This completes the proof of Proposition \ref{convergence.to.Bryant.soliton}. \\

We will need the following basic fact about the Bryant soliton:

\begin{lemma}
\label{property.of.bryant.soliton}
Consider the Bryant soliton, normalized so that the scalar curvature at the tip is equal to $1$. Let $\gamma$ be a geodesic ray emanating from the tip of the Bryant soliton which is parametrized by arclength. Then $\int_0^\infty \text{\rm Ric}(\gamma'(s),\gamma'(s)) \, ds = 1$.
\end{lemma}

\textbf{Proof.} 
On the Bryant soliton, we may write $\text{\rm Ric} = D^2 f$. This implies 
\[\frac{d}{ds} \langle \nabla f(\gamma(s)),\gamma'(s) \rangle = (D^2 f)(\gamma'(s),\gamma'(s)) = \text{\rm Ric}(\gamma'(s),\gamma'(s)).\] 
Clearly, $\nabla f=0$ at the tip. Moreover, the identity $R+|\nabla f|^2=1$ implies that $|\nabla f| \to 1$ at infinity. Consequently, $\langle \nabla f(\gamma(s)),\gamma'(s) \rangle = |\nabla f(\gamma(s))| \to 1$ as $s \to \infty$. Thus, $\int_0^\infty \text{\rm Ric}(\gamma'(s),\gamma'(s)) \, ds = 1$. This completes the proof of Lemma \ref{property.of.bryant.soliton}. \\

We now continue with the the analysis of our ancient solution. As in Section \ref{asymptotics}, we define 
\[d(t) = \int_0^{\bar{r}(t)} u(r,t)^{-\frac{1}{2}} \, dr.\] 
Equivalently, we may write $d(t) = d_{g(t)}(p,q)$, where $p$ denotes the tip and $q$ denotes the reference point introduced in Section \ref{asymptotics}. Clearly, $-d'(t) > 0$. 

\begin{lemma}
\label{derivative.of.d(t)}
Let $\delta>0$ be given. Then $(1-\delta) \, \mathcal{R}^{\frac{1}{2}} \leq -d'(t) \leq (1+\delta) \, \mathcal{R}^{\frac{1}{2}}$ if $-t$ is sufficiently large. 
\end{lemma}

\textbf{Proof.} 
Let $p$ denote the tip, and let $\gamma$ denote the unit-speed geodesic in $(M,g(t))$ from the tip $p$ to our reference point $q$, so that $\gamma(0)=p$ and $\gamma(d(t))=q$. In view of Lemma \ref{property.of.bryant.soliton} and Proposition \ref{convergence.to.Bryant.soliton}, we can find a large constant $A$ (depending on $\delta$) such that $A \geq 8\delta^{-1}$ and 
\[(1-\delta) \, \mathcal{R}^{\frac{1}{2}} \leq \int_0^{A \mathcal{R}^{-\frac{1}{2}}} \text{\rm Ric}_{g(t)}(\gamma'(s),\gamma'(s)) \, ds \leq (1+\frac{\delta}{2}) \, \mathcal{R}^{\frac{1}{2}}\] 
if $-t$ is sufficiently large (depending on $\delta$ and $A$). 

We next observe that $\gamma$ is part of a minimizing geodesic ray emanating from the tip $p$. Hence, we may apply Theorem 17.4(a) in \cite{Hamilton4} with $\sigma = A \mathcal{R}^{-\frac{1}{2}}$ and $L = d(t) + A \mathcal{R}^{-\frac{1}{2}}$. This gives 
\[0 \leq \int_{A \mathcal{R}^{-\frac{1}{2}}}^{d(t)} \text{\rm Ric}_{g(t)}(\gamma'(s),\gamma'(s)) \, ds \leq 4A^{-1} \, \mathcal{R}^{\frac{1}{2}}.\] 
Putting these facts together, we obtain 
\[(1-\delta) \, \mathcal{R}^{\frac{1}{2}} \leq \int_0^{d(t)} \text{\rm Ric}_{g(t)}(\gamma'(s),\gamma'(s)) \, ds \leq (1+\frac{\delta}{2}+4A^{-1}) \, \mathcal{R}^{\frac{1}{2}}\] 
if $-t$ is sufficiently large (depending on $\delta$ and $A$). Since $d'(t) = -\int_0^{d(t)} \text{\rm Ric}_{g(t)}(\gamma'(s),\gamma'(s)) \, ds$, it follows that 
\[(1-\delta) \, \mathcal{R}^{\frac{1}{2}} \leq -d'(t) \leq (1+\frac{\delta}{2}+4A^{-1}) \, \mathcal{R}^{\frac{1}{2}}\] 
if $-t$ is sufficiently large (depending on $\delta$ and $A$). Since $4A^{-1} \leq \frac{\delta}{2}$, the assertion follows. This completes the proof of Lemma \ref{derivative.of.d(t)}. \\

In the next step, we state a consequence of Hamilton's Harnack inequality. In the following, we view the scalar curvature $R$ as a function of $r$ and $t$. We denote by $R_t$ the partial derivative of $R$ with respect to $t$ (keeping $r$ fixed). 

\begin{proposition}
\label{harnack}
We have $R_t - \frac{2}{r} \, u^{-1} u_t v \geq 0$.
\end{proposition}

\textbf{Proof.} 
Hamilton's trace Harnack inequality \cite{Hamilton2} implies 
\[R_t - R_r v + 2 R_r w + 2 \, \text{\rm Ric}_{rr} \, w^2 \, \geq 0,\] 
where $v = \frac{1}{r} \, (1-u - \frac{1}{2} \, ru_r)$ and $w$ is arbitrary. The extra term $-R_r v$ arises because we compute the time derivative of the scalar curvature at a fixed radius $r$, whereas Hamilton computes the time derivative at a fixed point on the manifold. Indeed, if we fix a point on the manifold, then the radius $r$ shrinks at a rate given by $-v$, and the scalar curvature changes at a rate of $R_t-R_r v$. 

Applying the Harnack inequality with $w := v$ gives 
\[R_t + R_r v + 2 \, \text{\rm Ric}_{rr} \, v^2 \geq 0.\] 
Note that 
\[R_r = -\frac{4}{r^3} \, (1-u+\frac{1}{2} \, r^2 u_{rr}) = -\frac{2}{r} \, u^{-1} u_t + \frac{2}{r} \, u^{-1} u_r v\] 
and 
\[\text{\rm Ric}_{rr} = -\frac{1}{r} \, u^{-1} u_r,\] 
hence 
\[R_r + 2 \, \text{\rm Ric}_{rr} \, v = -\frac{2}{r} \, u^{-1} u_t.\] 
Putting these facts together, the assertion follows. \\

We next consider the quantity $R+|V|^2 = R+u^{-1} v^2$. Note that this function is smooth across the tip. 

\begin{remark}
On the Bryant soliton, the function $u_t$ vanishes identically, and the function $R+u^{-1} v^2$ is equal to $1$.
\end{remark}

\begin{proposition} 
\label{monotonicity.of.R+u^{-1}v^2.in.spacetime}
We have 
\[(R+u^{-1} v^2)_t + \frac{v}{2} \, \Big ( 1 + \frac{r}{2} \, u^{-1} v \Big )^{-1} \, (R+u^{-1} v^2)_r \geq 0.\] 
\end{proposition}

\textbf{Proof.} 
We observe that 
\[v = \frac{r}{4} \, R + \frac{1}{2r} \, (1-u),\] 
hence 
\[v_t = \frac{r}{4} \, R_t - \frac{1}{2r} \, u_t.\] 
This gives  
\begin{align*} 
(R+u^{-1} v^2)_t 
&= R_t + 2u^{-1} vv_t - u^{-2} u_t v^2 \\ 
&= \Big ( 1 + \frac{r}{2} \, u^{-1} v \Big ) \, R_t - \frac{1}{r} \, u^{-1} u_t v - u^{-2} u_t v^2 \\ 
&= \Big ( 1 + \frac{r}{2} \, u^{-1} v \Big ) \, \Big ( R_t - \frac{2}{r} \, u^{-1} u_t v \Big ) + \frac{1}{r} \, u^{-1} u_t v.
\end{align*} 
Moreover, using the relations 
\[R_r = -\frac{2}{r} \, u^{-1} u_t + \frac{2}{r} \, u^{-1} u_r v\] 
and 
\[u^{-2} u_t = -\frac{2}{r} \, u^{-1} u_r + u^{-2} u_r v - 2u^{-1} v_r,\] 
we obtain 
\begin{align*} 
(R+u^{-1} v^2)_r 
&= R_r - u^{-2} u_r v^2  + 2u^{-1} vv_r \\ 
&= -\frac{2}{r} \, u^{-1} u_t + \frac{2}{r} \, u^{-1} u_r v - u^{-2} u_r v^2 + 2u^{-1} vv_r \\ 
&= -\frac{2}{r} \, u^{-1} u_t - u^{-2}  u_t v \\ 
&= -\frac{2}{r} \, \Big (1+\frac{r}{2} \, u^{-1} v \Big ) \, u^{-1} u_t. 
\end{align*} 
Consequently, 
\begin{align*} 
&(R+u^{-1} v^2)_t + \frac{v}{2} \, \Big ( 1 + \frac{r}{2} \, u^{-1} v \Big )^{-1} \, (R+u^{-1} v^2)_r \\ 
&= \Big ( 1 + \frac{r}{2} \, u^{-1} v \Big ) \, \Big ( R_t - \frac{2}{r} \, u^{-1} u_t v \Big ), 
\end{align*} 
and the right hand side is nonnegative by Proposition \ref{harnack}. \\

\begin{proposition}
\label{pde.for.R+u^{-1}v^2}
The function $R+u^{-1}v^2$ satisfies 
\[(R+u^{-1}v^2)_t = u \, (R+u^{-1}v^2)_{rr} + \frac{2}{r} \, u \, (R+u^{-1}v^2)_r + \Xi(r,t) \, (R+u^{-1}v^2)_r,\] 
where 
\[\Xi := (1+u-\frac{1}{2} ru_r)^{-1} \, \Big [ \frac{1}{r} \, (1-\frac{1}{2} ru_r) \, (1-u-\frac{1}{2} ru_r) - u^3 \, \partial_r(u^{-2} (1+u-\frac{1}{2} ru_r)) \Big ].\] 
For each $t$, we have $\Xi(r,t) = O(r)$ near the tip.
\end{proposition}

\textbf{Proof.} 
Differentiating the identity 
\[(R+u^{-1} v^2)_r = -\frac{2}{r} \, \Big (1+\frac{r}{2} \, u^{-1} v \Big ) \, u^{-1} u_t = -\frac{1}{r} \, u^{-2} \, (1+u-\frac{1}{2} ru_r) \, u_t\] 
with respect to $r$ gives 
\begin{align*} 
(R+u^{-1} v^2)_{rr} 
&= -\frac{1}{r} \, u^{-2} \, (1+u-\frac{1}{2} ru_r) \, u_{tr} + \frac{1}{r^2} \, u^{-2} \, (1+u-\frac{1}{2} ru_r) \, u_t \\ 
&- \frac{1}{r} \, \partial_r(u^{-2} (1+u-\frac{1}{2} ru_r)) \, u_t. 
\end{align*} 
On the other hand, differentiating the identity 
\[R+u^{-1} v^2 = \frac{1}{r^2} \, u^{-1} \, (1+u-\frac{1}{2} ru_r)^2 - \frac{2}{r^2} \, (1+u)\] 
with respect to $t$ gives 
\begin{align*} 
(R+u^{-1} v^2)_t 
&= -\frac{1}{r} \, u^{-1} \, (1+u-\frac{1}{2} ru_r) \, u_{rt} - \frac{2}{r^2} \, u_t \\ 
&- \frac{1}{r^2} \, u^{-2} \, (1+u-\frac{1}{2} ru_r) (1-u-\frac{1}{2}ru_r) \, u_t.
\end{align*}
Putting these facts together, we obtain 
\begin{align*} 
&(R+u^{-1} v^2)_t - u \, (R+u^{-1}v^2)_{rr} - \frac{2}{r} \, u \, (R+u^{-1}v^2)_r \\ 
&= -\Big [ \frac{1}{r} \, (1-\frac{1}{2} ru_r) (1-u-\frac{1}{2}ru_r) - u^3 \, \partial_r(u^{-2} (1+u-\frac{1}{2} ru_r)) \Big ] \\ 
&\hspace{10mm} \cdot \frac{1}{r} \, u^{-2} \, u_t \\
&= \Big [ \frac{1}{r} \, (1-\frac{1}{2} ru_r) (1-u-\frac{1}{2}ru_r) - u^3 \, \partial_r(u^{-2} (1+u-\frac{1}{2} ru_r)) \Big ] \\ 
&\hspace{10mm} \cdot (1+u-\frac{1}{2} ru_r)^{-1} \,  (R+u^{-1} v^2)_r,
\end{align*} 
as claimed. \\

\begin{corollary} 
\label{lower.bound.for.R+u^{1}v^2}
We have $R+u^{-1} v^2 \geq \mathcal{R}$ at each point in spacetime.
\end{corollary}

\textbf{Proof.} 
Let us fix a point $(r_0,t_0)$ in space-time such that $r_0 \in [0,r_{\text{\rm max}}(t_0))$. Let $\hat{r}(t)$ denote the solution of the ODE 
\[\frac{d}{dt} \hat{r}(t) = \frac{v(\hat{r}(t),t)}{2} \, \Big ( 1 + \frac{\hat{r}(t)}{2} \, u(\hat{r}(t),t)^{-1} v(\hat{r}(t),t) \Big )^{-1}\] 
with initial condition $\hat{r}(t_0)=r_0$. Since $v$ is a nonnegative function, we obtain $\hat{r}(t) \leq r_0$ for $t \leq t_0$. Consequently, the function $r \mapsto \hat{r}(t)$ is defined for all $t \in (-\infty,t_0]$, and $\hat{r}(t) \in [0,r_{\text{\rm max}}(t))$ for all $t \leq t_0$.

By Proposition \ref{monotonicity.of.R+u^{-1}v^2.in.spacetime}, the function $t \mapsto R(\hat{r}(t),t)+u(\hat{r}(t),t)^{-1} v(\hat{r}(t),t)^2$ is monotone increasing. On the other hand, by Proposition \ref{convergence.to.Bryant.soliton}, we can find a sequence of times $t_k \to -\infty$ such that the rescaled manifolds $(M,\mathcal{R} \, g(t_k))$ converge to the Bryant soliton in the Cheeger-Gromov sense. 
Since $R+u^{-1} v^2 = 1$ on the Bryant soliton, we conclude that 
\[\lim_{k \to \infty} \sup_{r \in (0,r_0]} |R(r,t_k) + u(r,t_k)^{-1} v(r,t_k)^2 - \mathcal{R}| = 0.\] 
Consequently, 
\begin{align*} 
\mathcal{R} 
&= \lim_{k \to \infty} R(\hat{r}(t_k),t_k)+u(\hat{r}(t_k),t_k)^{-1} v(\hat{r}(t_k),t_k)^2 \\ 
&\leq R(r_0,t_0)+u(r_0,t_0)^{-1} v(r_0,t_0)^2. 
\end{align*}
This completes the proof of Corollary \ref{lower.bound.for.R+u^{1}v^2}. \\

\begin{corollary} 
\label{monotonicity.of.R+u^{-1}v^2.in.space}
We have $(R+u^{-1}v^2)_r \geq 0$ at each point in spacetime. 
\end{corollary}

\textbf{Proof.} 
Let us fix a point $(r_0,t_0)$ in space-time such that $r_0 \in [0,r_{\text{\rm max}}(t_0))$. Let $\hat{r}(t)$ denote the solution of the ODE 
\[\frac{d}{dt} \hat{r}(t) = \frac{v(\hat{r}(t),t)}{2} \, \Big ( 1 + \frac{\hat{r}(t)}{2} \, u(\hat{r}(t),t)^{-1} v(\hat{r}(t),t) \Big )^{-1}\] 
with initial condition $\hat{r}(t_0)=r_0$. Clearly, $\hat{r}(t) \leq r_0$ for $t \leq t_0$. Consequently, the function $r \mapsto \hat{r}(t)$ is defined for all $t \in (-\infty,t_0]$, and $\hat{r}(t) \in [0,r_{\text{\rm max}}(t))$ for all $t \leq t_0$.

Let us consider an arbitrary sequence of times $t_k \to -\infty$. For $k$ large, we define $Q_k = \{(r,t): t_k \leq t \leq t_0, \, r \leq \hat{r}(t)\}$. By Proposition \ref{pde.for.R+u^{-1}v^2}, the function $R+u^{-1} v^2$ attains its maximum on the parabolic boundary of $Q_k$. Therefore, 
\begin{align*} 
&\sup_{r \leq r_0} R(r,t_0)+u(r,t_0)^{-1} v(r,t_0)^2 \\ 
&\leq \max \Big \{ \sup_{t_k \leq t \leq t_0} R(\hat{r}(t),t)+u(\hat{r}(t),t)^{-1} v(\hat{r}(t),t)^2, \\ 
&\hspace{20mm} \sup_{r \leq \hat{r}(t_k)} R(r,t_k)+u(r,t_k)^{-1}v(r,t_k)^2 \Big \} 
\end{align*} 
for $k$ large. By Proposition \ref{monotonicity.of.R+u^{-1}v^2.in.spacetime}, the function $t \mapsto R(\hat{r}(t),t)+u(\hat{r}(t),t)^{-1} v(\hat{r}(t),t)^2$ is monotone increasing. This implies 
\[\sup_{t_k \leq t \leq t_0} R(\hat{r}(t),t)+u(\hat{r}(t),t)^{-1} v(\hat{r}(t),t)^2 \leq R(r_0,t_0)+u(r_0,t_0)^{-1} v(r_0,t_0)^2\] 
for $k$ large. This gives
\begin{align*} 
&\sup_{r \leq r_0} R(r,t_0)+u(r,t_0)^{-1} v(r,t_0)^2 \\ 
&\leq \max \Big \{ R(r_0,t_0)+u(r_0,t_0)^{-1} v(r_0,t_0)^2, \\ 
&\hspace{20mm} \sup_{r \leq \hat{r}(t_k)} R(r,t_k)+u(r,t_k)^{-1}v(r,t_k)^2 \Big \} 
\end{align*} 
for $k$ large. We now send $k \to \infty$. Recall that $\hat{r}(t_k) \leq r_0$ for $k$ large. Since the solution looks like the Bryant soliton near the tip, we obtain  
\[\lim_{k \to \infty} \sup_{r \leq \hat{r}(t_k)} |R(r,t_k)+u(r,t_k)^{-1}v(r,t_k)^2 - \mathcal{R}| = 0.\] 
This gives 
\[\sup_{r \leq r_0} R(r,t_0)+u(r,t_0)^{-1} v(r,t_0)^2 \leq \max \{R(r_0,t_0)+u(r_0,t_0)^{-1} v(r_0,t_0)^2,\mathcal{R}\}.\] 
Since $R(r_0,t_0)+u(r_0,t_0)^{-1} v(r_0,t_0)^2 \geq \mathcal{R}$ by Corollary \ref{lower.bound.for.R+u^{1}v^2}, we conclude that 
\[\sup_{r \leq r_0} R(r,t_0)+u(r,t_0)^{-1} v(r,t_0)^2 \leq R(r_0,t_0)+u(r_0,t_0)^{-1} v(r_0,t_0)^2,\] 
which implies the claim. \\

\begin{lemma}
\label{points.of.large.radius.lie.on.necks}
Given $\varepsilon_0>0$, there exists a large constant $C_0$ with the following property. If $r \geq C_0$ at some point in space-time, then that point lies at the center of an $\varepsilon_0$-neck.
\end{lemma}

\textbf{Proof.} 
By work of Perelman \cite{Perelman1}, the set of all points in $(M,g(t))$ which do not lie at the center of an $\varepsilon_0$-neck has diameter less than $C(\varepsilon_0) \, R_{\text{\rm max}}(t)^{-\frac{1}{2}}$ (see Theorem \ref{canonical.neighborhood.theorem} and Corollary \ref{bound.for.R_max}). Hence, if $r > C(\varepsilon_0) \, R_{\text{\rm max}}(t)^{-\frac{1}{2}}$ at some point in spacetime, then that point lies at the center of an $\varepsilon_0$-neck. On the other hand, $R_{\text{\rm max}}(t)$ is uniformly bounded from below by Proposition \ref{lower.bounds.for.d(t).and.R_max(t)}. From this, the assertion follows. \\

\begin{lemma}
\label{bound.for.u}
On an $\varepsilon_0$-neck, we have $r^2u \leq (1+100\varepsilon_0) \mathcal{R}^{-1}$.  
\end{lemma}

\textbf{Proof.} 
On an $\varepsilon_0$-neck, we have $u \leq \varepsilon_0$. Moreover, on an $\varepsilon_0$-neck, the radial Ricci curvature is smaller than $\frac{10\varepsilon_0}{r^2}$. This gives $0 \leq -ru_r \leq 10\varepsilon_0$. Using Corollary \ref{lower.bound.for.R+u^{1}v^2}, we obtain 
\begin{align*} 
\mathcal{R} 
&\leq R+u^{-1}v^2 \\ 
&= \frac{1}{r^2} \, u^{-1} \, (1+u-\frac{1}{2} ru_r)^2 - \frac{2}{r^2} \, (1+u) \\ 
&\leq \frac{1}{r^2} \, u^{-1} \, (1+100\varepsilon_0). 
\end{align*} 
This proves the assertion. \\

\begin{lemma} 
\label{bounds.for.derivatives.of.F}
There exists a large constant $C_1$ with the following property. If $F \geq C_0$, then we have $F \, |F_z| \leq C_1$ and $F^2 \, |F_{zz}|+F^3 \, |F_{zzz}| \leq C_1 \, F^{\frac{1}{100}}$.
\end{lemma} 

\textbf{Proof.} 
By Lemma \ref{points.of.large.radius.lie.on.necks}, every point with $F \geq C_0$ lies at the center of an $\varepsilon_0$-neck. Using Lemma \ref{bound.for.u}, we obtain $F^2 F_z^2 \leq (1+100\varepsilon_0) \mathcal{R}^{-1}$ on an $\varepsilon_0$-neck. We next observe that $F^m \, |\partial_z^{m+1} F| \leq C(m)$ on an $\varepsilon_0$-neck. Using standard interpolation inequalities, we obtain $F^2 \, |F_{zz}|+F^3 \, |F_{zzz}| \leq C \, F^{\frac{1}{100}}$ whenever $F \geq C_0$. \\

\begin{lemma}
\label{bound.for.F_zz}
There exist large constants $C_2 \geq 4C_0$ and $C_3$ with the following property. If $F \geq C_2$, then $0 \leq -F_{zz} \leq C_3 F^{-\frac{5}{2}+\frac{1}{100}}$.
\end{lemma}

\textbf{Proof.} 
Let us fix a point $(z_0,t_0)$ in spacetime, and let $r_0=F(z_0,t_0) \in [0,r_{\text{\rm max}}(t_0))$. We assume that $r_0 \geq \max \{10C_0,100C_1^2\}$. Let $\tilde{r}(t)$ denote the solution of the ODE
\[\frac{d}{dt} \tilde{r}(t) = -v(\tilde{r}(t),t) = -\frac{1}{\tilde{r}(t)} \, (1-u(\tilde{r}(t),t)-\frac{1}{2} \, \tilde{r}(t) \, u_r(\tilde{r}(t),t))\] 
with initial condition $\tilde{r}(t_0)=r_0$. Note that $\tilde{r}(t)$ can be interpreted as the radius, at time $t$, of a sphere of symmetry passing through a fixed point on the manifold. In particular, $\tilde{r}(t) \in [0,r_{\text{\rm max}}(t))$ for $t \leq t_0$.

We define a function $\tilde{F}(z,t)$ by 
\[\tilde{F} \bigg ( \int_{\tilde{r}(t)}^\rho u(r,t)^{-\frac{1}{2}} \, dr,t \bigg ) = \rho.\] 
Clearly, $\tilde{F}(0,t) = \tilde{r}(t)$. Note that $F(z,t)$ and $\tilde{F}(z,t)$ differ only by a translation in $z$: 
\[\tilde{F}(z,t) = F \bigg ( z + \int_{\bar{r}(t)}^{\tilde{r}(t)} u(r,t)^{-\frac{1}{2}} \, dr,t \bigg ).\] 
Since $\tilde{F}(0,t_0)=F(z_0,t_0)=r_0$, we obtain $\tilde{F}(z,t_0) = F(z+z_0,t_0)$ for all $z$.

It follows from Lemma \ref{behavior.at.tip} that $-\frac{d}{dt} \tilde{r}(t) \geq 0$ for each $t$. Integrating this inequality over $t$ gives $\tilde{r}(t) \geq r_0$ for all $t \leq t_0$. Equivalently, $\tilde{F}(0,t) \geq r_0$ for all $t \leq t_0$. Moreover, Lemma \ref{bounds.for.derivatives.of.F} implies $\tilde{F} \, |\tilde{F}_z| \leq C_1$ whenever $\tilde{F} \geq C_0$. Hence, if $r_0 \geq \max \{10C_0,100C_1^2\}$, then we obtain 
\[\tilde{F}(z,t) \geq \sqrt{\tilde{F}(0,t)^2 - 2C_1 |z|} \geq \sqrt{r_0^2 - 2C_1 r_0^{\frac{3}{2}}} \geq \frac{1}{2} \, r_0\] 
for all $t \leq t_0$ and all $z \in [-r_0^{\frac{3}{2}},r_0^{\frac{3}{2}}]$. 

The function $\tilde{F}$ satisfies the evolution equation 
\begin{align*} 
0 
&= \tilde{F}_t(z,t) - \tilde{F}_{zz}(z,t) + \tilde{F}(z,t)^{-1} \, (1+\tilde{F}_z(z,t)^2) \\ 
&+ 2 \, \tilde{F}_z(z,t) \, \bigg [ -\tilde{F}(0,t)^{-1} \, \tilde{F}_z(0,t) + \int_{\tilde{F}(0,t)}^{\tilde{F}(z,t)} \frac{1}{r^2} \, u(r,t)^{\frac{1}{2}} \, dr \bigg ]. 
\end{align*}
This implies 
\begin{align*} 
0 
&= (\tilde{F}^2)_t(z,t) - (\tilde{F}^2)_{zz}(z,t) + 2+4\tilde{F}_z(z,t)^2 \\ 
&+ 4\tilde{F}(z,t) \tilde{F}_z(z,t) \, \bigg [ -\tilde{F}(0,t)^{-1} \, \tilde{F}_z(0,t) + \int_{\tilde{F}(0,t)}^{\tilde{F}(z,t)} \frac{1}{r^2} \, u(r,t)^{\frac{1}{2}} \, dr \bigg ]. 
\end{align*} 
Consequently, if we define $\tilde{Q}(z,t) := \frac{1}{2} \, (\tilde{F}^2)_z(z,t) = \tilde{F}(z,t)\tilde{F}_z(z,t)$, then 
\begin{align*} 
0 
&= \tilde{Q}_t(z,t) - \tilde{Q}_{zz}(z,t) + 4\tilde{F}_z(z,t) \, \tilde{F}_{zz}(z,t) + 2\tilde{F}(z,t)^{-1} \tilde{F}_z(z,t)^3 \\ 
&+ 2 \, (\tilde{F}(z,t) \tilde{F}_{zz}(z,t)+\tilde{F}_z(z,t)^2) \, \bigg [ -\tilde{F}(0,t)^{-1} \, \tilde{F}_z(0,t) + \int_{\tilde{F}(0,t)}^{\tilde{F}(z,t)} \frac{1}{r^2} \, u(r,t)^{\frac{1}{2}} \, dr \bigg ]. 
\end{align*} 
Using Lemma \ref{bounds.for.derivatives.of.F}, we obtain $|\tilde{Q}(z,t)| \leq C$ and $|\tilde{Q}_t(z,t) - \tilde{Q}_{zz}(z,t)| \leq C r_0^{-3+\frac{1}{100}}$ for $t \in [t_0-r_0^3,t_0]$ and $z \in [-r_0^{\frac{3}{2}},r_0^{\frac{3}{2}}]$. Using standard interior estimates for linear parabolic equations, we conclude that $|\tilde{Q}_z(0,t_0)| \leq C r_0^{-\frac{3}{2}+\frac{1}{100}}$. Consequently, $|\tilde{F}_{zz}(0,t_0)| \leq C r_0^{-\frac{5}{2}+\frac{1}{100}}$. This finally implies $|F_{zz}(z_0,t_0)| \leq C r_0^{-\frac{5}{2}+\frac{1}{100}}$. This completes the proof of Lemma \ref{bound.for.F_zz}. \\

\begin{lemma}
\label{control.of.F_t}
There exist large constants $C_4$ and $C_5$ with the following property. If $-t \geq C_4$ and $F(z,t) \geq C_4$, then 
\[|F(z,t) F_t(z,t) + 1| \leq C_5 \, F(z,t)^{-\frac{3}{2}+\frac{1}{100}} + C_5 \, (-t)^{-1}.\]
\end{lemma} 

\textbf{Proof.}
Recall that 
\begin{align*} 
0 
&= F_t(z,t) - F_{zz}(z,t) + F(z,t)^{-1} \, (1+F_z(z,t)^2) \\ 
&+ 2 \, F_z(z,t) \, \bigg [ -F(0,t)^{-1} \, F_z(0,t) + \int_{F(0,t)}^{F(z,t)} \frac{1}{r^2} \, u(r,t)^{\frac{1}{2}} \, dr \bigg ]. 
\end{align*}
Note that $F(0,t)^{-1} \, |F_z(0,t)| \leq C \, (-t)^{-1}$. Moreover, if $F(z,t) \geq \max \{C_0,C_2\}$, then Lemma \ref{bounds.for.derivatives.of.F} and Lemma \ref{bound.for.F_zz} imply $|F_z(z,t)| \leq C \, F(z,t)^{-1}$ and $|F_{zz}(z,t)| \leq C \, F(z,t)^{-\frac{5}{2}+\frac{1}{100}}$. Putting these facts together, the assertion follows. \\

\begin{lemma}
\label{lower.bound.for.F_t}
There exists a large constant $C_6$ with the following property. If $-t \geq C_6$ and $F(z,t) \geq C_6$, then $-F(z,t) F_t(z,t) \geq \frac{1}{2}$.
\end{lemma}

\textbf{Proof.} 
This follows immediately from Lemma \ref{control.of.F_t}. \\

In the following, $p$ will denote the tip. By Proposition \ref{convergence.to.Bryant.soliton}, we can find a large constant $C_7$ with the following property. If $-t \geq C_7$ and $x$ is a point in $(M,g(t))$ with $d_{g(t)}(p,x) \geq C_7$, then the sphere of symmetry passing through $x$ has radius $r > \max \{C_0,C_2,C_4,C_6\}$. Moreover, let us fix a large constant $\Lambda \geq C_7$ such that $d(t) = d_{g(t)}(p,q) < \Lambda$ for all $t \in [-\max\{C_4,C_6,C_7\},0]$.

For each $z \in (-\infty,0]$, we define a time $\mathcal{T}(z) \in (-\infty,-\max\{C_4,C_6,C_7\}]$ by 
\[t = \mathcal{T}(z) \iff d(t)=\Lambda-z.\] 
In other words, at time $\mathcal{T}(z)$, the reference point $q$ has distance $\Lambda-z$ from the tip. 

\begin{lemma} 
\label{radius.for.t.less.than.T(z)}
Let $z \leq 0$ and $t \leq \mathcal{T}(z)$. Then $F(z,t) > \max \{C_0,C_2,C_4,C_6\}$.
\end{lemma}

\textbf{Proof.} 
By assumption, $d(t)=d_{g(t)}(p,q) \geq \Lambda-z$. Let $x$ be a point in $(M,g(t))$ which has signed distance $z$ from the reference point $q$. Then $d_{g(t)}(p,x) \geq \Lambda$. In particular, $d_{g(t)}(p,x) \geq C_7$. Moreover, $-t \geq -\mathcal{T}(z) \geq C_7$. By our choice of $C_7$, the sphere of symmetry passing through $x$ has radius greater than $\max \{C_0,C_2,C_4,C_6\}$. This completes the proof of Lemma \ref{radius.for.t.less.than.T(z)}. \\

\begin{lemma} 
\label{radius.for.t.equal.to.T(z)}
Let $z \leq 0$ and $t = \mathcal{T}(z)$. Then $F(z,t) \leq \Lambda$. 
\end{lemma}

\textbf{Proof.}
By assumption, $d(t)=d_{g(t)}(p,q) = \Lambda-z$. Let $x$ be a point in $(M,g(t))$ which has signed distance $z$ from the reference point $q$. Then $d_{g(t)}(p,x) = \Lambda$. Hence, the sphere of symmetry passing through $x$ has radius at most $\Lambda$. This completes the proof of Lemma \ref{radius.for.t.equal.to.T(z)}. \\

\begin{lemma} 
\label{relate.radius.to.T(z)-t}
There exists a large constant $C_8$ such that 
\[|F(z,t)^2 - 2 \, (\mathcal{T}(z)-t)| \leq C_8 \, (\mathcal{T}(z)-t)^{\frac{1}{4}+\frac{1}{200}} + C_8\] 
whenever $z \leq 0$ and $t \leq \mathcal{T}(z)$.
\end{lemma}

\textbf{Proof.} 
Using Lemma \ref{lower.bound.for.F_t} and Lemma \ref{radius.for.t.less.than.T(z)}, we obtain $-F(z,t) F_t(z,t) \geq \frac{1}{2}$ whenever $z \leq 0$ and $t \leq \mathcal{T}(z)$. Integrating this inequality over $t$ gives $F(z,t)^2 \geq \mathcal{T}(z)-t$ whenever $z \leq 0$ and $t \leq \mathcal{T}(z)$. 

Using Lemma \ref{control.of.F_t} and Lemma \ref{radius.for.t.less.than.T(z)}, we obtain 
\[|F(z,t) F_t(z,t) + 1| \leq C \, F(z,t)^{-\frac{3}{2}+\frac{1}{100}} + C \, (-t)^{-1}\] 
whenever $z \leq 0$ and $t \leq \mathcal{T}(z)$. Using the inequality $F(z,t)^2 \geq \mathcal{T}(z)-t$, we deduce that 
\[|F(z,t) F_t(z,t) + 1| \leq C \, (\mathcal{T}(z)-t)^{-\frac{3}{4}+\frac{1}{200}}\] 
whenever $z \leq 0$ and $t \leq \mathcal{T}(z)$. Integrating this inequality over $t$ gives 
\[|F(z,t)^2 - 2 \, (\mathcal{T}(z)-t)| \leq C \, (\mathcal{T}(z)-t)^{\frac{1}{4}+\frac{1}{200}} + C\] 
whenever $z \leq 0$ and $t \leq \mathcal{T}(z)$. In the last step, we have used the fact that, by Lemma \ref{radius.for.t.equal.to.T(z)}, $F(z,t) \leq \Lambda$ whenever $z \leq 0$ and $t = \mathcal{T}(z)$. This completes the proof of Lemma \ref{relate.radius.to.T(z)-t}. \\

\begin{lemma}
\label{lower.bound.for.FF_z.at.waist}
Let $\delta>0$ be given. Then $F(0,t) F_z(0,t) \geq (1+4\delta)^{-1} \, \mathcal{R}^{-\frac{1}{2}}$ if $-t$ is sufficiently large.
\end{lemma}

\textbf{Proof.} 
By Lemma \ref{derivative.of.d(t)}, $(1-\delta) \, \mathcal{R}^{\frac{1}{2}} \leq -d'(t) \leq (1+\delta) \, \mathcal{R}^{\frac{1}{2}}$ if $-t$ is sufficiently large (depending on $\delta$). Integrating over $t$, we obtain 
\[(1-2\delta) \, \mathcal{R}^{\frac{1}{2}} \, (-t) \leq d(t) \leq (1+2\delta) \, \mathcal{R}^{\frac{1}{2}} \, (-t)\] 
if $-t$ is sufficiently large (depending on $\delta$). Putting $t = \mathcal{T}(z)$ gives 
\[(1+2\delta)^{-1} \, \mathcal{R}^{-\frac{1}{2}} \, (\Lambda-z) \leq -\mathcal{T}(z) \leq (1-2\delta)^{-1} \, \mathcal{R}^{-\frac{1}{2}} \, (\Lambda-z)\] 
if $-z$ is sufficiently large (depending on $\delta$). 

In the following, we assume that $-t$ is sufficiently large, so that $t \leq \mathcal{T}(-\sqrt{-t}) \leq \mathcal{T}(0)$. We apply Lemma \ref{relate.radius.to.T(z)-t} with $z = 0$ and, separately, with $z = -\sqrt{-t}$. This gives 
\[|F(0,t)^2 - 2 \, (\mathcal{T}(0)-t)| \leq C_8 \, (\mathcal{T}(0)-t)^{\frac{1}{4}+\frac{1}{200}} + C_8.\] 
and 
\[|F(-\sqrt{-t},t)^2 - 2 \, (\mathcal{T}(-\sqrt{-t})-t)| \leq C_8 \, (\mathcal{T}(-\sqrt{-t})-t)^{\frac{1}{4}+\frac{1}{200}} + C_8.\] 
This implies 
\[F(0,t)^2-F(-\sqrt{-t},t)^2 - 2 \, (\mathcal{T}(0)-\mathcal{T}(-\sqrt{-t})) \geq -2C_8 \, (-t)^{\frac{1}{4}+\frac{1}{200}} - 2C_8.\] 
Moreover, 
\[\mathcal{T}(0) - \mathcal{T}(-\sqrt{-t}) \geq (1+2\delta)^{-1} \, \mathcal{R}^{-\frac{1}{2}} \, (\Lambda+\sqrt{-t}) + \mathcal{T}(0)\] 
if $-t$ is sufficiently large (depending on $\delta$). Putting these facts together, we obtain 
\begin{align*} 
F(0,t)^2-F(-\sqrt{-t},t)^2 
&\geq 2 \, (1+2\delta)^{-1} \, \mathcal{R}^{-\frac{1}{2}} \, (\Lambda+\sqrt{-t}) + 2 \, \mathcal{T}(0) \\ 
&- 2C_8 \, (-t)^{\frac{1}{4}+\frac{1}{200}} - 2C_8 
\end{align*} 
if $-t$ is sufficiently large (depending on $\delta$). Consequently, 
\[\sup_{z \in [-\sqrt{-t},0]} F(z,t) F_z(z,t) \geq \frac{1}{2\sqrt{-t}} \, (F(0,t)^2-F(-\sqrt{-t},t)^2) \geq (1+3\delta)^{-1} \, \mathcal{R}^{-\frac{1}{2}}\] 
if $-t$ is sufficiently large (depending on $\delta$). On the other hand, using Lemma \ref{bounds.for.derivatives.of.F} and Lemma \ref{bound.for.F_zz}, we obtain $|(FF_z)_z| = |FF_{zz}+F_z^2| \leq C F^{-\frac{3}{2}+\frac{1}{100}} \leq C \, (-t)^{-\frac{3}{4}+\frac{1}{200}}$ for all $z \in [-\sqrt{-t},0]$. This implies 
\[\sup_{z \in [-\sqrt{-t},0]} F(z,t) F_z(z,t) \leq F(0,t) F_z(0,t) + C \, (-t)^{-\frac{1}{4}+\frac{1}{200}}.\] 
Thus, we conclude that 
\[F(0,t) F_z(0,t) \geq (1+4\delta)^{-1} \, \mathcal{R}^{-\frac{1}{2}}\] 
if $-t$ is sufficiently large (depending on $\delta$). This completes the proof of Lemma \ref{lower.bound.for.FF_z.at.waist}. \\

The following lemma is similar to Proposition 6.10 in \cite{Brendle-Choi}:

\begin{lemma}
\label{lower.bound.for.FF_z}
Let $\delta>0$ be given. Then 
\[\inf_{z \geq 0} F(z,t) F_z(z,t) \geq (1+5\delta)^{-1} \, \mathcal{R}^{-\frac{1}{2}}\] 
if $-t$ is sufficiently large.
\end{lemma}

\textbf{Proof.} 
Let us define $Q(z,t) := \frac{1}{2} \, (F^2)_z(z,t) = F(z,t) F_z(z,t)$. Note that $Q(z,t) \geq 0$. Moreover, 
\begin{align*} 
0 
&= Q_t(z,t) - Q_{zz}(z,t) + 4F_z(z,t) \, F_{zz}(z,t) + 2F(z,t)^{-1} F_z(z,t)^3 \\ 
&+ 2 \, (F(z,t) F_{zz}(z,t)+F_z(z,t)^2) \, \bigg [ -F(0,t)^{-1} \, F_z(0,t) + \int_{F(0,t)}^{F(z,t)} \frac{1}{r^2} \, u(r,t)^{\frac{1}{2}} \, dr \bigg ]. 
\end{align*} 
If $-t$ is sufficiently large, then we have $F(z,t) \geq F(0,t) \geq \sqrt{-t}$ for all $z \geq 0$. Using Lemma \ref{bounds.for.derivatives.of.F}, we obtain 
\begin{align*} 
|Q_t(z,t) - Q_{zz}(z,t)| 
&\leq C \, F(z,t)^{-3+\frac{1}{100}} + C \, (-t)^{-1} \, F(z,t)^{-1+\frac{1}{100}} \\ 
&\leq C \, (-t)^{-\frac{3}{2}+\frac{1}{200}} 
\end{align*}
for all $z \geq 0$. Let $\hat{Q}(z,t) := Q(z,t) + (-t)^{-\frac{1}{2}+\frac{1}{100}}$. Clearly, $\inf_{z \geq 0} \hat{Q}(z,t) > 0$ for each $t$. Moreover, if $-t$ is sufficiently large, then 
\[\hat{Q}_t(z,t) - \hat{Q}_{zz}(z,t) > 0\] 
for all $z \geq 0$. Finally, Lemma \ref{lower.bound.for.FF_z.at.waist} implies that $\hat{Q}(0,t) > Q(0,t) \geq (1+4\delta)^{-1} \, \mathcal{R}^{-\frac{1}{2}}$ if $-t$ is sufficiently large. Basic facts about the one-dimensional heat equation on the half-line with Dirichlet boundary condition imply that $\inf_{z \geq 0} \hat{Q}(z,t) \geq (1+4\delta)^{-1} \, \mathcal{R}^{-\frac{1}{2}}$ if $-t$ is sufficiently large. Consequently, $\inf_{z \geq 0} Q(z,t) \geq (1+5\delta)^{-1} \, \mathcal{R}^{-\frac{1}{2}}$ if $-t$ is sufficiently large. This completes the proof of Lemma \ref{lower.bound.for.FF_z}. \\

\begin{lemma}
\label{paraboloid.asymptotics}
If $-t$ is sufficiently large, then $F(z,t)^2 \geq \frac{1}{C_9} \, (z-t)$ and $|(FF_z)_t| \leq C_{10} \, F^{-1+\frac{1}{100}}$ for all $z \geq 0$.
\end{lemma} 

\textbf{Proof.} 
If $-t$ is sufficiently large, then Lemma \ref{lower.bound.for.FF_z} implies that $F(z,t) F_z(z,t) \geq \frac{1}{C}$ for all $z \geq 0$. Integrating this inequality over $z$ gives $F(z,t)^2 \geq \frac{1}{C} \, (z-t)$ for all $z \geq 0$. This proves the first statement. To prove the second statement, we consider the function $Q(z,t) := F(z,t) F_z(z,t)$. We have shown above that $|Q_t(z,t) - Q_{zz}(z,t)| \leq C \, F(z,t)^{-1+\frac{1}{100}}$ for all $z \geq 0$. Moreover, Lemma \ref{bounds.for.derivatives.of.F} implies $|Q_{zz}(z,t)| = |F(z,t) F_{zzz}(z,t) + F_z(z,t) F_{zz}(z,t)| \leq C \, F(z,t)^{-2+\frac{1}{100}}$ for all $z \geq 0$. Consequently, $|Q_t(z,t)| \leq C \, F(z,t)^{-1+\frac{1}{100}}$ for all $z \geq 0$. This completes the proof of Lemma \ref{paraboloid.asymptotics}. \\

\begin{lemma}
\label{limit.of.FF_z}
We have $\lim_{z \to \infty} F(z,t) F_z(z,t) = \mathcal{R}^{-\frac{1}{2}}$ if $-t$ is sufficiently large.
\end{lemma}

\textbf{Proof.} 
Lemma \ref{paraboloid.asymptotics} implies that $\liminf_{z \to \infty} F(z,t) F_z(z,t)$ is independent of $t$, provided that $-t$ is sufficiently large. On the other hand, it follows from Lemma \ref{lower.bound.for.FF_z} that $\liminf_{z \to \infty} F(z,t) F_z(z,t) \geq (1+5\delta)^{-1} \, \mathcal{R}^{-\frac{1}{2}}$ if $-t$ is sufficiently large (depending on $\delta$). Since $\delta$ is arbitrary, we conclude that $\liminf_{z \to \infty} F(z,t) F_z(z,t) \geq \mathcal{R}^{-\frac{1}{2}}$ if $-t$ is sufficiently large. On the other hand, since $(M,g(t))$ is neck-like at spatial infinity, Lemma \ref{bound.for.u} implies that $\limsup_{z \to \infty} F(z,t) F_z(z,t) \leq \mathcal{R}^{-\frac{1}{2}}$ for each $t$. This completes the proof of Lemma \ref{limit.of.FF_z}. \\

After these preparations, we now complete the proof of Theorem \ref{main.thm.a}. Lemma \ref{paraboloid.asymptotics} implies that $r_{\text{\rm max}}(t) = \infty$ if $-t$ is sufficiently large. By Lemma \ref{limit.of.FF_z}, we have $\lim_{z \to \infty} F(z,t) F_z(z,t) = \mathcal{R}^{-\frac{1}{2}}$ if $-t$ is sufficiently large. Equivalently, $\lim_{r \to \infty} r^2 u(r,t) = \mathcal{R}^{-1}$ if $-t$ is sufficiently large. Moreover, since $(M,g(t))$ is neck-like at spatial infinity, we know that $\lim_{r \to \infty} u(r,t) = 0$ and $\lim_{r \to \infty} r u_r(r,t) = 0$. Using the identity 
\[R+u^{-1}v^2 = \frac{1}{r^2} \, u^{-1} \, (1+u-\frac{1}{2} ru_r)^2 - \frac{2}{r^2} \, (1+u),\] 
we obtain $\lim_{r \to \infty} R+u^{-1} v^2 = \mathcal{R}$ if $-t$ is sufficiently large. Corollary \ref{monotonicity.of.R+u^{-1}v^2.in.space} then implies $R+u^{-1} v^2 \leq \mathcal{R}$ if $-t$ is sufficiently large. Using Corollary \ref{lower.bound.for.R+u^{1}v^2}, we conclude that $R+u^{-1} v^2=\mathcal{R}$ if $-t$ is sufficiently large. In view of the identity 
\[(R+u^{-1} v^2)_r = -\frac{2}{r} \, \Big (1+\frac{r}{2} \, u^{-1} v \Big ) \, u^{-1} u_t,\] 
it follows that $u_t=0$ if $-t$ is sufficiently large. Consequently, $(M,g(t))$ is a steady gradient Ricci soliton if $-t$ is sufficiently large. By the uniqueness result in \cite{Chen-Zhu}, $(M,g(t))$ is a steady gradient Ricci soliton for all $t$. This completes the proof of Theorem \ref{main.thm.a}.

\part{Proof of Theorem \ref{main.thm.b}}

\section{A PDE for the Lie derivative of the metric along a vector field}

\label{lie.derivatives}

We now study general solutions to the Ricci flow which are not necessarily rotationally symmetric. Given a Riemannian metric $g$ and a symmetric $(0,2)$-tensor $h$, we define the Lichnerowicz Laplacian of $h$ by 
\[\Delta_{L,g} h_{ik} = \Delta h_{ik} + 2 \, R_{ijkl} \, h^{jl} - \text{\rm Ric}_i^l \, h_{kl} - \text{\rm Ric}_k^l \, h_{il}.\] 
Moreover, the divergence of $h$ is defined by 
\[(\text{\rm div} \, h)^k = D_i h^{ik}.\] 
The following fact plays a key role in our analysis: 

\begin{proposition}
\label{prep}
Let $g$ be a Riemannian metric on a manifold $M$, and let $V$ be a vector field. We define $h := \mathscr{L}_V(g)$ and $Z := \text{\rm div} \, h - \frac{1}{2} \, \nabla(\text{\rm tr} \, h)$. Then 
\[Z = \Delta V + \text{\rm Ric}(V),\] 
where $\text{\rm Ric}$ is viewed as a $(1,1)$-tensor. Moreover, 
\[\mathscr{L}_V(\text{\rm Ric}) = -\frac{1}{2} \, \Delta_{L,g} h + \frac{1}{2} \, \mathscr{L}_Z(g),\] 
where $\text{\rm Ric}$ is viewed as a $(0,2)$-tensor.
\end{proposition}

\textbf{Proof.} 
Using the identity $h_{ij} = D_i V_j + D_j V_i$, we obtain 
\begin{align*} 
Z_k 
&= g^{ij} \, D_i h_{jk} - \frac{1}{2} \, g^{ij} \, D_k h_{ij} \\ 
&= g^{ij} \, D_{i,j}^2 V_k + g^{ij} \, D_{i,k}^2 V_j - g^{ij} \, D_{k,i}^2 V_j \\ 
&= \Delta V_k + \text{\rm Ric}_k^l \, V_l. 
\end{align*} 
This proves the first statement. 

To prove the second statement, let $\varphi_s: M \to M$ denote the one-parameter family of diffeomorphisms generated by $V$. Then $\frac{\partial}{\partial s} \varphi_s^*(g) \big |_{s=0} = h$. Using Proposition 2.3.7 in \cite{Topping}, we obtain 
\[\mathscr{L}_V(\text{\rm Ric}) = \frac{\partial}{\partial s} \text{\rm Ric}_{\varphi_s^*(g)} \Big |_{s=0} = -\frac{1}{2} \, \Delta_{L,g} h + \frac{1}{2} \, \mathscr{L}_Z(g).\] 
This completes the proof of Proposition \ref{prep}. \\

We now state the main result of this section:

\begin{corollary} 
\label{pde.for.lie.derivative}
Suppose that $g(t)$ is a solution to the Ricci flow on a manifold $M$. Moreover, suppose that $V(t)$ is a family of vector fields satisfying 
\[\frac{\partial}{\partial t} V(t) = \Delta_{g(t)} V(t) + \text{\rm Ric}_{g(t)}(V(t)).\] 
Then the Lie derivative $h(t) := \mathscr{L}_{V(t)}(g(t))$ satisfies the parabolic Lichnerowicz equation 
\[\frac{\partial}{\partial t} h(t) = \Delta_{L,g(t)} h(t).\] 
\end{corollary} 

\textbf{Proof.} 
As above, let $Z := \text{\rm div} \, h - \frac{1}{2} \, \nabla(\text{\rm tr} \, h)$. Proposition \ref{prep} implies that $\frac{\partial}{\partial t} V = \Delta V + \text{\rm Ric}(V) = Z$, where $\text{\rm Ric}$ is viewed as a $(1,1)$-tensor. Moreover, $\frac{\partial}{\partial t} g = -2 \, \text{\rm Ric}$, where $\text{\rm Ric}$ is viewed as a $(0,2)$-tensor. Using Proposition \ref{prep} again, we obtain 
\begin{align*} 
\frac{\partial}{\partial t} h
&= \mathscr{L}_V \Big ( \frac{\partial}{\partial t} g \Big ) + \mathscr{L}_{\frac{\partial}{\partial t} V}(g) \\ 
&= -2 \, \mathscr{L}_V(\text{\rm Ric}) + \mathscr{L}_Z(g) \\ 
&= \Delta_{L,g} h, 
\end{align*} 
where $\text{\rm Ric}$ is viewed as a $(0,2)$-tensor. This completes the proof of Corollary \ref{pde.for.lie.derivative}. \\

\begin{proposition} 
\label{estimate.for.norm.of.vector.field}
Suppose that $g(t)$ is a solution to the Ricci flow on a manifold $M$. Moreover, suppose that $V(t)$ is a family of vector fields satisfying \
\[\frac{\partial}{\partial t} V(t) = \Delta_{g(t)} V(t) + \text{\rm Ric}_{g(t)}(V(t)) + Q(t).\] 
Then 
\[\frac{\partial}{\partial t} |V(t)|_{g(t)} \leq \Delta_{g(t)} |V(t)|_{g(t)} + |Q(t)|_{g(t)}\] 
on the set $\{V(t) \neq 0\}$.
\end{proposition}

\textbf{Proof.} 
We compute 
\begin{align*} 
\frac{1}{2} \, \frac{\partial}{\partial t}(|V|^2) 
&= \langle V,\frac{\partial}{\partial t} V \rangle - \text{\rm Ric}(V,V) \\ 
&= \langle V,\Delta V \rangle + \langle V,Q \rangle \\ 
&= \frac{1}{2} \, \Delta(|V|^2) - |DV|^2 + \langle V,Q \rangle \\ 
&\leq \frac{1}{2} \, \Delta(|V|^2) - \big | \nabla |V| \big |^2 + |V| \, |Q|. 
\end{align*}
From this, the assertion follows easily. \\

\section{The parabolic Lichnerowicz equation on shrinking cylinders}

\label{model.problem}

In this section, we study the parabolic Lichnerowicz equation in a model case where the background metrics are a family of shrinking cylinders. Let $(S^2 \times \mathbb{R},\bar{g}(t))$ be a family of shrinking cylinders evolving by Ricci flow, so that $\bar{g}(t) = (-2t) \, g_{S^2} + dz \otimes dz$ for $t<0$. 

\begin{proposition} 
\label{parabolic.lichnerowicz.equation.on.cylinder}
Let $h(t)$ be a one-parameter family of symmetric $(0,2)$-tensors on the cylinder which is defined in the region $\{|z| \leq \frac{L}{2}, \, -\frac{L}{2} \leq t \leq -1\}$ and satisfies the parabolic Lichnerowicz equation $\frac{\partial}{\partial t} h(t) = \Delta_{L,\bar{g}(t)} h(t)$. Assume that $|h(t)|_{\bar{g}(t)} \leq 1$ in the region $\{|z| \leq \frac{L}{2}, \, -\frac{L}{2} \leq t \leq -\frac{L}{4}\}$, and $|h(t)|_{\bar{g}(t)} \leq L^{10}$ in the region $\{|z| \leq \frac{L}{2}, \, -\frac{L}{4} \leq t \leq -1\}$. On each slice $S^2 \times \{z\}$, we may decompose the tensor $h(t)$ as 
\[h(t) = \omega(z,t) \, g_{S^2} + \chi(z,t) + dz \otimes \sigma(z,t) + \sigma(z,t) \otimes dz + \beta(z,t) \, dz \otimes dz,\] 
where $\omega(z,t)$ is a scalar function on $S^2$, $\chi(z,t)$ is a tracefree symmetric $(0,2)$-tensor on $S^2$, $\sigma(z,t)$ is a one-form on $S^2$, and $\beta(z,t)$ is a scalar function on $S^2$. Then there exists a function $\psi: S^2 \to \mathbb{R}$ (independent of $t$ and $z$) such that $\psi$ lies in the span of the first spherical harmonics on $S^2$, and 
\[\big | h(t) - \bar{\omega}(z,t) \, g_{S^2} - \bar{\beta}(z,t) \, dz \otimes dz - (-t) \, \psi \, g_{S^2} \big |_{\bar{g}(t)} \leq C L^{-\frac{1}{2}}\] 
in the region $\{|z| \leq 1000, \, -1000 \leq t \leq -1\}$. Here, $\bar{\omega}(z,t)$ and $\bar{\beta}(z,t)$ are rotationally invariant functions satisfying 
\[\int_{S^2 \times \{z\}} (\omega(z,t)-\bar{\omega}(z,t)) \, d\text{\rm vol}_{S^2} = \int_{S^2 \times \{z\}} (\beta(z,t)-\bar{\beta}(z,t)) \, d\text{\rm vol}_{S^2} = 0\] 
for $t \in [-1000,-1]$ and $z \in [-1000,1000]$. In other words, $\bar{\omega}(z,t)$ and $\bar{\beta}(z,t)$ are obtained from $\omega(z,t)$ and $\beta(z,t)$ by averaging over the individual two-spheres $S^2 \times \{z\}$.
\end{proposition}

\textbf{Proof.} 
The parabolic Lichnerowicz equation is equivalent to the following system of equations for $\omega(z,t)$, $\chi(z,t)$, $\sigma(z,t)$, and $\beta(z,t)$: 
\begin{align*} 
&\frac{\partial}{\partial t} \omega(z,t) = \frac{\partial^2}{\partial z^2} \omega(z,t) + \frac{1}{(-2t)} \, \Delta_{S^2} \omega(z,t) \\ 
&\frac{\partial}{\partial t} \chi(z,t) = \frac{\partial^2}{\partial z^2} \chi(z,t) + \frac{1}{(-2t)} \, (\Delta_{S^2} \chi(z,t) - 4 \, \chi(z,t)), \\ 
&\frac{\partial}{\partial t} \sigma(z,t) = \frac{\partial^2}{\partial z^2} \sigma(z,t) + \frac{1}{(-2t)} \, (\Delta_{S^2} \sigma(z,t) - \sigma(z,t)), \\ 
&\frac{\partial}{\partial t} \beta(z,t) = \frac{\partial^2}{\partial z^2} \beta(z,t) + \frac{1}{(-2t)} \, \Delta_{S^2} \beta(z,t). 
\end{align*}
By assumption, $|h(t)|_{\bar{g}(t)} \leq 1$ in the region $\{|z| \leq \frac{L}{2}, \, -\frac{L}{2} \leq t \leq -\frac{L}{4}\}$, and $|h(t)|_{\bar{g}(t)} \leq L^{10}$ in the region $\{|z| \leq \frac{L}{2}, \, -\frac{L}{4} \leq t \leq -1\}$. This implies 
\begin{align*} 
&|\omega(z,t)| \leq C \, (-t), \\ 
&|\chi(z,t)|_{g_{S^2}} \leq C \, (-t), \\ 
&|\sigma(z,t)|_{g_{S^2}} \leq C \, (-t)^{\frac{1}{2}}, \\ 
&|\beta(z,t)| \leq C 
\end{align*}
in the region $\{|z| \leq \frac{L}{2}, \, -\frac{L}{2} \leq t \leq -\frac{L}{4}\}$, and 
\begin{align*} 
&|\omega(z,t)| \leq C L^{10} \, (-t), \\ 
&|\chi(z,t)|_{g_{S^2}} \leq C L^{10} \, (-t), \\ 
&|\sigma(z,t)|_{g_{S^2}} \leq C L^{10} \, (-t)^{\frac{1}{2}}, \\ 
&|\beta(z,t)| \leq C L^{10}
\end{align*} 
in the region $\{|z| \leq \frac{L}{2}, \, -\frac{L}{4} \leq t \leq -1\}$.

\textit{Step 1:} We first analyze the equation for $\chi(z,t)$. Let $S_j$, $j = 1,2,\hdots$, denote the eigenfunctions of the Laplacian on tracefree symmetric $(0,2)$-tensors on $S^2$, so that $\Delta_{S^2} S_j = -\nu_j S_j$. Clearly, $\nu_j > 0$ for each $j$. We assume that the eigenfunctions $S_j$ are normalized so that $\int_{S^2} |S_j|_{g_{S^2}}^2 \, d\text{\rm vol}_{S^2} = 1$ for each $j$. Then $\sup_{S^2} |S_j|_{g_{S^2}} \leq C \, \|S_j\|_{H^2} \leq C\nu_j$ for each $j$. Moreover, $\nu_j \sim j$ as $j \to \infty$ (cf. \cite{Berline-Getzler-Vergne}, Corollary 2.43). Let us write $\chi(z,t) = \sum_{j=1}^\infty \chi_j(z,t) \, S_j$, where 
\[\chi_j(z,t) = \int_{S^2} \langle \chi(z,t),S_j \rangle_{g_{S^2}} \, d\text{\rm vol}_{S^2}.\]  
Note that $|\chi_j(z,t)| \leq C \sup_{S^2} |\chi(z,t)|_{g_{S^2}}$. Moreover, the function $\chi_j(z,t)$ satisfies 
\[\frac{\partial}{\partial t} \chi_j(z,t) = \frac{\partial^2}{\partial z^2} \chi_j(z,t) - \frac{\nu_j+4}{(-2t)} \, \chi_j(z,t).\] 
Hence, the function $\hat{\chi}_j(z,t) := (-t)^{-\frac{\nu_j+4}{2}} \, \chi_j(z,t)$ satisfies 
\[\frac{\partial}{\partial t} \hat{\chi}_j(z,t) = \frac{\partial^2}{\partial z^2} \hat{\chi}_j(z,t).\] 
Moreover, $|\hat{\chi}_j(z,t)| \leq C \, (-t)^{-\frac{\nu_j+2}{2}}$ in the region $\{|z| \leq \frac{L}{2}, \, -\frac{L}{2} \leq t \leq -\frac{L}{4}\}$, and $|\hat{\chi}_j(z,t)| \leq CL^{20} \, (-t)^{-\frac{\nu_j+2}{2}}$ in the region $\{|z| \leq \frac{L}{2}, \, -\frac{L}{4} \leq t \leq -1\}$. Using the solution formula for the Dirichlet problem for the one-dimensional heat equation on the rectangle $[-\frac{L}{4},\frac{L}{4}] \times [-\frac{L}{4},-1]$, we obtain 
\begin{align*}
|\hat{\chi}_j(z,t)| 
&\leq C \, \sup_{z \in [-\frac{L}{4},\frac{L}{4}]} |\hat{\chi}_j(z,-\frac{L}{4})| \\ 
&+ C \, L \int_{-\frac{L}{4}}^t e^{-\frac{L^2}{100(t-s)}} \, (t-s)^{-\frac{3}{2}} \, \Big ( |\hat{\chi}_j(\frac{L}{4},s)| + |\hat{\chi}_j(-\frac{L}{4},s)| \Big ) \, ds, 
\end{align*} 
hence 
\begin{align*} 
|\hat{\chi}_j(z,t)| 
&\leq C \, \Big ( \frac{L}{4} \Big )^{-\frac{\nu_j+2}{2}} + CL^{21} \int_{-\frac{L}{4}}^t e^{-\frac{L^2}{100(t-s)}} \, (t-s)^{-\frac{3}{2}} \, (-s)^{-\frac{\nu_j+2}{2}} \, ds \\ 
&\leq C \, \Big ( \frac{L}{4} \Big )^{-\frac{\nu_j+2}{2}} + CL^{20} \int_{-\frac{L}{4}}^t e^{-\frac{L^2}{200(t-s)}} \, (-s)^{-\frac{\nu_j+2}{2}} \, ds \\ 
&\leq C \, \Big ( \frac{L}{4} \Big )^{-\frac{\nu_j+2}{2}} + CL^{20} \int_{-\frac{L}{4}}^{(1+\frac{1}{\sqrt{\nu_j}})t} e^{-\frac{L^2}{200(t-s)}} \, (-s)^{-\frac{\nu_j+2}{2}} \, ds \\ 
&+ CL^{20} \int_{(1+\frac{1}{\sqrt{\nu_j}})t}^t e^{-\frac{L^2}{200(t-s)}} \, (-s)^{-\frac{\nu_j+2}{2}} \, ds \\ 
&\leq C \, \Big ( \frac{L}{4} \Big )^{-\frac{\nu_j+2}{2}} + CL^{20} \, e^{-\frac{L}{100}} \, \Big ( 1+\frac{1}{\sqrt{\nu_j}} \Big )^{-\frac{\nu_j}{2}} \, (-t)^{-\frac{\nu_j}{2}} \\ 
&+ CL^{20} \, e^{-\frac{L^2\sqrt{\nu_j}}{200(-t)}} \, (-t)^{-\frac{\nu_j}{2}} 
\end{align*}
for all $t \in [-1000,-1]$ and all $z \in [-1000,1000]$. Therefore, 
\begin{align*} 
|\chi_j(z,t)| 
&\leq C \, \Big ( \frac{L}{4(-t)} \Big )^{-\frac{\nu_j+2}{2}} + CL^{20} \, e^{-\frac{L}{100}} \, \Big ( 1+\frac{1}{\sqrt{\nu_j}} \Big )^{-\frac{\nu_j}{2}} \\ 
&+ CL^{20} \, e^{-\frac{L^2\sqrt{\nu_j}}{200(-t)}}  
\end{align*}
for all $t \in [-1000,-1]$ and all $z \in [-1000,1000]$. Summation over $j$ gives  
\[|\chi(z,t)|_{g_{S^2}} = \bigg | \sum_{j=1}^\infty \chi_j(z,t) \, S_j \bigg |_{g_{S^2}} \leq C \sum_{j=1}^\infty \nu_j \, |\chi_j(z,t)| \leq C L^{-1}\] 
in the region $\{|z| \leq 1000, \, -1000 \leq t \leq -1\}$. 

\textit{Step 2:} We next analyze the equation for $\sigma(z,t)$. Let $Q_j$, $j = 1,2,\hdots$, denote the eigenfunctions of the Laplacian on vector fields on $S^2$, so that $\Delta_{S^2} Q_j = -\mu_j Q_j$. By Proposition A.1 in \cite{Brendle}, the eigenvalues satisfy $\mu_j \geq 1$. We assume that the eigenfunctions $Q_j$ are normalized so that $\int_{S^2} |Q_j|_{g_{S^2}}^2 \, d\text{\rm vol}_{S^2} = 1$ for each $j$. Then $\sup_{S^2} |Q_j|_{g_{S^2}} \leq C \, \|Q_j\|_{H^2} \leq C\mu_j$ for each $j$. Moreover, $\mu_j \sim j$ as $j \to \infty$ (cf. \cite{Berline-Getzler-Vergne}, Corollary 2.43). Let us write $\sigma(z,t) = \sum_{j=1}^\infty \sigma_j(z,t) \, Q_j$, where 
\[\sigma_j(z,t) = \int_{S^2} \langle \sigma(z,t),Q_j \rangle_{g_{S^2}} \, d\text{\rm vol}_{S^2}.\] 
Note that $|\sigma_j(z,t)| \leq C \sup_{S^2} |\sigma(z,t)|_{g_{S^2}}$. Moreover, the function $\sigma_j(z,t)$ satisfies 
\[\frac{\partial}{\partial t} \sigma_j(z,t) = \frac{\partial^2}{\partial z^2} \sigma_j(z,t) - \frac{\mu_j+1}{(-2t)} \, \sigma_j(z,t).\] 
Hence, the function $\hat{\sigma}_j(z,t) := (-t)^{-\frac{\mu_j+1}{2}} \, \sigma_j(z,t)$ satisfies 
\[\frac{\partial}{\partial t} \hat{\sigma}_j(z,t) = \frac{\partial^2}{\partial z^2} \hat{\sigma}_j(z,t).\] 
Moreover, $|\hat{\sigma}_j(z,t)| \leq C \, (-t)^{-\frac{\mu_j}{2}}$ in the region $\{|z| \leq \frac{L}{2}, \, -\frac{L}{2} \leq t \leq -\frac{L}{4}\}$, and $|\hat{\sigma}_j(z,t)| \leq CL^{20} \, (-t)^{-\frac{\mu_j+2}{2}}$ in the region $\{|z| \leq \frac{L}{2}, \, -\frac{L}{4} \leq t \leq -1\}$. Using the solution formula for the Dirichlet problem for the one-dimensional heat equation on the rectangle $[-\frac{L}{4},\frac{L}{4}] \times [-\frac{L}{4},-1]$, we obtain 
\begin{align*}
|\hat{\sigma}_j(z,t)| 
&\leq C \, \sup_{z \in [-\frac{L}{4},\frac{L}{4}]} |\hat{\sigma}_j(z,-\frac{L}{4})| \\ 
&+ C \, L \int_{-\frac{L}{4}}^t e^{-\frac{L^2}{100(t-s)}} \, (t-s)^{-\frac{3}{2}} \, \Big ( |\hat{\sigma}_j(\frac{L}{4},s)| + |\hat{\sigma}_j(-\frac{L}{4},s)| \Big ) \, ds, 
\end{align*} 
hence 
\begin{align*} 
|\hat{\sigma}_j(z,t)| 
&\leq C \, \Big ( \frac{L}{4} \Big )^{-\frac{\mu_j}{2}} + CL^{21} \int_{-\frac{L}{4}}^t e^{-\frac{L^2}{100(t-s)}} \, (t-s)^{-\frac{3}{2}} \, (-s)^{-\frac{\mu_j+2}{2}} \, ds \\ 
&\leq C \, \Big ( \frac{L}{4} \Big )^{-\frac{\mu_j}{2}} + CL^{20} \int_{-\frac{L}{4}}^t e^{-\frac{L^2}{200(t-s)}} \, (-s)^{-\frac{\mu_j+2}{2}} \, ds \\ 
&\leq C \, \Big ( \frac{L}{4} \Big )^{-\frac{\mu_j}{2}} + CL^{20} \int_{-\frac{L}{4}}^{(1+\frac{1}{\sqrt{\mu_j}})t} e^{-\frac{L^2}{200(t-s)}} \, (-s)^{-\frac{\mu_j+2}{2}} \, ds \\ 
&+ CL^{20} \int_{(1+\frac{1}{\sqrt{\mu_j}})t}^t e^{-\frac{L^2}{200(t-s)}} \, (-s)^{-\frac{\mu_j+2}{2}} \, ds \\ 
&\leq C \, \Big ( \frac{L}{4} \Big )^{-\frac{\mu_j}{2}} + CL^{20} \, e^{-\frac{L}{100}} \, \Big ( 1+\frac{1}{\sqrt{\mu_j}} \Big )^{-\frac{\mu_j}{2}} \, (-t)^{-\frac{\mu_j}{2}} \\ 
&+ CL^{20} \, e^{-\frac{L^2\sqrt{\mu_j}}{200(-t)}} \, (-t)^{-\frac{\mu_j}{2}} 
\end{align*} 
for all $t \in [-1000,-1]$ and all $z \in [-1000,1000]$. Therefore,
\begin{align*} 
|\sigma_j(z,t)| 
&\leq C \, \Big ( \frac{L}{4(-t)} \Big )^{-\frac{\mu_j}{2}} + CL^{20} \, e^{-\frac{L}{100}} \, \Big ( 1+\frac{1}{\sqrt{\mu_j}} \Big )^{-\frac{\mu_j}{2}} \\ 
&+ CL^{20} \, e^{-\frac{L^2\sqrt{\mu_j}}{200(-t)}}  
\end{align*} 
for all $t \in [-1000,-1]$ and all $z \in [-1000,1000]$. Summation over $j$ gives   
\[|\sigma(z,t)|_{g_{S^2}} = \bigg | \sum_{j=1}^\infty \sigma_j(z,t) \, Q_j \bigg |_{g_{S^2}} \leq C \sum_{j=1}^\infty \mu_j \, |\sigma_j(z,t)| \leq C L^{-\frac{1}{2}}\] 
for all $t \in [-1000,-1]$ and all $z \in [-1000,1000]$. 


\textit{Step 3:} We next analyze the equation for $\beta(z,t)$. Let $Y_j$, $j=0,1,2,\hdots$, denote the eigenfunctions of the Laplacian on scalar functions on $S^2$, so that $\Delta_{S^2} Y_j = -\lambda_j Y_j$. Note that $\lambda_0 = 0$ and $\lambda_1=2$. We assume that the eigenfunctions $Y_j$ are normalized so that $\int_{S^2} Y_j^2 \, d\text{\rm vol}_{S^2} = 1$. Then $\sup_{S^2} |Y_j| \leq C \, \|Y_j\|_{H^2} \leq C\lambda_j$ for $j \geq 1$. Moreover, $\lambda_j \sim j$ as $j \to \infty$. Let us write $\beta(z,t) = \sum_{j=0}^\infty \beta_j(z,t) \, Y_j$, where 
\[\beta_j(z,t) = \int_{S^2} \beta(z,t) \, Y_j \, d\text{\rm vol}_{S^2}.\]  
Note that $|\beta_j(z,t)| \leq C \sup_{S^2} |\beta(z,t)|$. Moreover, the function $\beta_j(z,t)$ satisfies 
\[\frac{\partial}{\partial t} \beta_j(z,t) = \frac{\partial^2}{\partial z^2} \beta_j(z,t) - \frac{\lambda_j}{(-2t)} \, \beta_j(z,t).\] 
Hence, the function $\hat{\beta}_j(z,t) := (-t)^{-\frac{\lambda_j}{2}} \, \beta_j(z,t)$ satisfies 
\[\frac{\partial}{\partial t} \hat{\beta}_j(z,t) = \frac{\partial^2}{\partial z^2} \hat{\beta}_j(z,t).\] 
In the following, we consider modes with $j \geq 1$, so that $\lambda_j \geq 2$. By assumption, $|\hat{\beta}_j(z,t)| \leq C \, (-t)^{-\frac{\lambda_j}{2}}$ in the region $\{|z| \leq \frac{L}{2}, \, -\frac{L}{2} \leq t \leq -\frac{L}{4}\}$, and $|\hat{\beta}_j(z,t)| \leq CL^{20} \, (-t)^{-\frac{\lambda_j+2}{2}}$ in the region $\{|z| \leq \frac{L}{2}, \, -\frac{L}{4} \leq t \leq -1\}$. The solution formula for the Dirichlet problem for the one-dimensional heat equation on the rectangle $[-\frac{L}{4},\frac{L}{4}] \times [-\frac{L}{4},-1]$ implies 
\begin{align*}
|\hat{\beta}_j(z,t)| 
&\leq C \, \sup_{z \in [-\frac{L}{4},\frac{L}{4}]} |\hat{\beta}_j(z,-\frac{L}{4})| \\ 
&+ C \, L \int_{-\frac{L}{4}}^t e^{-\frac{L^2}{100(t-s)}} \, (t-s)^{-\frac{3}{2}} \, \Big ( |\hat{\beta}_j(\frac{L}{4},s)| + |\hat{\beta}_j(-\frac{L}{4},s)| \Big ) \, ds 
\end{align*} 
for all $t \in [-1000,-1]$, $z \in [-1000,1000]$, and $j \geq 1$. Therefore, we obtain 
\begin{align*} 
|\hat{\beta}_j(z,t)| 
&\leq C \, \Big ( \frac{L}{4} \Big )^{-\frac{\lambda_j}{2}} + CL^{21} \int_{-\frac{L}{4}}^t e^{-\frac{L^2}{100(t-s)}} \, (t-s)^{-\frac{3}{2}} \, (-s)^{-\frac{\lambda_j+2}{2}} \, ds \\ 
&\leq C \, \Big ( \frac{L}{4} \Big )^{-\frac{\lambda_j}{2}} + CL^{20} \int_{-\frac{L}{4}}^t e^{-\frac{L^2}{200(t-s)}} \, (-s)^{-\frac{\lambda_j+2}{2}} \, ds \\ 
&\leq C \, \Big ( \frac{L}{4} \Big )^{-\frac{\lambda_j}{2}} + CL^{20} \int_{-\frac{L}{4}}^{(1+\frac{1}{\sqrt{\lambda_j}})t} e^{-\frac{L^2}{200(t-s)}} \, (-s)^{-\frac{\lambda_j+2}{2}} \, ds \\ 
&+ CL^{20} \int_{(1+\frac{1}{\sqrt{\lambda_j}})t}^t e^{-\frac{L^2}{200(t-s)}} \, (-s)^{-\frac{\lambda_j+2}{2}} \, ds \\ 
&\leq C \, \Big ( \frac{L}{4} \Big )^{-\frac{\lambda_j}{2}} + CL^{20} \, e^{-\frac{L}{100}} \, \Big ( 1+\frac{1}{\sqrt{\lambda_j}} \Big )^{-\frac{\lambda_j}{2}} \, (-t)^{-\frac{\lambda_j}{2}} \\ 
&+ CL^{20} \, e^{-\frac{L^2\sqrt{\lambda_j}}{200(-t)}} \, (-t)^{-\frac{\lambda_j}{2}} 
\end{align*}
for all $t \in [-1000,-1]$, $z \in [-1000,1000]$, and $j \geq 1$. Consequently, 
\begin{align*} 
|\beta_j(z,t)| 
&\leq C \, \Big ( \frac{L}{4(-t)} \Big )^{-\frac{\lambda_j}{2}} + CL^{20} \, e^{-\frac{L}{100}} \, \Big ( 1+\frac{1}{\sqrt{\lambda_j}} \Big )^{-\frac{\lambda_j}{2}} \\ 
&+ CL^{20} \, e^{-\frac{L^2\sqrt{\lambda_j}}{200(-t)}}  
\end{align*}
for all $t \in [-1000,-1]$, $z \in [-1000,1000]$, and $j \geq 1$. Summation over $j \geq 1$ gives  
\[|\beta(z,t) - \bar{\beta}(z,t)| = \bigg | \sum_{j=1}^\infty \beta_j(z,t) \, Y_j \bigg | \leq C \sum_{j=1}^\infty \lambda_j \, |\beta_j(z,t)| \leq C L^{-1}\] 
for all $t \in [-1000,-1]$ and all $z \in [-1000,1000]$. 

\textit{Step 4:} We finally analyze the equation for $\omega(z,t)$. As above, let $Y_j$, $j=0,1,2,\hdots$, denote the eigenfunctions of the Laplacian on scalar functions on $S^2$, so that $\Delta_{S^2} Y_j = -\lambda_j Y_j$. Note that $\lambda_0 = 0$, $\lambda_1=\lambda_2=\lambda_3=2$, and $\lambda_4=6$. We write 
$\omega(z,t) = \sum_{j=0}^\infty \omega_j(z,t) \, Y_j$, where 
\[\omega_j(z,t) = \int_{S^2} \omega(z,t) \, Y_j \, d\text{\rm vol}_{S^2}.\] 
Note that $|\omega_j(z,t)| \leq C \sup_{S^2} |\omega(z,t)|$. Moreover, the function $\omega_j(z,t)$ satisfies 
\[\frac{\partial}{\partial t} \omega_j(z,t) = \frac{\partial^2}{\partial z^2} \omega_j(z,t) - \frac{\lambda_j}{(-2t)} \, \omega_j(z,t).\] 
Hence, the function $\hat{\omega}_j(z,t) := (-t)^{-\frac{\lambda_j}{2}} \, \omega_j(z,t)$ satisfies 
\[\frac{\partial}{\partial t} \hat{\omega}_j(z,t) = \frac{\partial^2}{\partial z^2} \hat{\omega}_j(z,t).\] 
In the following, we consider modes with $j \geq 1$. We break the discussion into two subcases: 
\begin{itemize}
\item Suppose first that $j \geq 4$, so that $\lambda_j \geq 6$. By assumption, $|\hat{\omega}_j(z,t)| \leq C \, (-t)^{-\frac{\lambda_j-2}{2}}$ in the region $\{|z| \leq \frac{L}{2}, \, -\frac{L}{2} \leq t \leq -\frac{L}{4}\}$, and $|\hat{\omega}_j(z,t)| \leq CL^{20} \, (-t)^{-\frac{\lambda_j+2}{2}}$ in the region $\{|z| \leq \frac{L}{2}, \, -\frac{L}{4} \leq t \leq -1\}$. The solution formula for the one-dimensional heat equation on the rectangle $[-\frac{L}{4},\frac{L}{4}] \times [-\frac{L}{4},-1]$ implies 
\begin{align*}
|\hat{\omega}_j(z,t)| 
&\leq C \, \sup_{z \in [-\frac{L}{4},\frac{L}{4}]} |\hat{\omega}_j(z,-\frac{L}{4})| \\ 
&+ C \, L \int_{-\frac{L}{4}}^t e^{-\frac{L^2}{100(t-s)}} \, (t-s)^{-\frac{3}{2}} \, \Big ( |\hat{\omega}_j(\frac{L}{4},s)| + |\hat{\omega}_j(-\frac{L}{4},s)| \Big ) \, ds
\end{align*} 
for all $t \in [-1000,-1]$, $z \in [-1000,1000]$, and $j \geq 4$. Therefore, 
\begin{align*} 
|\hat{\omega}_j(z,t)| 
&\leq C \, \Big ( \frac{L}{4} \Big )^{-\frac{\lambda_j-2}{2}} + CL^{20} \, e^{-\frac{L}{100}} \, \Big ( 1+\frac{1}{\sqrt{\lambda_j}} \Big )^{-\frac{\lambda_j}{2}} \, (-t)^{-\frac{\lambda_j}{2}} \\ 
&+ CL^{20} \, e^{-\frac{L^2\sqrt{\lambda_j}}{200(-t)}} \, (-t)^{-\frac{\lambda_j}{2}} 
\end{align*}
for all $t \in [-1000,-1]$, $z \in [-1000,1000]$, and $j \geq 4$. Consequently, 
\begin{align*} 
|\omega_j(z,t)| 
&\leq C \, \Big ( \frac{L}{4(-t)} \Big )^{-\frac{\lambda_j-2}{2}} + CL^{20} \, e^{-\frac{L}{100}} \, \Big ( 1+\frac{1}{\sqrt{\lambda_j}} \Big )^{-\frac{\lambda_j}{2}} \\ 
&+ CL^{20} \, e^{-\frac{L^2\sqrt{\lambda_j}}{200(-t)}}  
\end{align*}
for all $t \in [-1000,-1]$, $z \in [-1000,1000]$, and $j \geq 4$. Summation over $j \geq 4$ gives 
\[\bigg | \sum_{j=4}^\infty \omega_j(z,t) \, Y_j \bigg | \leq C \sum_{j=4}^\infty \lambda_j \, |\omega_j(z,t)| \leq C L^{-2}\] 
for all $t \in [-1000,-1]$ and all $z \in [-1000,1000]$. 

\item Suppose finally that $1 \leq j \leq 3$, so that $\lambda_j=2$. In this case, $|\hat{\omega}_j(z,t)| \leq C$ in the region $\{|z| \leq \frac{L}{2}, \, -\frac{L}{2} \leq t \leq -\frac{L}{4}\}$, and $|\hat{\omega}_j(z,t)| \leq CL^{20}$ in the region $\{|z| \leq \frac{L}{2}, \, -\frac{L}{4} \leq t \leq -1\}$. Using standard interior estimates for the linear heat equation, we obtain $|\frac{\partial}{\partial z} \hat{\omega}_j(z,t)| \leq C \, (-t)^{-\frac{1}{2}}$ in the region $\{|z| \leq \frac{L}{4}, \, t=-\frac{L}{4}\}$, and $|\frac{\partial}{\partial z} \hat{\omega}_j(z,t)| \leq CL^{20} \, (-t)^{-\frac{1}{2}}$ in the region $\{|z| \leq \frac{L}{4}, \, -\frac{L}{4} \leq t \leq -1\}$. The solution formula for the one-dimensional heat equation on the rectangle $[-\frac{L}{4},\frac{L}{4}] \times [-\frac{L}{4},-1]$ implies 
\begin{align*}
&\Big | \frac{\partial}{\partial z} \hat{\omega}_j(z,t) \Big | \\ 
&\leq C \, \sup_{z \in [-\frac{L}{4},\frac{L}{4}]} \Big | \frac{\partial}{\partial z} \hat{\omega}_j(z,-\frac{L}{4}) \Big | \\ 
&+ C \, L \int_{-\frac{L}{4}}^t e^{-\frac{L^2}{100(t-s)}} \, (t-s)^{-\frac{3}{2}} \, \Big ( \Big | \frac{\partial}{\partial z} \hat{\omega}_j(\frac{L}{4},s) \Big | + \Big | \frac{\partial}{\partial z} \hat{\omega}_j(-\frac{L}{4},s) \Big | \Big ) \, ds
\end{align*} 
for all $t \in [-2000,-1]$, $z \in [-2000,2000]$, and $1 \leq j \leq 3$. Therefore,  
\[\Big | \frac{\partial}{\partial z} \hat{\omega}_j(z,t) \Big | \leq C L^{-\frac{1}{2}}\] 
for all $t \in [-2000,-1]$, $z \in [-2000,2000]$, and $1 \leq j \leq 3$. Using standard interior estimates for the linear heat equation, we obtain 
\[\Big | \frac{\partial^2}{\partial z^2} \hat{\omega}_j(z,t) \Big | \leq C L^{-\frac{1}{2}}\] 
for all $t \in [-1000,-1]$, $z \in [-1000,1000]$, and $1 \leq j \leq 3$. This implies 
\[\Big | \frac{\partial}{\partial t} \hat{\omega}_j(z,t) \Big | \leq C L^{-\frac{1}{2}}\] 
for all $t \in [-1000,-1]$, $z \in [-1000,1000]$, and $1 \leq j \leq 3$. Consequently, for each $1 \leq j \leq 3$, there exists a constant $q_j$ such that 
\[|\hat{\omega}_j(z,t) - q_j| \leq  C L^{-\frac{1}{2}}\] 
for all $t \in [-1000,-1]$ and all $z \in [-1000,1000]$. Note that $q_j$ is independent of $t$ and $z$. Thus, we conclude that 
\[|\omega_j(z,t) - (-t) q_j| \leq C L^{-\frac{1}{2}}\] 
for all $t \in [-1000,-1]$, $z \in [-1000,1000]$, and $1 \leq j \leq 3$. 
\end{itemize}
Putting these facts together, we conclude that  
\begin{align*} 
&|\omega(z,t) - \bar{\omega}(z,t) - (-t) \, (q_1Y_1+q_2Y_2+q_3Y_3)| \\ 
&= \bigg | \sum_{j=1}^3 (\omega_j(z,t) - (-t) q_j) \, Y_j + \sum_{j=4}^\infty \omega_j(z,t) \, Y_j \bigg | \leq C L^{-\frac{1}{2}} 
\end{align*} 
for all $t \in [-1000,-1]$ and all $z \in [-1000,1000]$. 

To summarize, we have shown that 
\begin{align*} 
&\big | h(t) - \bar{\omega}(z,t) \, g_{S^2} - \bar{\beta}(z,t) \, dz \otimes dz - (-t) \, (q_1Y_1+q_2Y_2+q_3Y_3) \, g_{S^2} \big |_{\bar{g}(t)} \\ 
&\leq C \, |\omega(z,t) - \bar{\omega}(z,t) - (-t) \, (q_1Y_1+q_2Y_2+q_3Y_3)| \\ 
&+ C \, |\chi(z,t)|_{g_{S^2}} + C \, |\sigma(z,t)|_{g_{S^2}} + C \, |\beta(z,t) - \bar{\beta}(z,t)| \\ 
&\leq C \, L^{-\frac{1}{2}} 
\end{align*}
in the region $\{|z| \leq 1000, \, -1000 \leq t \leq -1\}$. Hence, if we define $\psi := q_1Y_1+q_2Y_2+q_3Y_3$, then the assertion follows. \\

\section{Gluing approximate Killing vector fields}

\label{approximate.killing.vector.fields}

\begin{lemma}
\label{ode.for.U}
Let $U$ be a vector field on a Riemannian manifold, and let $\gamma$ be a unit-speed geodesic. Then $|D_s D_s U + g^{km} \, R(\gamma'(s),\partial_k,\gamma'(s),U) \, \partial_m| \leq C \, |D(\mathscr{L}_U(g))|$ along $\gamma$.
\end{lemma}

\textbf{Proof.} 
We compute 
\begin{align*} 
&|D_i D_j U_k + R_{ikjl} \, U^l + D_j D_i U_k + R_{jkil} \, U^l| \\ 
&= |D_i (D_j U_k + D_k U_j) + D_j (D_i U_k + D_k U_i) - D_k (D_i U_j + D_j U_i)| \\ 
&\leq C \, |D(\mathscr{L}_U(g))|. 
\end{align*} 
From this, the assertion follows easily. \\

\begin{lemma}
\label{estimate.for.almost.Killing.vector.fields}
Let $\bar{g}$ denote the standard metric on the cylinder $S^2 \times [-20,20]$ with scalar curvature $1$, and let $g$ be a Riemannian metric which is close to $\bar{g}$ in $C^{10}$. Let $\bar{x}$ be a point on the center slice $S^2 \times \{0\}$. Suppose that $U$ is a vector field satisfying $\sup_{B_g(\bar{x},12)} |D(\mathscr{L}_U(g))| \leq 1$ and $|U|+|DU| \leq 1$ at $\bar{x}$. Then $\sup_{B_g(\bar{x},12)} |U| \leq C$.
\end{lemma}

\textbf{Proof.} 
Let $\gamma$ be a unit-speed geodesic emanating from $\bar{x}$ with length at most $12$. By Lemma \ref{ode.for.U}, $|D_s D_s U + g^{km} \, R(\gamma'(s),\partial_k,\gamma'(s),U) \, \partial_m| \leq C$ along $\gamma$. Since $|U|+|D U| \leq 1$ at $\bar{x}$, we conclude that $|U| \leq C$ along $\gamma$. \\

\begin{lemma}
\label{estimate.for.almost.Killing.vector.fields.2}
Let $\bar{g}$ denote the standard metric on the cylinder $S^2 \times [-20,20]$ with scalar curvature $1$, and let $g$ be a Riemannian metric which is close to $\bar{g}$ in $C^{10}$. Let $\bar{x}$ be a point on the center slice $S^2 \times \{0\}$, and let $\Sigma$ denote the leaf of the CMC foliation with respect to $g$ which passes through $\bar{x}$. Suppose that $U$ is a vector field satisfying $\sup_{B_g(\bar{x},12)} |\mathscr{L}_U(g)|+|D(\mathscr{L}_U(g))| \leq 1$ and $\int_\Sigma |U|^2 \, d\mu_g \leq 1$. Then $\sup_{B_g(\bar{x},12)} |U| \leq C$.
\end{lemma}

\textbf{Proof.} 
Suppose that the assertion is false. Then there exist a sequence of metrics $g^{(j)}$ on $S^2 \times [-20,20]$ and a sequence of vector fields $U^{(j)}$ such that $g^{(j)} \to \bar{g}$ in $C^{10}$, $\sup_{B_{g^{(j)}}(\bar{x},12)} |\mathscr{L}_{U^{(j)}}(g^{(j)})|+|D(\mathscr{L}_{U^{(j)}}(g^{(j)}))| \leq 1$, $\int_{\Sigma^{(j)}} |U^{(j)}|^2 \, d\mu_{g^{(j)}} \leq 1$, and $\sup_{B_{g^{(j)}}(\bar{x},12)} |U^{(j)}| \to \infty$. Here, $\Sigma^{(j)}$ denotes the slice of the CMC foliation with respect to $g^{(j)}$ which passes through $\bar{x}$. For each $j$, we define a real number $A_j$ so that $|U^{(j)}|+|DU^{(j)}| = A_j$ at $\bar{x}$. By Lemma \ref{estimate.for.almost.Killing.vector.fields}, $\sup_{B_{g^{(j)}}(\bar{x},12)} |U^{(j)}| \leq C A_j + C$. In particular, $A_j \to \infty$. Moreover, the estimate $\sup_{B_{g^{(j)}}(\bar{x},12)} |D(\mathscr{L}_{U^{(j)}}(g^{(j)}))| \leq 1$ implies $\sup_{B_{g^{(j)}}(\bar{x},12)} |\Delta_{g^{(j)}} U^{(j)} + \text{\rm Ric}_{g^{(j)}}(U^{(j)})| \leq C$. Consequently, the rescaled vector fields $A_j^{-1} \, U^{(j)}$ converge in $C^{1,\frac{1}{2}}(B_{\bar{g}}(\bar{x},10))$ to a vector field $U$. The limiting vector field $U$ satisfies $\mathscr{L}_U(\bar{g})=0$ and $\int_{\bar{\Sigma}} |U|^2 \, d\mu_{\bar{g}}=0$. In other words, $U$ is a Killing vector field on the cylinder which vanishes along $\bar{\Sigma}$. Consequently, $U$ vanishes identically. On the other hand, $|U|+|DU|=1$ at $\bar{x}$. This is a contradiction. \\

\begin{proposition} 
\label{rigidity.of.Killing.vector.fields.on.cylinder}
If $\varepsilon_0$ is sufficiently small, then the following holds. Let $\bar{g}$ denote the standard metric on the cylinder $S^2 \times [-20,20]$ with scalar curvature $1$, let $g$ be a Riemannian metric with $\|g-\bar{g}\|_{C^{10}} \leq \varepsilon_0$, and let $\varepsilon \leq \varepsilon_0$. Let $\bar{x}$ be a point on the center slice $S^2 \times \{0\}$, and let $\Sigma$ denote the leaf of the CMC foliation with respect to $g$ which passes through $\bar{x}$. Suppose that $U^{(1)},U^{(2)},U^{(3)}$ are vector fields with the following properties: 
\begin{itemize}
\item $\sup_{B_g(\bar{x},12)} \sum_{a=1}^3 |\mathscr{L}_{U^{(a)}}(g)|^2+|D(\mathscr{L}_{U^{(a)}}(g))|^2 \leq \varepsilon^2$.
\item $\sup_\Sigma \sum_{a=1}^3 |\langle U^{(a)},\nu \rangle|^2 \leq \varepsilon^2$.
\item $\sum_{a,b=1}^3 \big | \delta_{ab} - \text{\rm area}_g(\Sigma)^{-2} \int_\Sigma \langle U^{(a)},U^{(b)} \rangle \, d\mu_g \big |^2 \leq \varepsilon^2$.
\end{itemize}
Moreover, suppose that $\tilde{U}^{(1)},\tilde{U}^{(2)},\tilde{U}^{(3)}$ are vector fields with the following properties: 
\begin{itemize}
\item $\sup_{B_g(\bar{x},12)} \sum_{a=1}^3 |\mathscr{L}_{\tilde{U}^{(a)}}(g)|^2+|D(\mathscr{L}_{\tilde{U}^{(a)}}(g))|^2 \leq \varepsilon^2$.
\item $\sup_\Sigma \sum_{a=1}^3 |\langle \tilde{U}^{(a)},\nu \rangle|^2 \leq \varepsilon^2$.
\item $\sum_{a,b=1}^3 \big | \delta_{ab} - \text{\rm area}_g(\Sigma)^{-2} \int_\Sigma \langle \tilde{U}^{(a)},\tilde{U}^{(b)} \rangle \, d\mu_g \big |^2 \leq \varepsilon^2$.
\end{itemize}
Then there exists a $3 \times 3$ matrix $\omega \in O(3)$ such that 
\[\sup_{B_g(\bar{x},9)} \sum_{a=1}^3 \Big | \sum_{b=1}^3 \omega_{ab} \, U^{(b)} - \tilde{U}^{(a)} \Big |^2 \leq C\varepsilon^2.\]
\end{proposition}

\textbf{Proof.} 
Suppose that the assertion is false. Then we can find a sequence of metrics $g^{(j)}$ on $S^2 \times [-20,20]$, a collection of vector fields $U^{(1,j)},U^{(2,j)},U^{(3,j)}$, a collection of vector fields $\tilde{U}^{(1,j)},\tilde{U}^{(2,j)},\tilde{U}^{(3,j)}$, and a sequence of positive numbers $\varepsilon_j$ with the following properties: 
\begin{itemize}
\item $\|g^{(j)}-\bar{g}\|_{C^{10}} \leq j^{-1}$ and $\varepsilon_j \leq j^{-1}$.
\item $\sup_{B_{g^{(j)}}(\bar{x},12)} \sum_{a=1}^3 |\mathscr{L}_{U^{(a,j)}}(g^{(j)})|^2+|D(\mathscr{L}_{U^{(a,j)}}(g^{(j)}))|^2 \leq \varepsilon_j^2$.
\item $\sup_{\Sigma^{(j)}} \sum_{a=1}^3 |\langle U^{(a,j)},\nu \rangle|^2 \leq \varepsilon_j^2$.
\item $\sum_{a,b=1}^3 \big | \delta_{ab} - \text{\rm area}_g(\Sigma^{(j)})^{-2} \int_{\Sigma^{(j)}} \langle U^{(a,j)},U^{(b,j)} \rangle \, d\mu_{g^{(j)}} \big |^2 \leq \varepsilon_j^2$.
\item $\sup_{B_{g^{(j)}}(\bar{x},12)} \sum_{a=1}^3 |\mathscr{L}_{\tilde{U}^{(a,j)}}(g^{(j)})|^2+|D(\mathscr{L}_{\tilde{U}^{(a,j)}}(g^{(j)}))|^2 \leq \varepsilon_j^2$.
\item $\sup_{\Sigma^{(j)}} \sum_{a=1}^3 |\langle \tilde{U}^{(a,j)},\nu \rangle|^2 \leq \varepsilon_j^2$.
\item $\sum_{a,b=1}^3 \big | \delta_{ab} - \text{\rm area}_g(\Sigma^{(j)})^{-2} \int_{\Sigma^{(j)}} \langle \tilde{U}^{(a,j)},\tilde{U}^{(b,j)} \rangle \, d\mu_{g^{(j)}} \big |^2 \leq \varepsilon_j^2$.
\item $\delta_j^2 := \inf_{\omega \in O(3)} \sup_{B_{g^{(j)}}(\bar{x},9)} \sum_{a=1}^3 \big | \sum_{b=1}^3 \omega_{ab} \, U^{(b,j)} - \tilde{U}^{(a,j)} \big |^2 \geq (j \, \varepsilon_j)^2$.
\end{itemize} 
Here, $\Sigma^{(j)}$ denotes the leaf of the CMC foliation with respect to $g^{(j)}$ which passes through $\bar{x}$. 

Clearly, $\int_{\Sigma^{(j)}} \sum_{a=1}^3 |U^{(a,j)}|^2 \, d\mu_{g^{(j)}} \leq C$ and $\int_{\Sigma^{(j)}} \sum_{a=1}^3 |\tilde{U}^{(a,j)}|^2 \, d\mu_{g^{(j)}} \leq C$. Hence, Lemma \ref{estimate.for.almost.Killing.vector.fields.2} implies that $\sup_{B_{g^{(j)}}(\bar{x},12)} \sum_{a=1}^3 |U^{(a,j)}|^2 \leq C$ and $\sup_{B_{g^{(j)}}(\bar{x},12)} \sum_{a=1}^3 |\tilde{U}^{(a,j)}|^2 \leq C$. Moreover, $\sup_{B_{g^{(j)}}(\bar{x},12)} \sum_{a=1}^3 |\Delta_{g^{(j)}} U^{(a,j)} + \text{\rm Ric}_{g^{(j)}}(U^{(a,j)})|^2 \leq C \varepsilon_j^2$ and $\sup_{B_{g^{(j)}}(\bar{x},12)} \sum_{a=1}^3 |\Delta_{g^{(j)}} \tilde{U}^{(a,j)} + \text{\rm Ric}_{g^{(j)}}(\tilde{U}^{(a,j)})|^2 \leq C \varepsilon_j^2$. After passing to a subsequence, the vector fields $U^{(a,j)}$ converge in $C^{1,\frac{1}{2}}(B_{\bar{g}}(\bar{x},10))$ to a vector field $U^{(a)}$ which satisfies $\mathscr{L}_{U^{(a)}}(\bar{g})=0$ and is tangential along $\bar{\Sigma}$. Similarly, the vector fields $\tilde{U}^{(a,j)}$ converge in $C^{1,\frac{1}{2}}(B_{\bar{g}}(\bar{x},10))$ to a vector field $\tilde{U}^{(a)}$ which satisfies $\mathscr{L}_{\tilde{U}^{(a)}}(\bar{g})=0$ and is tangential along $\bar{\Sigma}$. Note that $\text{\rm area}_{\bar{g}}(\bar{\Sigma})^{-2} \int_{\bar{\Sigma}} \langle U^{(a)},U^{(b)} \rangle \, d\mu_{\bar{g}} = \text{\rm area}_{\bar{g}}(\bar{\Sigma})^{-2} \int_{\bar{\Sigma}} \langle \tilde{U}^{(a)},\tilde{U}^{(b)} \rangle \, d\mu_{\bar{g}} = \delta_{ab}$. Consequently, there exists a matrix $\bar{\omega} \in O(3)$ such that $\tilde{U}^{(a)} = \sum_{b=1}^3 \bar{\omega}_{ab} \, U^{(b)}$. This implies $\delta_j \to 0$.

For each $j$, we choose a $3 \times 3$ matrix $\omega^{(j)} \in O(3)$ such that  
\[\sup_{B_{g^{(j)}}(\bar{x},9)} \sum_{a=1}^3 \Big | \sum_{b=1}^3 \omega_{ab}^{(j)} \, U^{(b,j)} - \tilde{U}^{(a,j)} \Big |^2 = \delta_j^2.\] 
Clearly, $\omega^{(j)} \to \bar{\omega}$ as $j \to \infty$. We next define 
\[V^{(a,j)} := \delta_j^{-1} \, \Big ( \sum_{b=1}^3 \omega_{ab}^{(j)} \, U^{(b,j)} - \tilde{U}^{(a,j)} \Big ).\] 
The vector fields $V^{(a,j)}$ have the following properties: 
\begin{itemize}
\item $\sup_{B_{g^{(j)}}(\bar{x},9)} \sum_{a=1}^3 |V^{(a,j)}|^2 = 1$.
\item $\sup_{B_{g^{(j)}}(\bar{x},12)} \sum_{a=1}^3 |\mathscr{L}_{V^{(a,j)}}(g^{(j)})|^2+|D(\mathscr{L}_{V^{(a,j)}}(g^{(j)}))|^2 \leq C \, j^{-2}$.
\item $\sup_{\Sigma^{(j)}} \sum_{a=1}^3 |\langle V^{(a,j)},\nu \rangle|^2 \leq C \, j^{-2}$.
\end{itemize} 
Using Lemma \ref{estimate.for.almost.Killing.vector.fields.2}, we obtain $\sup_{B_{g^{(j)}}(\bar{x},12)} \sum_{a=1}^3 |V^{(a,j)}|^2 \leq C$. Moreover, $\sup_{B_{g^{(j)}}(\bar{x},12)} \sum_{a=1}^3 |\Delta_{g^{(j)}} V^{(a,j)} + \text{\rm Ric}_{g^{(j)}}(V^{(a,j)})|^2 \leq C \, j^{-2}$. Hence, after passing to a subsequence, the vector fields $V^{(a,j)}$ converge in $C^{1,\frac{1}{2}}(B_{\bar{g}}(\bar{x},10))$ to a vector field $V^{(a)}$ which satisfies $\mathscr{L}_{V^{(a)}}(\bar{g}) = 0$ and is tangential along $\bar{\Sigma}$. Consequently, $V^{(a)} = \sum_{b=1}^3 \sigma_{ab} \, \tilde{U}^{(b)}$ for some $3 \times 3$-matrix $\sigma$. 

We next observe that 
\begin{align*} 
&\bigg | \delta_j \int_{\Sigma^{(j)}} (\langle \tilde{U}^{(a,j)},V^{(b,j)} \rangle + \langle V^{(a,j)},\tilde{U}^{(b,j)} \rangle) \, d\mu_{g^{(j)}} + \delta_j^2 \int_{\Sigma^{(j)}} \langle V^{(a,j)},V^{(b,j)} \rangle \, d\mu_{g^{(j)}} \bigg | \\ 
&= \bigg | \int_{\Sigma^{(j)}} \langle \tilde{U}^{(a,j)} + \delta_j \, V^{(a,j)},\tilde{U}^{(b,j)} + \delta_j \, V^{(b,j)} \rangle \, d\mu_{g^{(j)}} - \int_{\Sigma^{(j)}} \langle \tilde{U}^{(a,j)},\tilde{U}^{(b,j)} \rangle \, d\mu_{g^{(j)}} \bigg | \\ 
&= \bigg | \sum_{c,d=1}^3 \omega_{ac}^{(j)} \omega_{bd}^{(j)} \int_{\Sigma^{(j)}} \langle U^{(c,j)},U^{(d,j)} \rangle \, d\mu_{g^{(j)}} - \int_{\Sigma^{(j)}} \langle \tilde{U}^{(a,j)},\tilde{U}^{(b,j)} \rangle \, d\mu_{g^{(j)}} \bigg | \\ 
&\leq C\varepsilon_j. 
\end{align*} 
Since $\delta_j \to 0$ and $\delta_j^{-1} \varepsilon_j \leq j^{-1}$, we conclude that 
\[\int_{\bar{\Sigma}} (\langle \tilde{U}^{(a)},V^{(b)} \rangle + \langle V^{(a)},\tilde{U}^{(b)} \rangle) \, d\mu_{\bar{g}} = 0.\] 
Consequently, $\sigma$ is an anti-symmetric matrix. Let $\tilde{\omega}^{(j)} := \exp(-\delta_j \sigma) \, \omega^{(j)} \in O(3)$. Since $V^{(a)} = \sum_{b=1}^3 \sigma_{ab} \, \tilde{U}^{(b)}$, we obtain 
\[\sup_{B_{g^{(j)}}(\bar{x},9)} \sum_{a=1}^3 \Big | V^{(a,j)} + \delta_j^{-1} \sum_{b=1}^3 (\tilde{\omega}_{ab}^{(j)}-\omega_{ab}^{(j)}) \, U^{(b,j)} \Big |^2 \to 0\] 
as $j \to \infty$. On the other hand, it follows from the definition of $\delta_j$ that 
\[\sup_{B_{g^{(j)}}(\bar{x},9)} \sum_{a=1}^3 \Big | V^{(a,j)} + \delta_j^{-1} \sum_{b=1}^3 (\omega_{ab}-\omega_{ab}^{(j)}) \, U^{(b,j)} \Big |^2 \geq 1\] 
for each $j$ and each $\omega \in O(3)$. This is a contradiction. \\

\begin{corollary}
\label{gluing}
Let $\bar{g}$ denote the standard metric on the cylinder $S^2 \times [-20,20]$ with scalar curvature $1$, let $g$ be a Riemannian metric with $\|g-\bar{g}\|_{C^{10}} \leq \varepsilon_0$, and let $\varepsilon \leq \varepsilon_0$. Moreover, suppose that $U^{(1)},U^{(2)},U^{(3)}$ and $\tilde{U}^{(1)},\tilde{U}^{(2)},\tilde{U}^{(3)}$ are vector fields satisfying the assumptions of Proposition \ref{rigidity.of.Killing.vector.fields.on.cylinder}. Let $\eta$ be a smooth cutoff function such that $\eta = 1$ on $S^2 \times [-1000,-1]$ and $\eta = 0$ on $S^2 \times [1,1000]$. Then there exists a $3 \times 3$ matrix $\omega \in O(3)$ with the property that the vector fields $V^{(a)} := \eta \sum_{b=1}^3 \omega_{ab} \, U^{(b)} + (1-\eta) \, \tilde{U}^{(a)}$ satisfy $\sum_{a=1}^3 |\mathscr{L}_{V^{(a)}}(g)|^2 + |D(\mathscr{L}_{V^{(a)}}(g))|^2 \leq C\varepsilon^2$ in the transition region $S^2 \times [-1,1]$. 
\end{corollary} 

\textbf{Proof.} 
By Proposition \ref{rigidity.of.Killing.vector.fields.on.cylinder}, we can find a $3 \times 3$ matrix $\omega \in O(3)$ with the property that the vector fields 
\[W^{(a)} := \sum_{b=1}^3 \omega_{ab} \, U^{(b)} - \tilde{U}^{(a)}.\] 
satisfy $\sup_{B_{\bar{g}}(\bar{x},9)} \sum_{a=1}^3 |W^{(a)}|^2 \leq C\varepsilon^2$. Moreover, $\sup_{B_{\bar{g}}(\bar{x},9)} \sum_{a=1}^3 |\Delta W^{(a)} + \text{\rm Ric}(W^{(a)})|^2 \leq C\varepsilon^2$. Using standard interior estimates for elliptic equations, we obtain $\sup_{B_{\bar{g}}(\bar{x},8)} |DW^{(a)}|^2 \leq C\varepsilon^2$. We now define 
\[V^{(a)} := \eta \sum_{b=1}^3 \omega_{ab} \, U^{(b)} + (1-\eta) \, \tilde{U}^{(a)}.\] 
Then 
\begin{align*} 
\mathscr{L}_{V^{(a)}}(g) 
&= \eta \sum_{b=1}^3 \omega_{ab} \, \mathscr{L}_{U^{(b)}}(g) + (1-\eta) \, \mathscr{L}_{\tilde{U}^{(a)}}(g) \\ 
&+ d\eta \otimes g(W^{(a)},\cdot) + g(W^{(a)},\cdot) \otimes d\eta. 
\end{align*}
Using the estimate $\sup_{B_{\bar{g}}(\bar{x},8)} \sum_{a=1}^3 |W^{(a)}|^2 +|DW^{(a)}|^2 \leq C\varepsilon^2$, we conclude that $\sup_{B_{\bar{g}}(\bar{x},8)} \sum_{a=1}^3 |\mathscr{L}_{V^{(a)}}(g)|^2 + |D(\mathscr{L}_{V^{(a)}}(g))|^2 \leq C\varepsilon^2$. Since the transition region $S^2 \times [-1,1]$ is contained in $B_{\bar{g}}(\bar{x},8)$, the assertion follows. \\

\section{The Neck Improvement Theorem}

\label{NIT}

\begin{definition}
\label{evolving.neck}
Let $(M,g(t))$ be a solution to the Ricci flow in dimension $3$, and let $(\bar{x},\bar{t})$ be a point in space-time with $R(\bar{x},\bar{t}) = r^{-2}$. We say that $(\bar{x},\bar{t})$ lies at the center of an evolving $\varepsilon$-neck if, after rescaling by the factor $r^{-1}$, the parabolic neighborhood $B_{g(\bar{t})}(\bar{x},\varepsilon^{-1} r) \times [\bar{t}-\varepsilon^{-1} r^2,\bar{t}]$ is $\varepsilon$-close in $C^{[\varepsilon^{-1}]}$ to a family of shrinking cylinders.
\end{definition}

\begin{definition} 
\label{symmetry.of.necks}
Let $(M,g(t))$ be a solution to the Ricci flow in dimension $3$, and let $(\bar{x},\bar{t})$ be a point in space-time with $R(\bar{x},\bar{t}) = r^{-2}$. We assume that $(\bar{x},\bar{t})$ lies at the center of an evolving $\varepsilon _0$-neck for some small number $\varepsilon_0$. We say that $(\bar{x},\bar{t})$ is $\varepsilon$-symmetric if there exist smooth, time-independent vector fields $U^{(1)},U^{(2)},U^{(3)}$ which are defined on an open set containing $\bar{B}_{g(\bar{t})}(\bar{x},100r)$ and satisfy the following conditions: 
\begin{itemize}
\item $\sup_{\bar{B}_{g(\bar{t})}(\bar{x},100r) \times [\bar{t}-100r^2,\bar{t}]} \sum_{l=0}^2 \sum_{a=1}^3 r^{2l} \, |D^l(\mathscr{L}_{U^{(a)}}(g(t)))|^2 \leq \varepsilon^2$.
\item If $t \in [\bar{t}-100r^2,\bar{t}]$ and $\Sigma \subset \bar{B}_{g(\bar{t})}(\bar{x},100r)$ is a leaf of the CMC foliation of $(M,g(t))$, then $\sup_\Sigma \sum_{a=1}^3 r^{-2} \, |\langle U^{(a)},\nu \rangle|^2 \leq \varepsilon^2$, where $\nu$ denotes the unit normal vector to $\Sigma$ in $(M,g(t))$.
\item If $t \in [\bar{t}-100r^2,\bar{t}]$ and $\Sigma \subset \bar{B}_{g(\bar{t})}(\bar{x},100r)$ is a leaf of the CMC foliation of $(M,g(t))$, then 
\[\sum_{a,b=1}^3 \bigg | \delta_{ab} - \text{\rm area}_{g(t)}(\Sigma)^{-2} \int_\Sigma \langle U^{(a)},U^{(b)} \rangle_{g(t)} \, d\mu_{g(t)} \bigg |^2 \leq \varepsilon^2.\]
\end{itemize}
\end{definition}

\begin{lemma} 
\label{openness.1}
Suppose that $(\bar{x},\bar{t})$ is a point in spacetime which is $\varepsilon$-symmetric. If $(\tilde{x},\tilde{t})$ is sufficiently close to $(\bar{x},\bar{t})$, then $(\tilde{x},\tilde{t})$ is $2\varepsilon$-symmetric.
\end{lemma} 

\textbf{Proof.} 
This follows immediately from the definition. \\

\begin{lemma}
\label{constructing.approximate.killing.vector.fields}
If $L$ is sufficiently large and $\varepsilon_0$ is sufficiently small depending on $L$, then the following holds. Let $(M,g(t))$ be a solution of the Ricci flow in dimension $3$, and let $(x_0,-1)$ be a point in space-time which lies at the center of an evolving $\varepsilon_0$-neck and satisfies $R(x_0,-1) = 1$. Moreover, we assume that every point in the parabolic neighborhood $B_{g(-1)}(x_0,L) \times [-L-1,-1)$ is $\varepsilon$-symmetric for some positive real number $\varepsilon \leq \varepsilon_0$. Then, given any $\bar{t} \in [-\frac{L}{10},-1]$, we can find time-independent vector fields $U^{(1)},U^{(2)},U^{(3)}$ on $B_{g(-1)}(x_0,\frac{127L}{128})$ with the following properties: 
\begin{itemize}
\item $\sum_{a=1}^3 |\mathscr{L}_{U^{(a)}}(g(t))|^2+(-t) \, |D(\mathscr{L}_{U^{(a)}}(g(t)))|^2 \leq C\varepsilon^2$ on $B_{g(-1)}(x_0,\frac{127L}{128}) \times [10\bar{t},\bar{t}]$.
\item $\sum_{a=1}^3 (-t)^{-1} \, |\langle U^{(a)},\nu \rangle|^2 \leq C\varepsilon^2$ on $B_{g(-1)}(x_0,\frac{127L}{128}) \times [10\bar{t},\bar{t}]$, where $\nu$ denotes the unit normal vector to the CMC foliation of $(M,g(t))$.
\item If $t \in [10\bar{t},\bar{t}]$ and $\Sigma \subset B_{g(-1)}(x_0,\frac{127L}{128})$ is a leaf of the CMC foliation of $(M,g(t))$, then 
\[\sum_{a,b=1}^3 \bigg | \delta_{ab} - \text{\rm area}_{g(t)}(\Sigma)^{-2} \int_\Sigma \langle U^{(a)},U^{(b)} \rangle_{g(t)} \, d\mu_{g(t)} \bigg |^2 \leq C \varepsilon^2.\]
\end{itemize}
Moreover, $U^{(1)},U^{(2)},U^{(3)}$ are $C(L)\varepsilon_0$-close to the standard rotation vector fields on the cylinder in the $C^2$-norm.
\end{lemma}

\textbf{Proof.} 
We proceed in two steps:

\textit{Step 1:} Suppose first that $\bar{t} \in [-\frac{L}{10},-1)$. By assumption, the point $(\bar{x},\bar{t})$ is $\varepsilon$-symmetric whenever $\bar{x} \in B_{g(-1)}(x_0,L)$. By a repeated application of Corollary \ref{gluing}, we can construct vector fields $U^{(1)},U^{(2)},U^{(3)}$ satisfying the conditions above. Moreover, in view of Definition \ref{symmetry.of.necks}, the Lie derivatives $\mathscr{L}_{U^{(1)}}(g),\mathscr{L}_{U^{(2)}}(g),\mathscr{L}_{U^{(3)}}(g)$ are bounded by $C(L)\varepsilon$ in the $C^2$-norm. Consequently, the vector fields $U^{(1)},U^{(2)},U^{(3)}$ are $C(L)\varepsilon_0$-close to the standard rotation vector fields on the cylinder in the $C^{2,\frac{1}{2}}$-norm.

\textit{Step 2:} Suppose next that $\bar{t} = -1$. In this case, the assertion follows from the result in Step 1 by passing to the limit. Since the vector fields constructed in Step 1 are bounded in $C^{2,\frac{1}{2}}$, we may take the limit in $C^2$. \\

\begin{lemma}
\label{comparing.approximate.killing.vector.fields}
If $L$ is sufficiently large and $\varepsilon_0$ is sufficiently small depending on $L$, then the following holds. Let $(M,g(t))$ be a solution of the Ricci flow in dimension $3$, and let $(x_0,-1)$ be a point in space-time which lies at the center of an evolving $\varepsilon_0$-neck and satisfies $R(x_0,-1) = 1$. Consider a time $\bar{t} \in [-L,-1]$ and a positive real number $\varepsilon \leq \varepsilon_0$. Suppose that $U^{(1)},U^{(2)},U^{(3)}$ are time-independent vector fields on $B_{g(-1)}(x_0,\frac{127L}{128})$ with the following properties: 
\begin{itemize}
\item $\sum_{a=1}^3 |\mathscr{L}_{U^{(a)}}(g(\bar{t}))|^2+(-\bar{t}) \, |D(\mathscr{L}_{U^{(a)}}(g(\bar{t})))|^2 \leq \varepsilon^2$ on $B_{g(-1)}(x_0,\frac{127L}{128})$.
\item $\sum_{a=1}^3 (-\bar{t})^{-1} \, |\langle U^{(a)},\nu \rangle|^2 \leq \varepsilon^2$ on $B_{g(-1)}(x_0,\frac{127L}{128})$, where $\nu$ denotes the unit normal vector to the CMC foliation of $(M,g(\bar{t}))$.
\item If $\Sigma \subset B_{g(-1)}(x_0,\frac{127L}{128})$ is a leaf of the CMC foliation of $(M,g(\bar{t}))$, then 
\[\sum_{a,b=1}^3 \bigg | \delta_{ab} - \text{\rm area}_{g(\bar{t})}(\Sigma)^{-2} \int_\Sigma \langle U^{(a)},U^{(b)} \rangle_{g(\bar{t})} \, d\mu_{g(\bar{t})} \bigg |^2 \leq \varepsilon^2.\]
\end{itemize}
Moreover, suppose that $\tilde{U}^{(1)},\tilde{U}^{(2)},\tilde{U}^{(3)}$ are time-independent vector fields on $B_{g(-1)}(x_0,\frac{127L}{128})$ with the following properties: 
\begin{itemize}
\item $\sum_{a=1}^3 |\mathscr{L}_{\tilde{U}^{(a)}}(g(\bar{t}))|^2+(-\bar{t}) \, |D(\mathscr{L}_{\tilde{U}^{(a)}}(g(\bar{t})))|^2 \leq \varepsilon^2$ on $B_{g(-1)}(x_0,\frac{127L}{128})$.
\item $\sum_{a=1}^3 (-\bar{t})^{-1} \, |\langle \tilde{U}^{(a)},\nu \rangle|^2 \leq \varepsilon^2$ on $B_{g(-1)}(x_0,\frac{127L}{128})$, where $\nu$ denotes the unit normal vector to the CMC foliation of $(M,g(\bar{t}))$.
\item If $\Sigma \subset B_{g(-1)}(x_0,\frac{127L}{128})$ is a leaf of the CMC foliation of $(M,g(\bar{t}))$, then 
\[\sum_{a,b=1}^3 \bigg | \delta_{ab} - \text{\rm area}_{g(\bar{t})}(\Sigma)^{-2} \int_\Sigma \langle \tilde{U}^{(a)},\tilde{U}^{(b)} \rangle_{g(\bar{t})} \, d\mu_{g(\bar{t})} \bigg |^2 \leq \varepsilon^2.\]
\end{itemize}
Then there exists a $3 \times 3$ matrix $\omega \in O(3)$ such that 
\[\sup_{B_{g(-1)}(x_0,\frac{31L}{32})} (-\bar{t})^{-1} \sum_{a=1}^3 \Big | \sum_{b=1}^3 \omega_{ab} \, U^{(b)} - \tilde{U}^{(a)} \Big |_{g(\bar{t})}^2 \leq C L^2 \varepsilon^2.\] 
\end{lemma}

\textbf{Proof.} 
By assumption, the flow is close to a family of shrinking cylinders. For each integer $m \in [-\frac{63L}{64},\frac{63L}{64}]$, Proposition \ref{rigidity.of.Killing.vector.fields.on.cylinder} implies that there exists a $3 \times 3$ matrix $\omega^{(m)} \in O(3)$ such that 
\[\sup_{S^2 \times [m-1,m+1]} (-\bar{t})^{-1} \sum_{a=1}^3 \Big | \sum_{b=1}^3 \omega_{ab}^{(m)} \, U^{(b)} - \tilde{U}^{(a)} \Big |_{g(\bar{t})}^2 \leq C \varepsilon^2.\] 
From this, we deduce that $|\omega^{(m)}-\omega^{(m+1)}| \leq C\varepsilon$ for every integer $m$. Consequently, there exists a $3 \times 3$ matrix $\omega \in O(3)$ such that $|\omega^{(m)}-\omega| \leq CL\varepsilon$ for every integer $m \in [-\frac{63L}{64},\frac{63L}{64}]$. This implies 
\[\sup_{S^2 \times [m-1,m+1]} (-\bar{t})^{-1} \sum_{a=1}^3 \Big | \sum_{b=1}^3 \omega_{ab} \, U^{(b)} - \tilde{U}^{(a)} \Big |_{g(\bar{t})}^2 \leq C L^2 \varepsilon^2\] 
for every integer $m \in [-\frac{63L}{64},\frac{63L}{64}]$. This completes the proof of Lemma \ref{comparing.approximate.killing.vector.fields}. \\

We now state the main result of this section:

\begin{theorem}[Neck Improvement Theorem] 
\label{neck.improvement.theorem}
There exist a large constant $L$ and small positive constant $\varepsilon_1$ with the following property. Let $(M,g(t))$ be a solution of the Ricci flow in dimension $3$, and let $(x_0,t_0)$ be a point in space-time which lies at the center of an evolving $\varepsilon_1$-neck and satisfies $R(x_0,t_0) = r^{-2}$. Moreover, suppose that every point in the parabolic neighborhood $B_{g(t_0)}(x_0,Lr) \times [t_0-Lr^2,t_0)$ is $\varepsilon$-symmetric, where $\varepsilon \leq \varepsilon_1$. Then $(x_0,t_0)$ is $\frac{\varepsilon}{2}$-symmetric.
\end{theorem}

\textbf{Proof.} 
Throughout the proof, we will assume that $L$ is sufficiently large, and $\varepsilon_1$ is sufficiently small depending on $L$. Without loss of generality, we may assume that $t_0=-1$ and $R(x_0,-1) = 1$. In the parabolic neighborhood $B_{g(-1)}(x_0,L) \times [-L-1,-1]$, the metric $g(t)$ is $\varepsilon_1$-close to a family of shrinking cylinders in the $C^{100}$-norm. Let $\bar{g}(t) = (-2t) \, g_{S^2} + dz \otimes dz$ denote the standard metric on the shrinking cylinders. 

\textit{Step 1:} Using Lemma \ref{constructing.approximate.killing.vector.fields} and Lemma \ref{comparing.approximate.killing.vector.fields}, we can construct time-dependent vector fields $U^{(1)},U^{(2)},U^{(3)}$, defined on $B_{g(-1)}(x_0,\frac{15L}{16}) \times [-L,-1]$, with the following properties: 
\begin{itemize}
\item $\frac{\partial}{\partial t} U^{(a)} = 0$ on $B_{g(-1)}(x_0,\frac{15L}{16}) \times [-L,-\frac{L}{4}]$, and $|\frac{\partial}{\partial t} U^{(a)}| \leq C L (-t)^{-\frac{1}{2}} \varepsilon$ on $B_{g(-1)}(x_0,\frac{15L}{16}) \times [-\frac{L}{4},-1]$. 
\item $|\mathscr{L}_{U^{(a)}}(g(t))|+(-t)^{\frac{1}{2}} \, |D(\mathscr{L}_{U^{(a)}}(g(t)))| \leq C\varepsilon$ on $B_{g(-1)}(x_0,\frac{15L}{16}) \times [-L,-1]$.
\end{itemize}
Moreover, we can arrange that $U^{(1)},U^{(2)},U^{(3)}$ are $C(L)\varepsilon_1$-close to the standard rotation vector fields on the cylinder in the $C^2$-norm. Note that $|\Delta U^{(a)} + \text{\rm Ric}(U^{(a)})| \leq C \, |D(\mathscr{L}_{U^{(a)}}(g))| \leq C(-t)^{-\frac{1}{2}} \varepsilon$ on $B_{g(-1)}(x_0,\frac{15L}{16}) \times [-L,-1]$.

\textit{Step 2:} Let $V^{(a)}$ denote the solution of the PDE 
\[\frac{\partial}{\partial t} V^{(a)} = \Delta V^{(a)} + \text{\rm Ric}(V^{(a)})\] 
in the region $\{|z| \leq \frac{7L}{8}, \, -L \leq t \leq -1\}$ with Dirichlet boundary condition $V^{(a)}=U^{(a)}$ on the parabolic boundary $\{|z| \leq \frac{7L}{8}, \, t=-L\} \cup \{|z|=\frac{7L}{8}, \, -L \leq t \leq -1\}$. The difference $V^{(a)} - U^{(a)}$ satisfies 
\begin{align*} 
&|\frac{\partial}{\partial t} (V^{(a)}-U^{(a)}) - \Delta (V^{(a)}-U^{(a)}) - \text{\rm Ric}(V^{(a)}-U^{(a)})| \\ 
&= |\Delta U^{(a)}+\text{\rm Ric}(U^{(a)})| \leq C(-t)^{-\frac{1}{2}}\varepsilon 
\end{align*}
in the region $\{|z| \leq \frac{7L}{8}, \, -L \leq t \leq -\frac{L}{4}\}$, and   
\begin{align*} 
&|\frac{\partial}{\partial t} (V^{(a)}-U^{(a)}) - \Delta (V^{(a)}-U^{(a)}) - \text{\rm Ric}(V^{(a)}-U^{(a)})| \\ 
&\leq |\frac{\partial}{\partial t} U^{(a)}| + |\Delta U^{(a)}+\text{\rm Ric}(U^{(a)})| \leq CL(-t)^{-\frac{1}{2}}\varepsilon 
\end{align*} 
in the region $\{|z| \leq \frac{7L}{8}, \, -\frac{L}{4} \leq t \leq -1\}$. Hence, Proposition \ref{estimate.for.norm.of.vector.field} implies 
\[\frac{\partial}{\partial t} |V^{(a)}-U^{(a)}| \leq \Delta |V^{(a)}-U^{(a)}| + C(-t)^{-\frac{1}{2}}\varepsilon\] 
in the region $\{|z| \leq \frac{7L}{8}, \, -L \leq t \leq -\frac{L}{4}\}$, and 
\[\frac{\partial}{\partial t} |V^{(a)}-U^{(a)}| \leq \Delta |V^{(a)}-U^{(a)}| + CL(-t)^{-\frac{1}{2}}\varepsilon\] 
in the region $\{|z| \leq \frac{7L}{8}, \, -\frac{L}{4} \leq t \leq -1\}$. Using the maximum principle, we obtain 
\[|V^{(a)}-U^{(a)}| \leq CL^{\frac{1}{2}}\varepsilon\] 
in the region $\{|z| \leq \frac{7L}{8}, \, -L \leq t \leq -\frac{L}{4}\}$, and 
\[|V^{(a)}-U^{(a)}| \leq CL^2\varepsilon\] 
in the region $\{|z| \leq \frac{7L}{8}, \, -\frac{L}{4} \leq t \leq -1\}$. Standard interior estimates for linear parabolic equations imply 
\[|D(V^{(a)}-U^{(a)})| \leq C\varepsilon\] 
in the region $\{|z| \leq \frac{3L}{4}, \, -\frac{3L}{4} \leq t \leq -\frac{L}{4}\}$, and 
\[|D(V^{(a)}-U^{(a)})| \leq CL^2 \varepsilon\] 
in the region $\{|z| \leq \frac{3L}{4}, \, -\frac{L}{4} \leq t \leq -1\}$. In particular, in the region $\{|z| \leq \frac{3L}{4}, \, -\frac{3L}{4} \leq t \leq -1\}$, the vector fields $V^{(1)},V^{(2)},V^{(3)}$ are $C(L)\varepsilon_1$-close to the standard rotation vector fields on the cylinder in the $C^1$-norm. Consequently, in the region $\{|z| \leq 1000, \, -1000 \leq t \leq -1\}$, the vector fields $V^{(1)},V^{(2)},V^{(3)}$ are $C(L)\varepsilon_1$-close to the standard rotation vector fields on the cylinder in the $C^{100}$-norm. 

\textit{Step 3:} We now define $h^{(a)}(t) := \mathscr{L}_{V^{(a)}(t)}(g(t))$. Since $\frac{\partial}{\partial t} V^{(a)} = \Delta V^{(a)}+\text{\rm Ric}(V^{(a)})$, we conclude that 
\[\frac{\partial}{\partial t} h^{(a)}(t) = \Delta_{L,g(t)} h^{(a)}(t)\] 
by Corollary \ref{pde.for.lie.derivative}. Using the estimate for $V^{(a)}-U^{(a)}$ in Step 2, we obtain  
\[|h^{(a)}| \leq |\mathscr{L}_{U^{(a)}}(g)| + C \, |D(V^{(a)}-U^{(a)})| \leq C\varepsilon\] 
in the region $\{|z| \leq \frac{3L}{4}, \, -\frac{3L}{4} \leq t \leq -\frac{L}{4}\}$, and 
\[|h^{(a)}| \leq |\mathscr{L}_{U^{(a)}}(g)| + C \, |D(V^{(a)}-U^{(a)})| \leq CL^2 \varepsilon\] 
in the region $\{|z| \leq \frac{3L}{4}, \, -\frac{L}{4} \leq t \leq -1\}$. Using standard interior estimates for linear parabolic equations, we deduce that 
\[\sum_{l=0}^{100} |D^l h^{(a)}| \leq C(L) \varepsilon\] 
in the region $\{|z| \leq \frac{L}{2}, \, -\frac{L}{2} \leq t \leq -1\}$.

\textit{Step 4:} By assumption, $\sum_{l=0}^{100} |D^l(\bar{g}(t)-g(t))| \leq C(L)\varepsilon_1$ in the region $\{|z| \leq \frac{L}{2}, \, -\frac{L}{2} \leq t \leq -1\}$, where $\bar{g}(t) = (-2t) \, g_{S^2} + dz \otimes dz$ denotes the standard metric on the shrinking cylinders. Let $\bar{h}^{(a)}$ denote the solution of the equation
\[\frac{\partial}{\partial t} \bar{h}^{(a)}(t) = \Delta_{L,\bar{g}(t)} \bar{h}^{(a)}(t)\] 
in the region $\{|z| \leq \frac{L}{2}, \, -\frac{L}{2} \leq t \leq -1\}$ with Dirichlet boundary condition $\bar{h}^{(a)}=h^{(a)}$ on the parabolic boundary $\{|z| \leq \frac{L}{2}, \, t=-\frac{L}{2}\} \cup \{|z|=\frac{L}{2}, \, -\frac{L}{2} \leq t \leq -1\}$. We compute 
\[\frac{\partial}{\partial t} (\bar{h}^{(a)}(t)-h^{(a)}(t)) - \Delta_{L,\bar{g}(t)} (\bar{h}^{(a)}(t)-h^{(a)}(t)) = E^{(a)}(t),\] 
where the error term $E^{(a)}(t)$ is defined by $E^{(a)}(t) := \Delta_{L,\bar{g}(t)} h^{(a)}(t) - \Delta_{L,g(t)} h^{(a)}(t)$. Using the estimates $\sum_{l=0}^{100} |D^l(\bar{g}(t)-g(t))| \leq C(L)\varepsilon_1$ and $\sum_{l=0}^{100} |D^l h^{(a)}| \leq C(L) \varepsilon$, we obtain $\sum_{l=0}^{90} |D^l E^{(a)}| \leq C(L) \varepsilon_1 \varepsilon$ in the region $\{|z| \leq \frac{L}{2}, \, -\frac{L}{2} \leq t \leq -1\}$. Using the maximum principle, we conclude that 
\[|\bar{h}^{(a)}-h^{(a)}| \leq C(L)\varepsilon_1\varepsilon\] 
in the region $\{|z| \leq \frac{L}{2}, \, -\frac{L}{2} \leq t \leq -1\}$. Standard interior estimates for linear parabolic equations imply 
\[\sum_{l=0}^{80} |D^l(\bar{h}^{(a)}-h^{(a)})| \leq C(L)\varepsilon_1\varepsilon\] 
in the region $\{|z| \leq 1000, \, -1000 \leq t \leq -1\}$.

\textit{Step 5:} We now apply Proposition \ref{parabolic.lichnerowicz.equation.on.cylinder} to $\bar{h}^{(a)}(t)$. Using the results in Step 3 and Step 4, we obtain 
\[|\bar{h}^{(a)}| \leq C\varepsilon+C(L)\varepsilon_1\varepsilon\] 
in the region $\{|z| \leq \frac{L}{2}, \, -\frac{L}{2} \leq t \leq -\frac{L}{4}\}$, and 
\[|\bar{h}^{(a)}| \leq CL^2 \varepsilon + C(L)\varepsilon_1\varepsilon\] 
in the region $\{|z| \leq \frac{L}{2}, \, -\frac{L}{4} \leq t \leq -1\}$. By Proposition \ref{parabolic.lichnerowicz.equation.on.cylinder}, we can find a function $\psi^{(a)}: S^2 \to \mathbb{R}$ (independent of $z$ and $t$) and rotationally invariant functions $\bar{\omega}^{(a)}(z,t)$ and $\bar{\beta}^{(a)}(z,t)$ with the following properties: 
\begin{itemize} 
\item $\psi^{(a)}$ lies in the span of the first spherical harmonics on $S^2$.
\item $\bar{\omega}^{(a)}(z,t)$ and $\bar{\beta}^{(a)}(z,t)$ are solutions of the one-dimensional heat equation.
\item $\big | \bar{h}^{(a)}(t) - \bar{\omega}^{(a)}(z,t) \, g_{S^2} - \bar{\beta}^{(a)}(z,t) \, dz \otimes dz - (-t) \, \psi^{(a)} \, g_{S^2} \big | \leq CL^{-\frac{1}{2}} \varepsilon + C(L)\varepsilon_1\varepsilon$ 
in the region $\{|z| \leq 1000, \, -1000 \leq t \leq -1\}$. 
\end{itemize} 
Note that $\bar{\omega}^{(a)}$, $\bar{\beta}^{(a)}$, and $\psi^{(a)}$ are bounded by $C(L)\varepsilon$. Moreover, the tensor $\bar{h}^{(a)}(t) - \bar{\omega}^{(a)}(z,t) \, g_{S^2} - \bar{\beta}^{(a)}(z,t) \, dz \otimes dz - (-t) \, \psi^{(a)} \, g_{S^2}$ satisfies the parabolic Lichnerowicz equation with respect to the background metrics $\bar{g}(t)$. Using standard interior estimates for linear parabolic equations, we obtain 
\begin{align*} 
&\sum_{l=0}^{80} \big | D^l(\bar{h}^{(a)}(t) -  \bar{\omega}^{(a)}(z,t) \, g_{S^2} - \bar{\beta}^{(a)}(z,t) \, dz \otimes dz - (-t) \, \psi^{(a)} \, g_{S^2}) \big | \\ 
&\leq CL^{-\frac{1}{2}} \varepsilon + C(L)\varepsilon_1\varepsilon 
\end{align*}
in the region $\{|z| \leq 800, \, -400 \leq t \leq -1\}$. Combining this estimate with the estimate in Step 4, we conclude that 
\begin{align*} 
&\sum_{l=0}^{80} \big | D^l(h^{(a)}(t) -  \bar{\omega}^{(a)}(z,t) \, g_{S^2} - \bar{\beta}^{(a)}(z,t) \, dz \otimes dz - (-t) \, \psi^{(a)} \, g_{S^2}) \big | \\ 
&\leq CL^{-\frac{1}{2}} \varepsilon + C(L)\varepsilon_1\varepsilon 
\end{align*}
in the region $\{|z| \leq 800, \, -400 \leq t \leq -1\}$.

\textit{Step 6:} Let us define a vector field $\xi^{(a)}$ on $S^2$ by $g_{S^2}(\xi^{(a)},\cdot) = -\frac{1}{4} \, d\psi^{(a)}$. Note that $\xi^{(a)}$ is independent of $z$ and $t$, and $|\xi^{(a)}| \leq C(L)\varepsilon$. Since $\psi^{(a)}$ lies in the span of the first spherical harmonic on $S^2$, we obtain $\mathscr{L}_{\xi^{(a)}}(g_{S^2}) = \frac{1}{2} \, \psi^{(a)} \, g_{S^2}$, hence $\mathscr{L}_{\xi^{(a)}}(\bar{g}(t)) = (-t) \, \psi^{(a)} \, g_{S^2}$. We now define $W^{(a)} := V^{(a)} - \xi^{(a)}$. In the region $\{|z| \leq 1000, \, -1000 \leq t \leq -1\}$, the vector fields $W^{(1)},W^{(2)},W^{(3)}$ are $C(L)\varepsilon_1$-close to the standard rotation vector fields on the cylinder in the $C^{80}$-norm. Using the identity  
\[\mathscr{L}_{W^{(a)}}(g(t)) = h^{(a)}(t) - (-t) \, \psi^{(a)} \, g_{S^2} + \mathscr{L}_{\xi^{(a)}}(\bar{g}(t)-g(t))\] 
and the estimates in Step 5, we obtain 
\begin{align*} 
&\sum_{l=0}^{60} |D^l(\mathscr{L}_{W^{(a)}}(g(t)) -  \bar{\omega}^{(a)}(z,t) \, g_{S^2} - \bar{\beta}^{(a)}(z,t) \, dz \otimes dz)| \\ 
&\leq CL^{-\frac{1}{2}} \varepsilon + C(L)\varepsilon_1\varepsilon 
\end{align*}
in the region $\{|z| \leq 800, \, -400 \leq t \leq -1\}$. 

We next estimate the time derivative of $W^{(a)}$. We compute 
\[\frac{\partial}{\partial t} W^{(a)} = \frac{\partial}{\partial t} V^{(a)} = \Delta V^{(a)} + \text{\rm Ric}(V^{(a)}) = \text{\rm div} \, h^{(a)} - \frac{1}{2} \, \nabla(\text{\rm tr} \, h^{(a)}).\] 
Using the estimates in Step 5, we obtain 
\[\sum_{l=0}^{60} \Big | D^l \Big ( \text{\rm div} \, h^{(a)} - \frac{\partial \bar{\beta}^{(a)}}{\partial z} \, \frac{\partial}{\partial z} - \frac{1}{(-t)} \, \xi^{(a)} \Big ) \Big | \leq CL^{-\frac{1}{2}} \varepsilon + C(L)\varepsilon_1\varepsilon\] 
and 
\begin{align*} 
&\sum_{l=0}^{60} \Big | D^l \Big ( \nabla(\text{\rm tr} \, h^{(a)}) - \Big ( \frac{\partial \bar{\beta}^{(a)}}{\partial z} + \frac{1}{(-t)} \, \frac{\partial \bar{\omega}^{(a)}}{\partial z} \Big ) \, \frac{\partial}{\partial z} - \frac{2}{(-t)} \, \xi^{(a)} \Big ) \Big | \\ 
&\leq CL^{-\frac{1}{2}} \varepsilon + C(L)\varepsilon_1\varepsilon 
\end{align*}
in the region $\{|z| \leq 800, \, -400 \leq t \leq -1\}$. Putting these facts together, we conclude that 
\[\sum_{l=0}^{60} \Big | D^l \Big (\frac{\partial}{\partial t} W^{(a)} - \Big ( \frac{1}{2} \, \frac{\partial \bar{\beta}^{(a)}}{\partial z} - \frac{1}{(-2t)} \, \frac{\partial \bar{\omega}^{(a)}}{\partial z} \Big ) \, \frac{\partial}{\partial z} \Big ) \Big | \leq CL^{-\frac{1}{2}} \varepsilon + C(L)\varepsilon_1\varepsilon\] 
in the region $\{|z| \leq 800, \, -400 \leq t \leq -1\}$. 

\textit{Step 7:} We now define $X^{(1)} := [W^{(2)},W^{(3)}]$, $X^{(2)} := [W^{(3)},W^{(1)}]$, $W^{(3)} := [W^{(1)},W^{(2)}]$. In the region $\{|z| \leq 800, \, -400 \leq t \leq -1\}$, the vector fields $X^{(1)},X^{(2)},X^{(3)}$ agree with the standard rotation vector fields on the cylinder up to constant factors and errors of order $C(L)\varepsilon_1$. Moreover,  
\begin{align*} 
\mathscr{L}_{X^{(1)}}(g)
&= \mathscr{L}_{W^{(2)}}(\mathscr{L}_{W^{(3)}}(g)) - \mathscr{L}_{W^{(3)}}(\mathscr{L}_{W^{(2)}}(g)) \\ 
&= \mathscr{L}_{W^{(2)}}(\mathscr{L}_{W^{(3)}}(g) -  \bar{\omega}^{(3)}(z,t) \, g_{S^2} - \bar{\beta}^{(3)}(z,t) \, dz \otimes dz) \\ 
&- \mathscr{L}_{W^{(3)}}(\mathscr{L}_{W^{(2)}}(g) - \bar{\omega}^{(2)}(z,t) \, g_{S^2} - \bar{\beta}^{(2)}(z,t) \, dz \otimes dz) \\ 
&+ \mathscr{L}_{W^{(2)}}(\bar{\omega}^{(3)}(z,t) \, g_{S^2} + \bar{\beta}^{(3)}(z,t) \, dz \otimes dz) \\ 
&- \mathscr{L}_{W^{(3)}}(\bar{\omega}^{(2)}(z,t) \, g_{S^2} + \bar{\beta}^{(2)}(z,t) \, dz \otimes dz). 
\end{align*} 
Analogous identities hold for $\mathscr{L}_{X^{(2)}}(g)$ and $\mathscr{L}_{X^{(3)}}(g)$. Since $W^{(1)},W^{(2)},W^{(3)}$ agree with the standard rotation vector fields on the cylinder up to errors of order $C(L)\varepsilon_1$, we obtain 
\[\sum_{l=0}^{40} |D^l(\mathscr{L}_{W^{(a)}}(\bar{\omega}^{(b)}(z,t) \, g_{S^2} + \bar{\beta}^{(b)}(z,t) \, dz \otimes dz))| \leq C(L)\varepsilon_1 \varepsilon\] 
in the region $\{|z| \leq 800, \, -400 \leq t \leq -1\}$. Combining this with the estimates in Step 6, we conclude that 
\[\sum_{l=0}^{40} |D^l (\mathscr{L}_{X^{(a)}}(g))| \leq CL^{-\frac{1}{2}} \varepsilon + C(L)\varepsilon_1\varepsilon\]
in the region $\{|z| \leq 800, \, -400 \leq t \leq -1\}$. 

We now estimate the time derivative of $X^{(a)}$. We compute 
\begin{align*} 
\frac{\partial}{\partial t} X^{(1)} 
&= \Big [ \frac{\partial}{\partial t} W^{(2)},W^{(3)} \Big ] + \Big [ W^{(2)},\frac{\partial}{\partial t} W^{(3)} \Big ] \\ 
&= \Big [ \frac{\partial}{\partial t} W^{(2)} - \Big ( \frac{1}{2} \, \frac{\partial \bar{\beta}^{(2)}}{\partial z} - \frac{1}{(-2t)} \, \frac{\partial \bar{\omega}^{(2)}}{\partial z} \Big ) \, \frac{\partial}{\partial z},W^{(3)} \Big ] \\ 
&+ \Big [ W^{(2)},\frac{\partial}{\partial t} W^{(3)} - \Big ( \frac{1}{2} \, \frac{\partial \bar{\beta}^{(3)}}{\partial z} - \frac{1}{(-2t)} \, \frac{\partial \bar{\omega}^{(3)}}{\partial z} \Big ) \, \frac{\partial}{\partial z} \Big ] \\ 
&+ \Big [ \Big ( \frac{1}{2} \, \frac{\partial \bar{\beta}^{(2)}}{\partial z} - \frac{1}{(-2t)} \, \frac{\partial \bar{\omega}^{(2)}}{\partial z} \Big ) \, \frac{\partial}{\partial z},W^{(3)} \Big ] \\ 
&+ \Big [ W^{(2)},\Big ( \frac{1}{2} \, \frac{\partial \bar{\beta}^{(3)}}{\partial z} - \frac{1}{(-2t)} \, \frac{\partial \bar{\omega}^{(3)}}{\partial z} \Big ) \, \frac{\partial}{\partial z} \Big ]. 
\end{align*} 
Analogous identities hold for $\frac{\partial}{\partial t} X^{(2)}$ and $\frac{\partial}{\partial t} X^{(3)}$. Since $W^{(1)},W^{(2)},W^{(3)}$ agree with the standard rotation vector fields on the cylinder up to errors of order $C(L)\varepsilon_1$, we obtain 
\[\sum_{l=0}^{40} \Big| D^l \Big [ W^{(a)},\Big ( \frac{1}{2} \, \frac{\partial \bar{\beta}^{(b)}}{\partial z} - \frac{1}{(-2t)} \, \frac{\partial \bar{\omega}^{(b)}}{\partial z} \Big ) \, \frac{\partial}{\partial z} \Big ] \Big | \leq C(L)\varepsilon_1\varepsilon\]
in the region $\{|z| \leq 800, \, -400 \leq t \leq -1\}$. Combining this with the estimates in Step 6, we conclude that 
\[\sum_{l=0}^{40} \Big | D^l \Big ( \frac{\partial}{\partial t} X^{(a)} \Big ) \Big | \leq CL^{-\frac{1}{2}} \varepsilon + C(L)\varepsilon_1\varepsilon\]
in the region $\{|z| \leq 800, \, -400 \leq t \leq -1\}$. 

\textit{Step 8:} Let $Y^{(a)}$ be a time-independent vector field such that $Y^{(a)}=X^{(a)}$ at time $-1$. In the region $\{|z| \leq 800, \, -400 \leq t \leq -1\}$, the vector fields $Y^{(1)},Y^{(2)},Y^{(3)}$ agree with the standard rotation vector fields on the cylinder up to constant factors and errors of order $C(L)\varepsilon_1$. The estimates for $\frac{\partial}{\partial t} X^{(a)}$ in Step 7 imply 
\[\sum_{l=0}^{40} |D^l(Y^{(a)}-X^{(a)})| \leq CL^{-\frac{1}{2}} \varepsilon + C(L)\varepsilon_1\varepsilon\]
in the region $\{|z| \leq 800, \, -400 \leq t \leq -1\}$. Using the estimates for $\mathscr{L}_{X^{(a)}}(g)$ in Step 6, we obtain 
\[\sum_{l=0}^{30} |D^l (\mathscr{L}_{Y^{(a)}}(g))| \leq CL^{-\frac{1}{2}} \varepsilon + C(L)\varepsilon_1\varepsilon\]
in the region $\{|z| \leq 800, \, -400 \leq t \leq -1\}$. 

\textit{Step 9:} In the following, we fix a time $t \in [-200,-1]$. We denote by $\Sigma_s$ the leaves of the CMC foliation of $(M,g(t))$. Note that the foliation depends on $t$, but we suppress this dependence in the notation. Let $\nu$ denote the unit normal vector field to the foliation $\Sigma_s$. For each $s$, we denote by $v: \Sigma_s \to \mathbb{R}$ the lapse function associated with this foliation. We assume that the foliation $\Sigma_s$ is parametrized so that $x_0 \in \Sigma_0$ and $\int_{\Sigma_s} v = 1$ for all $s$. Since $\Sigma_s$ is a CMC surface for each $s$, the function $v$ satisfies the Jacobi equation
\[\Delta_{\Sigma_s} v + (|A|^2+\text{\rm Ric}(\nu,\nu)) \, v = \text{\rm constant}\] 
on $\Sigma_s$, where $|A|$ denotes the norm of the second fundamental form of $\Sigma_s$ in $(M,g(t))$. The Jacobi operator $\Delta_{\Sigma_s} +  (|A|^2+\text{\rm Ric}(\nu,\nu))$ is a small perturbation of the Laplacian $\Delta_{\Sigma_s}$. Hence, for each $s$, the Jacobi operator $\Delta_{\Sigma_s} + (|A|^2+\text{\rm Ric}(\nu,\nu))$ is an invertible operator from the space $\{f \in C^{2,\frac{1}{2}}(\Sigma_s): \int_{\Sigma_s} f = 0\}$ to the space $\{f \in C^{\frac{1}{2}}(\Sigma_s): \int_{\Sigma_s} f v = 0\}$, and we have a uniform bound for the norm of its inverse.

In the following, we only consider those leaves of the foliation $\Sigma_s$ which are contained in the region $\{|z| \leq 700\}$. Let us define a function $F^{(a)}: \Sigma_s \to \mathbb{R}$ by $F^{(a)} := \langle Y^{(a)},\nu \rangle$. The quantity 
\[\Delta_{\Sigma_s} F^{(a)} +  (|A|^2+\text{\rm Ric}(\nu,\nu)) \, F^{(a)} =: H^{(a)}\] 
can be expressed in terms of $\mathscr{L}_{Y^{(a)}}(g)$ and the first derivatives of $\mathscr{L}_{Y^{(a)}}(g)$. Using the estimate for $\mathscr{L}_{Y^{(a)}}(g)$ in Step 8, we deduce that $\sum_{l=0}^{20} |D^l H^{(a)}| \leq CL^{-\frac{1}{2}} \varepsilon + C(L)\varepsilon_1\varepsilon$ in the region $\{|z| \leq 700\}$. We next define $G^{(a)}(s) := \int_{\Sigma_s} F^{(a)}$ and $\tilde{F}^{(a)} := F^{(a)} - G^{(a)}(s) \, v$. Then $\int_{\Sigma_s} \tilde{F}^{(a)} = 0$ and 
\[\Delta_{\Sigma_s} \tilde{F}^{(a)} +  (|A|^2+\text{\rm Ric}(\nu,\nu)) \, \tilde{F}^{(a)} = H^{(a)} - \int_{\Sigma_s} H^{(a)} v\] 
on $\Sigma_s$. Using the estimate $\sum_{l=0}^{20} |D^l H^{(a)}| \leq CL^{-\frac{1}{2}} \varepsilon + C(L)\varepsilon_1\varepsilon$, we conclude that $\sum_{l=0}^{10} |D^l \tilde{F}^{(a)}| \leq CL^{-\frac{1}{2}} \varepsilon + C(L)\varepsilon_1\varepsilon$ in the region $\{|z| \leq 600\}$. Since $v^{-1} \, \tilde{F}^{(a)} = v^{-1} \, \langle Y^{(a)},\nu \rangle - G^{(a)}(s)$, it follows that 
\[\sum_{l=0}^{10} |D^l (v^{-1} \, \langle Y^{(a)},\nu \rangle - G^{(a)}(s))| \leq CL^{-\frac{1}{2}} \varepsilon + C(L)\varepsilon_1\varepsilon\] 
in the region $\{|z| \leq 600\}$.

By the divergence theorem, the quantity $G^{(a)}(s) - G^{(a)}(0) = \int_{\Sigma_s} \langle Y^{(a)},\nu \rangle - \int_{\Sigma_0} \langle Y^{(a)},\nu \rangle$ can be expressed as an integral of $\text{\rm div} \, Y^{(a)}$ over the region bounded by $\Sigma_0$ and $\Sigma_s$. Differentiating this identity with respect to $s$ gives 
\[\frac{d}{ds} G^{(a)}(s) = \int_{\Sigma_s} v \, \text{\rm div} \, Y^{(a)}.\] 
Using the estimate for $\mathscr{L}_{Y^{(a)}}(g)$ in Step 8, we obtain 
\[\sum_{l=1}^{10} \big | \frac{d^l}{ds^l} G^{(a)}(s) \big | \leq CL^{-\frac{1}{2}} \varepsilon + C(L)\varepsilon_1\varepsilon.\] 
Putting these facts together, we conclude that 
\[\sum_{l=1}^{10} |D^l (v^{-1} \, \langle Y^{(a)},\nu \rangle)| \leq CL^{-\frac{1}{2}} \varepsilon + C(L)\varepsilon_1\varepsilon\] 
in the region $\{|z| \leq 600\}$.

\textit{Step 10:} Finally, we define $Z^{(1)} := [Y^{(2)},Y^{(3)}]$, $Z^{(2)} := [Y^{(3)},Y^{(1)}]$, $Z^{(3)} := [Y^{(1)},Y^{(2)}]$. Clearly, $Z^{(1)},Z^{(2)},Z^{(3)}$ are time-independent vector fields. In the region $\{|z| \leq 800, \, -400 \leq t \leq -1\}$, the vector fields $Z^{(1)},Z^{(2)},Z^{(3)}$ agree with the standard rotation vector fields on the cylinder up to constant factors and errors of order $C(L)\varepsilon_1$. Note that 
\[\mathscr{L}_{Z^{(1)}}(g) = \mathscr{L}_{Y^{(2)}}(\mathscr{L}_{Y^{(3)}}(g)) - \mathscr{L}_{Y^{(3)}}(\mathscr{L}_{Y^{(2)}}(g)).\] 
We have shown in Step 8 that $\sum_{l=0}^{30} |D^l (\mathscr{L}_{Y^{(a)}}(g))| \leq CL^{-\frac{1}{2}} \varepsilon + C(L)\varepsilon_1\varepsilon$ in the region $\{|z| \leq 500, \,  t \in [-200,-1]\}$. This gives 
\[\sum_{l=0}^{20} |D^l (\mathscr{L}_{Z^{(a)}}(g))| \leq CL^{-\frac{1}{2}} \varepsilon + C(L)\varepsilon_1\varepsilon\]
in the region $\{|z| \leq 500, \, -200 \leq t \leq -1\}$. Now, let us fix a time $t \in [-200,-1]$, and let $\nu$ and $v$ denote the normal vector field and the lapse function of the CMC foliation at time $t$. Since the vector field $T := v^{-1} \, \nu$ is a gradient vector field, we have 
\[\langle Z^{(1)},T \rangle = \big \langle Y^{(2)},\nabla(\langle Y^{(3)},T \rangle) \big \rangle - \big \langle Y^{(3)},\nabla(\langle Y^{(2)},T \rangle) \big \rangle.\] 
Using the estimates in Step 9, we obtain 
\[\sum_{l=1}^{10} |D^l(\langle Y^{(a)},T \rangle)| \leq CL^{-\frac{1}{2}} \varepsilon + C(L)\varepsilon_1\varepsilon\] 
in the region $\{|z| \leq 500, \, -200 \leq t \leq -1\}$. Consequently, 
\[\sum_{l=0}^8 |D^l (\langle Z^{(a)},T \rangle)| \leq CL^{-\frac{1}{2}} \varepsilon + C(L)\varepsilon_1\varepsilon\] 
in the region $\{|z| \leq 500, \, -200 \leq t \leq -1\}$. Using again the fact that $T$ is a gradient vector field, we compute
\[\langle \nabla(|T|^2),Z^{(a)} \rangle = -(\mathscr{L}_{Z^{(a)}}(g))(T,T) + 2 \, \langle \nabla(\langle Z^{(a)},T \rangle),T \rangle\] 
and 
\[g([T,Z^{(a)}],\cdot) = (\mathscr{L}_{Z^{(a)}}(g))(T,\cdot) - d(\langle Z^{(a)},T \rangle).\] 
Using our estimates for $\mathscr{L}_{Z^{(a)}}(g)$ and $\langle Z^{(a)},T \rangle$, we finally obtain 
\[\sum_{l=0}^6 |D^l(\langle \nabla(|T|^2),Z^{(a)} \rangle)| \leq CL^{-\frac{1}{2}} \varepsilon + C(L)\varepsilon_1\varepsilon\] 
and 
\[\sum_{l=0}^6 |D^l([T,Z^{(a)}])| \leq CL^{-\frac{1}{2}} \varepsilon + C(L)\varepsilon_1\varepsilon\] 
in the region $\{|z| \leq 500, \, -200 \leq t \leq -1\}$. 

To summarize, we have shown that 
\[\sum_{l=0}^8 |D^l(\langle Z^{(a)},\nu \rangle)| \leq CL^{-\frac{1}{2}} \varepsilon + C(L)\varepsilon_1\varepsilon,\] 
\[\sum_{l=0}^6 |D^l(\langle \nabla v,Z^{(a)} \rangle)| \leq CL^{-\frac{1}{2}} \varepsilon + C(L)\varepsilon_1\varepsilon,\] 
\[\sum_{l=0}^6 |D^l([\nu,Z^{(a)}])| \leq CL^{-\frac{1}{2}} \varepsilon + C(L)\varepsilon_1\varepsilon,\] 
\[\sum_{l=0}^6 |D^l([v \, \nu,Z^{(a)}])| \leq CL^{-\frac{1}{2}} \varepsilon + C(L)\varepsilon_1\varepsilon\] 
in the region $\{|z| \leq 500, \, -200 \leq t \leq -1\}$. In particular, if $t \in [-200,-1]$ and $\Sigma \subset \{|z| \leq 400\}$ is a leaf of the CMC foliation in $(M,g(t))$, then the lapse function $v$ satisfies 
\[\sup_\Sigma \big | v - \text{\rm area}_{g(t)}(\Sigma)^{-1} \big | \leq CL^{-\frac{1}{2}} \varepsilon + C(L)\varepsilon_1\varepsilon.\] 

\textit{Step 11:} In the next step, we obtain information on the Ricci tensor and the second fundamental form of the CMC foliation. To that end, let us consider an arbitrary point $(\bar{x},\bar{t})$ in the region $\{|z| \leq 400, \, -200 \leq t \leq -1\}$. Let $\{e_1,e_2\}$ denote an orthonormal basis for the tangent space to the CMC foliation at $(\bar{x},\bar{t})$. Since the vector fields $Z^{(1)},Z^{(2)},Z^{(3)}$ are close to the standard rotation vector fields on the cylinder up to some constant factor, we can find a vector $\lambda = (\lambda_1,\lambda_2,\lambda_3) \in \mathbb{R}^3$ such that $\sum_{a=1}^3 \lambda_a \, \langle Z^{(a)},e_1 \rangle = \sum_{a=1}^3 \lambda_a \, \langle Z^{(a)},e_2 \rangle = 0$ at the point $(\bar{x},\bar{t})$ and 
\[\sum_{a=1}^3 \lambda_a \, \langle D_{e_1} Z^{(a)},e_2 \rangle = 1\] 
at the point $(\bar{x},\bar{t})$. Note that $|\lambda| \leq C$.

Using the estimate for $\langle Z^{(a)},\nu \rangle$ in Step 10, we obtain $\big | \sum_{a=1}^3 \lambda_a \, \langle Z^{(a)},\nu \rangle \big | \leq CL^{-\frac{1}{2}} \varepsilon + C(L)\varepsilon_1\varepsilon$. Consequently, 
\[\Big | \sum_{a=1}^3 \lambda_a \, Z^{(a)} \Big | \leq CL^{-\frac{1}{2}} \varepsilon + C(L)\varepsilon_1\varepsilon\] 
at the point $(\bar{x},\bar{t})$. Using the estimate for $\mathscr{L}_{Z^{(a)}}(g)$ in Step 10, we obtain $|\langle D_{e_1} Z^{(a)},e_1 \rangle| + |\langle D_{e_2} Z^{(a)},e_2 \rangle| \leq CL^{-\frac{1}{2}} \varepsilon + C(L)\varepsilon_1\varepsilon$ and $|\langle D_{e_1} Z^{(a)},e_2 \rangle + \langle D_{e_2} Z^{(a)},e_1 \rangle| \leq CL^{-\frac{1}{2}} \varepsilon + C(L)\varepsilon_1\varepsilon$. This implies 
\[\Big | 1 + \sum_{a=1}^3 \lambda_a \, \langle D_{e_2} Z^{(a)},e_1 \rangle \Big | \leq CL^{-\frac{1}{2}} \varepsilon + C(L)\varepsilon_1\varepsilon\] 
at the point $(\bar{x},\bar{t})$. Moreover, the estimate for the derivatives of $\langle Z^{(a)},\nu \rangle$ in Step 10 gives $|\langle D_{e_i} Z^{(a)},\nu \rangle + \langle Z^{(a)},D_{e_i} \nu \rangle| \leq CL^{-\frac{1}{2}} \varepsilon + C(L)\varepsilon_1\varepsilon$ for $i \in \{1,2\}$. Hence, for each $i \in \{1,2\}$, we obtain 
\[\Big | \sum_{a=1}^3 \lambda_a \, \langle D_{e_i} Z^{(a)},\nu \rangle \Big | \leq CL^{-\frac{1}{2}} \varepsilon + C(L)\varepsilon_1\varepsilon\] 
at the point $(\bar{x},\bar{t})$. 

In view of the estimates in Step 10, the Ricci tensor satisfies $|\mathscr{L}_{Z^{(a)}} \text{\rm Ric}| \leq CL^{-\frac{1}{2}} \varepsilon + C(L)\varepsilon_1\varepsilon$ for each $a \in \{1,2,3\}$. A straightforward calculation gives 
\begin{align*} 
(\mathscr{L}_{Z^{(a)}} \text{\rm Ric})(e_1,e_2) 
&= (D_{Z^{(a)}} \text{\rm Ric})(e_1,e_2) + \text{\rm Ric}(D_{e_1} Z^{(a)},e_2) + \text{\rm Ric}(e_1,D_{e_2} Z^{(a)}) \\ 
&= (D_{Z^{(a)}} \text{\rm Ric})(e_1,e_2) + (\langle D_{e_1} Z^{(a)},e_1 \rangle+\langle D_{e_2} Z^{(a)},e_2 \rangle) \, \text{\rm Ric}(e_1,e_2) \\ 
&+ \langle D_{e_1} Z^{(a)},e_2 \rangle \, \text{\rm Ric}(e_2,e_2) + \langle D_{e_2} Z^{(a)},e_1 \rangle \, \text{\rm Ric}(e_1,e_1) \\ 
&+ \langle D_{e_1} Z^{(a)},\nu \rangle \, \text{\rm Ric}(\nu,e_2) + \langle D_{e_2} Z^{(a)},\nu \rangle \, \text{\rm Ric}(e_1,\nu). 
\end{align*} 
If we multiply this identity by $\lambda_a$ and sum over $a \in \{1,2,3\}$, we conclude that 
\[|\text{\rm Ric}(e_2,e_2)-\text{\rm Ric}(e_1,e_1)| \leq CL^{-\frac{1}{2}} \varepsilon + C(L)\varepsilon_1\varepsilon\] 
at the point $(\bar{x},\bar{t})$. Therefore, $\big | \text{\rm Ric}(e_i,e_j) - \frac{1}{2} \, \text{\rm tr}_\Sigma(\text{\rm Ric}) \, \delta_{ij} \big | \leq CL^{-\frac{1}{2}} \varepsilon + C(L)\varepsilon_1\varepsilon$ at the point $(\bar{x},\bar{t})$, where $\text{\rm tr}_\Sigma(\text{\rm Ric}) = \text{\rm Ric}(e_1,e_1)+\text{\rm Ric}(e_2,e_2)$.

Let $A$ denote the second fundamental form of the CMC foliation. We can think of $A$ as a $(0,2)$-tensor on $M$, which vanishes in the normal direction. The estimates in Step 10 imply $|\mathscr{L}_{Z^{(a)}} A| \leq CL^{-\frac{1}{2}} \varepsilon + C(L)\varepsilon_1\varepsilon$ for each $a \in \{1,2,3\}$. A straightforward calculation gives 
\begin{align*} 
(\mathscr{L}_{Z^{(a)}} A)(e_1,e_2) 
&= (D_{Z^{(a)}} A)(e_1,e_2) + A(D_{e_1} Z^{(a)},e_2) + A(e_1,D_{e_2} Z^{(a)}) \\ 
&= (D_{Z^{(a)}} A)(e_1,e_2) + (\langle D_{e_1} Z^{(a)},e_1 \rangle+\langle D_{e_2} Z^{(a)},e_2 \rangle) \, A(e_1,e_2) \\ 
&+ \langle D_{e_1} Z^{(a)},e_2 \rangle \, A(e_2,e_2) + \langle D_{e_2} Z^{(a)},e_1 \rangle \, A(e_1,e_1).
\end{align*} 
If we multiply this identity by $\lambda_a$ and sum over $a \in \{1,2,3\}$, we conclude that 
\[|A(e_2,e_2)-A(e_1,e_1)| \leq CL^{-\frac{1}{2}} \varepsilon + C(L)\varepsilon_1\varepsilon\] 
at the point $(\bar{x},\bar{t})$. Therefore, $\big | A(e_i,e_j) - \frac{1}{2} \, H \, \delta_{ij} \big | \leq CL^{-\frac{1}{2}} \varepsilon + C(L)\varepsilon_1\varepsilon$ at the point $(\bar{x},\bar{t})$, where $H$ denotes the mean curvature of the CMC foliation.

Finally, the estimates in Step 10 imply 
\[\inf_\rho \sup_\Sigma \big | \frac{1}{2} \, \text{\rm tr}_\Sigma(\text{\rm Ric}) - \rho \big | \leq CL^{-\frac{1}{2}} \varepsilon + C(L)\varepsilon_1\varepsilon\] 
if $t \in [-200,-1]$ and $\Sigma \subset \{|z| \leq 400\}$ is a leaf of the CMC foliation in $(M,g(t))$. To summarize, if $t \in [-200,-1]$ and $\Sigma \subset \{|z| \leq 400\}$ is a leaf of the CMC foliation in $(M,g(t))$, then 
\[\inf_\rho \sup_\Sigma \big | (\text{\rm Ric}-\rho g)|_{T\Sigma} \big | \leq CL^{-\frac{1}{2}} \varepsilon + C(L)\varepsilon_1\varepsilon\] 
and 
\[\sup_\Sigma \big | (A-\frac{1}{2} \, Hg)_{T\Sigma} \big | \leq CL^{-\frac{1}{2}} \varepsilon + C(L)\varepsilon_1\varepsilon,\] 
where $H$ denotes the mean curvature of $\Sigma$ (which is constant). 

\textit{Step 12:} Let us fix a time $t \in [-200,-1]$. By Step 10, the vector fields $Z^{(1)},Z^{(2)},Z^{(3)}$ are tangential to the CMC foliation of $(M,g(t))$, up to errors of order $CL^{-\frac{1}{2}} \varepsilon + C(L)\varepsilon_1\varepsilon$. Moreover, the vector fields $Z^{(1)},Z^{(2)},Z^{(3)}$ are tangential to the CMC foliation of $(M,g(-1))$, up to errors of order $CL^{-\frac{1}{2}} \varepsilon + C(L)\varepsilon_1\varepsilon$. Since the vector fields $Z^{(1)},Z^{(2)},Z^{(3)}$ are close to the standard rotation vector fields on the cylinder, we conclude that every leaf of the CMC foliation of $(M,g(t))$ which is contained in the region $\{|z| \leq 400\}$ is $(CL^{-\frac{1}{2}}+C(L)\varepsilon_1\varepsilon)$-close in the $C^1$-norm to a leaf of the CMC foliation of $(M,g(-1))$. 

\textit{Step 13:} We again fix a time $t \in [-200,-1]$. Let $\Sigma_s$ denote the CMC foliation of $(M,g(t))$, and let $\nu$ and $v$ denote the normal vector field and the lapse function associated with this foliation. In the following, we only consider those leaves of the foliation which are contained in the region $\{|z| \leq 300\}$. Our goal is to show that the quantity $\text{\rm area}_{g(t)}(\Sigma_s)^{-2} \int_{\Sigma_s} \langle Z^{(a)},Z^{(b)} \rangle_{g(t)} \, d\mu_{g(t)}$ is nearly constant in $s$, up to errors of order $CL^{-\frac{1}{2}} \varepsilon + C(L)\varepsilon_1\varepsilon$. Recall that the surfaces $\Sigma_s$ move with normal velocity $v$. This implies $\frac{d}{ds} \text{\rm area}_{g(t)}(\Sigma_s) = \int_{\Sigma_s} Hv \, d\mu_{g(t)} = H$, where $H$ denotes the mean curvature of $\Sigma_s$ with respect to the metric $g(t)$. We next compute 
\begin{align*} 
&\frac{d}{ds} \bigg ( \int_{\Sigma_s} \langle Z^{(a)},Z^{(b)} \rangle_{g(t)} \, d\mu_{g(t)} \bigg ) \\ 
&= \int_{\Sigma_s} Hv \, \langle Z^{(a)},Z^{(b)} \rangle_{g(t)} \, d\mu_{g(t)} + \int_{\Sigma_s} (\mathscr{L}_{v \, \nu}(g))(Z^{(a)},Z^{(b)})  \, d\mu_{g(t)} \\ 
&+ \int_{\Sigma_s} \langle [v \, \nu,Z^{(a)}],Z^{(b)} \rangle_{g(t)} \, d\mu_{g(t)} + \int_{\Sigma_s} \langle Z^{(a)},[v \, \nu,Z^{(b)}] \rangle_{g(t)} \, d\mu_{g(t)}. 
\end{align*} 
The estimate $\big | (A-\frac{1}{2} \, Hg)|_{T\Sigma} \big | \leq CL^{-\frac{1}{2}} \varepsilon + C(L)\varepsilon_1\varepsilon$ in Step 11 implies 
\[\big | (\mathscr{L}_{v \, \nu}(g))(Z^{(a)},Z^{(b)}) - Hv \, \langle Z^{(a)},Z^{(b)} \rangle_{g(t)} \big | \leq CL^{-\frac{1}{2}} \varepsilon + C(L)\varepsilon_1\varepsilon.\] 
Moreover, the estimate for $[v \, \nu,Z^{(a)}]$ in Step 10 gives 
\[|\langle [v \, \nu,Z^{(a)}],Z^{(b)} \rangle_{g(t)}| \leq CL^{-\frac{1}{2}} \varepsilon + C(L)\varepsilon_1\varepsilon.\] 
An analogous argument yields 
\[|\langle Z^{(a)},[v \, \nu,Z^{(b)}] \rangle_{g(t)}| \leq CL^{-\frac{1}{2}} \varepsilon + C(L)\varepsilon_1\varepsilon.\] 
Putting these facts together, we obtain 
\begin{align*} 
&\bigg | \frac{d}{ds} \bigg ( \int_{\Sigma_s} \langle Z^{(a)},Z^{(b)} \rangle_{g(t)} \, d\mu_{g(t)} \bigg ) - 2 \int_{\Sigma_s} Hv \, \langle Z^{(a)},Z^{(b)} \rangle_{g(t)} \, d\mu_{g(t)} \bigg | \\ 
&\leq (CL^{-\frac{1}{2}} \varepsilon + C(L)\varepsilon_1\varepsilon), 
\end{align*} 
hence 
\begin{align*} 
&\Big | \frac{d}{ds} \bigg ( \int_{\Sigma_s} \langle Z^{(a)},Z^{(b)} \rangle_{g(t)} \, d\mu_{g(t)} \bigg ) - \frac{2H}{\text{\rm area}_{g(t)}(\Sigma_s)} \int_{\Sigma_s} \langle Z^{(a)},Z^{(b)} \rangle_{g(t)} \, d\mu_{g(t)} \Big | \\ 
&\leq (CL^{-\frac{1}{2}} \varepsilon + C(L)\varepsilon_1\varepsilon). 
\end{align*}
Thus, we conclude that  
\[\bigg | \frac{d}{ds} \bigg ( \text{\rm area}_{g(t)}(\Sigma_s)^{-2} \int_{\Sigma_s} \langle Z^{(a)},Z^{(b)} \rangle_{g(t)} \, d\mu_{g(t)} \bigg ) \bigg | \leq CL^{-\frac{1}{2}} \varepsilon + C(L)\varepsilon_1\varepsilon.\] 

\textit{Step 14:} The estimate in Step 13 implies that there exists a symmetric $3 \times 3$ matrix $Q_{ab}$ (independent of $\Sigma$) such that 
\[\bigg | Q_{ab} - \text{\rm area}_{g(-1)}(\Sigma)^{-2} \int_\Sigma \langle Z^{(a)},Z^{(b)} \rangle_{g(-1)} \, d\mu_{g(-1)} \bigg | \leq CL^{-\frac{1}{2}} \varepsilon + C(L)\varepsilon_1\varepsilon\]
whenever $\Sigma \subset \{|z| \leq 300\}$ is a leaf of the CMC foliation of $(M,g(-1))$. Moreover, since the vector fields $Z^{(1)},Z^{(2)},Z^{(3)}$ are close to the standard rotation vector fields on the cylinder up to some constant factor, the eigenvalues of the matrix $Q_{ab}$ lie in the interval $[\frac{1}{C},C]$ for some fixed constant $C$. The estimate for the Ricci tensor in Step 11 gives 
\[\bigg | \frac{d}{dt} \bigg ( \text{\rm area}_{g(t)}(\Sigma)^{-2} \int_\Sigma \langle Z^{(a)},Z^{(b)} \rangle_{g(t)} \, d\mu_{g(t)} \bigg ) \bigg | \leq CL^{-\frac{1}{2}} \varepsilon + C(L)\varepsilon_1\varepsilon\] 
whenever $t \in [-200,-1]$ and $\Sigma \subset \{|z| \leq 300\}$ is a fixed leaf of the CMC foliation of $(M,g(-1))$. Consequently, 
\[\bigg | Q_{ab} - \text{\rm area}_{g(t)}(\Sigma)^{-2} \int_\Sigma \langle Z^{(a)},Z^{(b)} \rangle_{g(t)} \, d\mu_{g(t)} \bigg | \leq CL^{-\frac{1}{2}} \varepsilon + C(L)\varepsilon_1\varepsilon\]
whenever $t \in [-200,-1]$ and $\Sigma \subset \{|z| \leq 300\}$ is a leaf of the CMC foliation of $(M,g(-1))$. 

By Step 12, every leaf of the CMC foliation of $(M,g(t))$ which is contained in the region $\{|z| \leq 200\}$ is $(CL^{-\frac{1}{2}}+C(L)\varepsilon_1\varepsilon)$-close in the $C^1$-norm to a leaf of the CMC foliation of $(M,g(-1))$. This finally implies 
\[\bigg | Q_{ab} - \text{\rm area}_{g(t)}(\Sigma)^{-2} \int_\Sigma \langle Z^{(a)},Z^{(b)} \rangle_{g(t)} \, d\mu_{g(t)} \bigg | \leq CL^{-\frac{1}{2}} \varepsilon + C(L)\varepsilon_1\varepsilon\]
whenever $t \in [-200,-1]$ and $\Sigma \subset \{|z| \leq 200\}$ is a leaf of the CMC foliation of $(M,g(t))$. Note that the matrix $Q_{ab}$ is independent of $t$ and independent of $\Sigma$.

By considering the vector fields $\sum_{b=1}^3 (Q^{-\frac{1}{2}})_{ab} \, Z^{(b)}$, we see that the point $(x_0,-1)$ is $(CL^{-\frac{1}{2}} \varepsilon + C(L)\varepsilon_1\varepsilon)$-symmetric. Hence, if we choose $L$ sufficiently large and $\varepsilon_1$ sufficiently small (depending on $L$), then $(x_0,-1)$ is $\frac{\varepsilon}{2}$-symmetric. This completes the proof of Theorem \ref{neck.improvement.theorem}. \\

\section{Rotational symmetry of ancient $\kappa$-solutions in dimension $3$}

\label{rot.sym}

In this section, we give the proof of Theorem \ref{main.thm.b}. Throughout this section, we assume that $(M,g(t))$, $t \in (-\infty,0]$, is a three-dimensional ancient $\kappa$-solution which is noncompact and has positive sectional curvature. Our goal is to show that $(M,g(t))$ is rotationally symmetric. For each $t$, we denote by $R_{\text{\rm max}}(t)$ the supremum of the scalar curvature of $(M,g(t))$. By Perelman's pointwise derivative estimate \cite{Perelman1}, the function $t \mapsto R_{\text{\rm max}}(t)^{-1}$ is uniformly Lipschitz continuous.

Let us fix a large constant $L$ and a small constant $\varepsilon_1$ such that the conclusion of the Neck Improvement Theorem holds. We assume that $\varepsilon_1$ is chosen small enough so that the results in Section \ref{approximate.killing.vector.fields} can be applied on every $\varepsilon_1$-neck. For each point $(x,t)$ in spacetime, we denote by $\lambda_1(x,t)$ the smallest eigenvalue of the Ricci tensor at $(x,t)$. The following is a direct consequence of Perelman's work:

\begin{proposition}
\label{choice.of.theta}
Given $\varepsilon_1$, we can find a small positive constant $\theta$ (depending on $\varepsilon_1$) with the following property. Suppose that $(\bar{x},\bar{t})$ is a point in spacetime satisfying $\lambda_1(\bar{x},\bar{t}) \leq \theta R(\bar{x},\bar{t})$. Then $(\bar{x},\bar{t})$ lies at the center of an evolving $\varepsilon_1$-neck. Moreover, if $x$ lies outside the compact domain bounded by the leaf of the CMC foliation passing through $(\bar{x},\bar{t})$, then $(x,\bar{t})$ lies at the center of an evolving $\varepsilon_1$-neck.
\end{proposition}

\textbf{Proof.} 
By Theorem \ref{canonical.neighborhood.theorem} and Corollary \ref{bound.for.R_max}, we can find a domain $\Omega_{\bar{t}}$ with the following properties: 
\begin{itemize}
\item If $x \in M \setminus \Omega_{\bar{t}}$, then $(x,\bar{t})$ lies at the center of an evolving $\varepsilon_1$-neck.
\item $\text{\rm diam}_{g(\bar{t})}(\Omega_{\bar{t}}) \leq C_0 \, R_{\text{\rm max}}(\bar{t})^{-\frac{1}{2}}$, where $C_0$ is a large constant that depends on $\varepsilon_1$.
\end{itemize}
Now, if we choose $\theta$ sufficiently small, then every point $(\bar{x},\bar{t})$ satisfying $\lambda_1(\bar{x},\bar{t}) \leq \theta R(\bar{x},\bar{t})$ lies at the center of a neck $N$ of length $10 C_0 \, R(\bar{x},\bar{t})^{-\frac{1}{2}}$, and furthermore every point on $N$ lies at the center of an evolving $\varepsilon_1$-neck. The complement $M \setminus N$ has two connected components, one of which is bounded and one of which is unbounded. The bounded connected component of $M \setminus N$ must contain a point which does not lie at the center of an evolving $\varepsilon_1$-neck. If the unbounded connected component of $M \setminus N$ also contains a point which does not lie at the center of an evolving $\varepsilon_1$-neck, then $\text{\rm diam}_{g(\bar{t})}(\Omega_{\bar{t}}) \geq 4C_0 \, R(\bar{x},\bar{t})^{-\frac{1}{2}} \geq 4C_0 \, R_{\text{\rm max}}(\bar{t})^{-\frac{1}{2}}$, which is a contradiction. Consequently, every point in the unbounded connected component of $M \setminus N$ must lie at the center of an evolving $\varepsilon_1$-neck. From this, the assertion follows easily. \\


In the following, we fix $\theta$ so that the conclusion of Proposition \ref{choice.of.theta} holds.

\begin{definition}
\label{symmetry.of.cap}
We say that the flow is $\varepsilon$-symmetric at time $\bar{t}$ if there exist a compact domain $D \subset M$ and time-independent vector fields $U^{(1)},U^{(2)},U^{(3)}$ which are defined on an open set containing $D$ such that the following statements hold: 
\begin{itemize} 
\item There exists a point $x \in \partial D$ such that $\lambda_1(x,\bar{t}) < \theta R(x,\bar{t})$.
\item For each $x \in D$, we have $\lambda_1(x,\bar{t}) > \frac{1}{2} \theta R(x,\bar{t})$.
\item The boundary $\partial D$ is a leaf of the CMC foliation of $(M,g(\bar{t}))$.
\item For each $x \in M \setminus D$, the point $(x,\bar{t})$ is $\varepsilon$-symmetric in the sense of Definition \ref{symmetry.of.necks}.
\item $\sup_{D \times [\bar{t}-R_{\text{\rm max}}(\bar{t})^{-1},\bar{t}]} \sum_{l=0}^2 \sum_{a=1}^3 R_{\text{\rm max}}(\bar{t})^{-l} \, |D^l(\mathscr{L}_{U^{(a)}}(g(t)))|^2 \leq \varepsilon^2$.
\item If $\Sigma \subset D$ is a leaf of the CMC foliation of $(M,g(\bar{t}))$ satisfying $\sup_{x \in \Sigma} d_{g(\bar{t})}(x,\partial D) \leq 10 \, \text{\rm area}_{g(\bar{t})}(\partial D)^{\frac{1}{2}}$, then $\sup_\Sigma \sum_{a=1}^3 R_{\text{\rm max}}(\bar{t}) \, |\langle U^{(a)},\nu \rangle|^2 \leq \varepsilon^2$, where $\nu$ denotes the unit normal vector to $\Sigma$ in $(M,g(\bar{t}))$.
\item If $\Sigma \subset D$ is a leaf of the CMC foliation of $(M,g(\bar{t}))$ satisfying $\sup_{x \in \Sigma} d_{g(\bar{t})}(x,\partial D) \leq 10 \, \text{\rm area}_{g(\bar{t})}(\partial D)^{\frac{1}{2}}$, then 
\[\sum_{a,b=1}^3 \bigg | \delta_{ab} - \text{\rm area}_{g(\bar{t})}(\Sigma)^{-2} \int_\Sigma \langle U^{(a)},U^{(b)} \rangle_{g(\bar{t})} \, d\mu_{g(\bar{t})} \bigg |^2 \leq \varepsilon^2.\]
\end{itemize}
\end{definition}

\begin{remark}
For each $x \in M \setminus D$, the point $(x,\bar{t})$ lies at the center of an evolving $\varepsilon_1$-neck by Proposition \ref{choice.of.theta}. 
\end{remark}

\begin{remark}
Since $D \subset \{x \in M: \lambda_1(x,\bar{t}) > \frac{1}{2} \theta R(x,\bar{t})\}$, Corollary \ref{bound.for.R_max} implies that $\text{\rm diam}_{g(\bar{t})}(D) \leq C \, R_{\text{\rm max}}(\bar{t})^{-\frac{1}{2}}$ and $\frac{1}{C} \, R_{\text{\rm max}}(\bar{t}) \leq R(x,\bar{t}) \leq R_{\text{\rm max}}(\bar{t})$ for all $x \in D$. Here, $C$ is a large constant that depends on $\theta$. By Lemma \ref{ode.for.U}, the vector fields $U^{(1)},U^{(2)},U^{(3)}$ satisfy $\sup_D \sum_{a=1}^3 |U^{(a)}|^2 \leq C \, R_{\text{\rm max}}(\bar{t})^{-1}$, where the norm is computed with respect to $g(\bar{t})$.
\end{remark}

\begin{lemma}
\label{openness.2} 
Suppose that the flow is $\varepsilon$-symmetric at time $\bar{t}$. If $\tilde{t}$ is sufficiently close to $\bar{t}$, then the flow is $2\varepsilon$-symmetric at time $\bar{t}$.
\end{lemma}

\textbf{Proof.} 
Since the flow is $\varepsilon$-symmetric at time $\bar{t}$, we can find a compact domain $D \subset M$ and time-independent vector fields $U^{(1)},U^{(2)},U^{(3)}$ which satisfy the conditions in Definition \ref{symmetry.of.cap}. In particular, every point in $(M \setminus D) \times \{\bar{t}\}$ is $\varepsilon$-symmetric. 

By continuity, we can find a slightly larger domain $D_1$ with the following properties: 
\begin{itemize} 
\item There exists a point $x \in \partial D_1$ such that $\lambda_1(x,\bar{t}) < \theta R(x,\bar{t})$.
\item For each $x \in D_1$, we have $\lambda_1(x,\bar{t}) > \frac{1}{2} \theta R(x,\bar{t})$.
\item The boundary $\partial D_1$ is a leaf of the CMC foliation of $(M,g(\bar{t}))$.
\item $\sup_{D_1 \times [\bar{t}-R_{\text{\rm max}}(\bar{t})^{-1},\bar{t}]} \sum_{l=0}^2 \sum_{a=1}^3 R_{\text{\rm max}}(\bar{t})^{-l} \, |D^l(\mathscr{L}_{U^{(a)}}(g(t)))|^2 \leq 2\varepsilon^2$.
\item If $\Sigma \subset D_1$ is a leaf of the CMC foliation of $(M,g(\bar{t}))$ satisfying $\sup_{x \in \Sigma} d_{g(\bar{t})}(x,\partial D_1) \leq 10 \, \text{\rm area}_{g(\bar{t})}(\partial D_1)^{\frac{1}{2}}$, then $\sup_\Sigma \sum_{a=1}^3 R_{\text{\rm max}}(\bar{t}) \, |\langle U^{(a)},\nu \rangle|^2 \leq 2\varepsilon^2$, where $\nu$ denotes the unit normal vector to $\Sigma$ in $(M,g(\bar{t}))$.
\item If $\Sigma \subset D_1$ is a leaf of the CMC foliation of $(M,g(\bar{t}))$ satisfying $\sup_{x \in \Sigma} d_{g(\bar{t})}(x,\partial D_1) \leq 10 \, \text{\rm area}_{g(\bar{t})}(\partial D_1)^{\frac{1}{2}}$, then 
\[\sum_{a,b=1}^3 \bigg | \delta_{ab} - \text{\rm area}_{g(\bar{t})}(\Sigma)^{-2} \int_\Sigma \langle U^{(a)},U^{(b)} \rangle_{g(\bar{t})} \, d\mu_{g(\bar{t})} \bigg |^2 \leq 2\varepsilon^2.\]
\end{itemize}
Let $D_0$ be a compact domain with the property that $\partial D_0$ is a leaf of the CMC foliation and $\partial D_0$ lies in between $\partial D$ and $\partial D_1$. Using Lemma \ref{openness.1}, we can find an open interval $I$ containing $\bar{t}$ such that every point on $(M \setminus D_0) \times I$ is $2\varepsilon$-symmetric. Moreover, if $\tilde{t}$ is sufficiently close to $\bar{t}$, then we can find a domain $\tilde{D}$ close to $D_1$ with the following properties: 
\begin{itemize} 
\item $D_0 \subset \tilde{D}$.
\item There exists a point $x \in \partial \tilde{D}$ such that $\lambda_1(x,\tilde{t}) < \theta R(x,\tilde{t})$.
\item For each $x \in \tilde{D}$, we have $\lambda_1(x,\tilde{t}) > \frac{1}{2} \theta R(x,\tilde{t})$.
\item The boundary $\partial \tilde{D}$ is a leaf of the CMC foliation of $(M,g(\tilde{t}))$.
\item $\sup_{\tilde{D} \times [\tilde{t}-R_{\text{\rm max}}(\tilde{t})^{-1},\tilde{t}]} \sum_{l=0}^2 \sum_{a=1}^3 R_{\text{\rm max}}(\tilde{t})^{-l} \, |D^l(\mathscr{L}_{U^{(a)}}(g(t)))|^2 \leq 4\varepsilon^2$.
\item If $\Sigma \subset \tilde{D}$ is a leaf of the CMC foliation of $(M,g(\tilde{t}))$ satisfying $\sup_{x \in \Sigma} d_{g(\tilde{t})}(x,\partial \tilde{D}) \leq 10 \, \text{\rm area}_{g(\tilde{t})}(\partial \tilde{D})^{\frac{1}{2}}$, then $\sup_\Sigma \sum_{a=1}^3 R_{\text{\rm max}}(\tilde{t}) \, |\langle U^{(a)},\nu \rangle|^2 \leq 4\varepsilon^2$, where $\nu$ denotes the unit normal vector to $\Sigma$ in $(M,g(\tilde{t}))$.
\item If $\Sigma \subset \tilde{D}$ is a leaf of the CMC foliation of $(M,g(\tilde{t}))$ satisfying $\sup_{x \in \Sigma} d_{g(\tilde{t})}(x,\partial \tilde{D}) \leq 10 \, \text{\rm area}_{g(\tilde{t})}(\partial \tilde{D})^{\frac{1}{2}}$, then 
\[\sum_{a,b=1}^3 \bigg | \delta_{ab} - \text{\rm area}_{g(\tilde{t})}(\Sigma)^{-2} \int_\Sigma \langle U^{(a)},U^{(b)} \rangle_{g(\tilde{t})} \, d\mu_{g(\tilde{t})} \bigg |^2 \leq 4\varepsilon^2.\]
\end{itemize}
Therefore, if $\tilde{t}$ is sufficiently close to $\bar{t}$, then the flow is $2\varepsilon$-symmetric at time $\tilde{t}$. \\

\begin{lemma}
\label{exact.symmetry}
Let us fix a time $\bar{t}$. Suppose that, for each $\varepsilon>0$, the flow is $\varepsilon$-symmetric at time $\bar{t}$. Then the manifold $(M,g(\bar{t}))$ is rotationally symmetric.
\end{lemma}

\textbf{Proof.} 
We consider the vector fields in Definition \ref{symmetry.of.cap}, and pass to the limit as $\varepsilon \to 0$. Hence, we can find a compact domain $D \subset M$ and vector fields $U^{(1)},U^{(2)},U^{(3)}$ on $D$ with the following properties: 
\begin{itemize} 
\item The boundary $\partial D$ is a leaf of the CMC foliation of $(M,g(\bar{t}))$.
\item The metric $g(\bar{t})$ is rotationally symmetric on $M \setminus D$. 
\item The vector fields $U^{(a)}$ satisfy $\mathscr{L}_{U^{(a)}}(g(\bar{t})) = 0$ in $D$.
\item The vector fields $U^{(a)}$ are tangential along $\partial D$.
\item $\text{\rm area}_{g(\bar{t})}(\partial D)^{-2} \int_{\partial D} \langle U^{(a)},U^{(b)} \rangle_{g(\bar{t})} \, d\mu_{g(\bar{t})} = \delta_{ab}$.
\end{itemize}
In particular, $(\partial D,g(\bar{t}))$ is a round sphere, and every Killing vector field on $(\partial D,g(\bar{t}))$ can be extended to a Killing vector field on $(D,g(\bar{t}))$. This implies that the metric $g(\bar{t})$ is rotationally symmetric in $D$. This completes the proof of Lemma \ref{exact.symmetry}. \\


We now proceed with the proof of Theorem \ref{main.thm.b}. We first show that we can find a sequence of times where the solution is arbitrarily close to the Bryant soliton. This argument relies on the Harnack inequality together with the classification of steady gradient Ricci solitons in \cite{Brendle}.

\begin{proposition} 
\label{solution.is.close.to.Bryant.soliton.at.time.hat.t_k}
We can find a sequence of times $\hat{t}_k \to -\infty$ and a sequence of points $\hat{p}_k \in M$ with the following property. If we perform a parabolic rescaling around the point $(\hat{p}_k,\hat{t}_k)$ by the factor $R_{\text{\rm max}}(\hat{t}_k)^{\frac{1}{2}}$, then the rescaled flows converge to the Bryant soliton in the Cheeger-Gromov sense. Moreover, the points $\hat{p}_k$ converge to the tip of the Bryant soliton, and we have $\frac{R(\hat{p}_k,\hat{t}_k)}{R_{\text{\rm max}}(\hat{t}_k)} \to 1$ as $k \to \infty$.
\end{proposition}

\textbf{Proof.} 
By \cite{Zhang}, $(M,g(t))$ is a Type II ancient solution, i.e. $\sup_{(x,t) \in M \times (-\infty,0]} (-t) \, R(x,t) = \infty$. We now argue as in Section 16 of \cite{Hamilton4} to extract a Type II blow-up limit. For $k$ large, we choose a point $(\hat{p}_k,\hat{t}_k) \in M \times (-k,0)$ with the property that 
\[\sup_{(x,t) \in M \times (-k,0)} \Big ( 1+\frac{t}{k} \Big ) \, (-t) \, R(x,t) \leq \Big ( 1+\frac{1}{k} \Big ) \, \Big ( 1+\frac{\hat{t}_k}{k} \Big ) \, (-\hat{t}_k) \, R(\hat{p}_k,\hat{t}_k).\] 
In particular, $R_{\text{\rm max}}(\hat{t}_k) \leq (1+\frac{1}{k}) \, R(\hat{p}_k,\hat{t}_k)$. Since $(M,g(t))$ is a Type II ancient solution, we know that 
\[\sup_{(x,t) \in M \times (-k,0)} \Big ( 1+\frac{t}{k} \Big ) \, (-t) \, R(x,t) \to \infty,\] 
hence 
\[\Big ( 1+\frac{\hat{t}_k}{k} \Big ) \, (-\hat{t}_k) \, R(\hat{p}_k,\hat{t}_k) \to \infty.\] 
This implies $(-\hat{t}_k) \, R(\hat{p}_k,\hat{t}_k) \to \infty$, $(k+\hat{t}_k) \, R(\hat{p}_k,\hat{t}_k) \to \infty$, and 
\[\limsup_{k \to \infty} \sup_{(x,t) \in M \times [\hat{t}_k - A R(\hat{p}_k,\hat{t}_k)^{-1},\hat{t}_k + A R(\hat{p}_k,\hat{t}_k)^{-1}]} \frac{R(x,t)}{R(\hat{p}_k,\hat{t}_k)} \leq 1\] 
for every fixed $A$. 

We now rescale around the point $(\hat{p}_k,\hat{t}_k)$ by the factor $R(\hat{p}_k,\hat{t}_k)^{\frac{1}{2}}$. Passing to the limit as $k \to \infty$, we obtain an eternal solution to the Ricci flow which is complete; $\kappa$-noncollapsed; has nonnegative sectional curvature; and has scalar curvature at most $1$ at each point in space-time. Moreover, there exists a point on the limiting solution where the scalar curvature is equal to $1$. Therefore, the limiting solution attains equality in Hamilton's Harnack inequality \cite{Hamilton2}, and consequently must be a steady gradient Ricci soliton \cite{Hamilton3}. By \cite{Brendle}, the limit flow must be the Bryant soliton. This completes the proof of Proposition \ref{solution.is.close.to.Bryant.soliton.at.time.hat.t_k}. \\

\begin{corollary}
\label{solution.is.almost.symmetric.at.time.hat.t_k}
There exists a sequence $\hat{\varepsilon}_k \to 0$ with the following properties. For each $t \in [\hat{t}_k-\hat{\varepsilon}_k^{-2} R_{\text{\rm max}}(\hat{t}_k)^{-1},\hat{t}_k]$, we have $(1-\hat{\varepsilon}_k) \, R_{\text{\rm max}}(\hat{t}_k) \leq R(\hat{p}_k,t) \leq R_{\text{\rm max}}(t) \leq R_{\text{\rm max}}(\hat{t}_k)$. Moreover, for each $t \in [\hat{t}_k-\hat{\varepsilon}_k^{-2} R_{\text{\rm max}}(\hat{t}_k)^{-1},\hat{t}_k]$, the flow is $\hat{\varepsilon}_k$-symmetric at time $t$. 
\end{corollary} 

\textbf{Proof.} 
The Harnack inequality (cf. \cite{Hamilton2}) implies that $R_{\text{\rm max}}(t) \leq R_{\text{\rm max}}(\hat{t}_k)$ for each $t \leq \hat{t}_k$. The remaining statements follow by combining Proposition \ref{solution.is.close.to.Bryant.soliton.at.time.hat.t_k} and Theorem \ref{canonical.neighborhood.theorem}. \\

\textbf{From now on, we assume that the ancient solution $(M,g(t))$ is not rotationally symmetric.} In view of Corollary \ref{bound.for.R_max}, we can find a sequence of positive real numbers $\varepsilon_k$ with the following properties: 
\begin{itemize}
\item $\varepsilon_k \to 0$. 
\item $\varepsilon_k \geq 2\hat{\varepsilon}_k$. 
\item If a point $(x,t)$ in spacetime satisfies $R(x,t) \leq \hat{\varepsilon}_k R_{\text{\rm max}}(t)$, then $(x,t)$ lies at the center of an evolving $\varepsilon_k^2$-neck. 
\end{itemize}
For each $k$, we define 
\[t_k = \inf \{t \in [\hat{t}_k,0]: \text{\rm The flow is not $\varepsilon_k$-symmetric at time $t$}\}.\] 
For abbreviation, let $R_{\text{\rm max}}(t_k) = r_k^{-2}$. The Harnack inequality \cite{Hamilton2} implies $R_{\text{\rm max}}(t) \leq r_k^{-2}$ for all $t \leq t_k$. \\

\begin{lemma}
\label{varepsilon_k.symmetry.at.earlier.times}
If $t \in [\hat{t}_k-\hat{\varepsilon}_k^{-2} R_{\text{\rm max}}(\hat{t}_k)^{-1},t_k)$, then the flow is $\varepsilon_k$-symmetric at time $t$. In particular, if $(x,t) \in M \times [\hat{t}_k-\hat{\varepsilon}_k^{-2} R_{\text{\rm max}}(\hat{t}_k)^{-1},t_k)$ is a point in spacetime satisfying $\lambda_1(x,t) < \frac{1}{2} \theta R(x,t)$, then the point $(x,t)$ is $\varepsilon_k$-symmetric.
\end{lemma}

\textbf{Proof.} 
The first statement follows immediately from the definition of $t_k$. The second statement follows from the first statement, keeping in mind Definition \ref{symmetry.of.cap}. \\

\begin{lemma}
\label{limit.of.t_k}
The sequence $t_k$ satisfies $\lim_{k \to \infty} t_k = -\infty$. 
\end{lemma}

\textbf{Proof.} 
Suppose that $\limsup_{k \to \infty} t_k > -\infty$. Let us consider an arbitrary time $\bar{t} < \limsup_{k \to \infty} t_k$. Then there exist arbitrarily large integers $k$ with the property that $\bar{t} \in [\hat{t}_k,t_k)$. By Lemma \ref{varepsilon_k.symmetry.at.earlier.times}, there exist arbitrarily large integers $k$ with the property that the flow is $\varepsilon_k$-symmetric at time $\bar{t}$. Since $\varepsilon_k \to 0$, Lemma \ref{exact.symmetry} implies that $(M,g(\bar{t}))$ is rotationally symmetric. 

To summarize, we have shown that the solution $(M,g(t))$ is rotationally symmetric for all $t < \limsup_{k \to \infty} t_k$. By the uniqueness result in \cite{Chen-Zhu}, the solution is rotationally symmetric for all $t$, contrary to our assumption. This completes the proof of Lemma \ref{limit.of.t_k}. \\

In the next step, we show that, at time $t_k$, the solution is close to the Bryant soliton. This argument relies in a crucial way on Theorem \ref{main.thm.a}.

\begin{proposition}
\label{solution.is.close.to.Bryant.soliton.at.time.t_k} 
There exists a sequence of points $p_k \in M$ with the following properties. If we perform a parabolic rescaling around the point $(p_k,t_k)$ by the factor $R_{\text{\rm max}}(t_k)^{\frac{1}{2}} = r_k^{-1}$, then the rescaled flows converge to the Bryant soliton in the Cheeger-Gromov sense. Moreover, the points $p_k$ converge to the tip of the Bryant soliton, and we have $r_k^2 \, R(p_k,t_k) \to 1$ as $k \to \infty$.
\end{proposition}

\textbf{Proof.} 
For each $k$, the manifold $(M,g(t_k))$ contains a point which does not lie on a neck. Hence, we can find a sequence of points $q_k \in M$ such that $\liminf_{k \to \infty} \frac{\lambda_1(q_k,t_k)}{R(q_k,t_k)} > 0$. By Corollary \ref{bound.for.R_max}, $\liminf_{k \to \infty} r_k^2 \, R(q_k,t_k) > 0$. This implies $\liminf_{k \to \infty} r_k^2 \, \lambda_1(q_k,t_k) > 0$. We now rescale the flow $(M,g(t))$ around the point $(q_k,t_k)$ by the factor $r_k^{-1}$. Passing to the limit as $k \to \infty$, we obtain a noncompact ancient $\kappa$-solution $(M^\infty,g^\infty(s))$. Since $\liminf_{k \to \infty} r_k^2 \, \lambda_1(q_k,t_k) > 0$, the limit manifold $(M^\infty,g^\infty(0))$ does not split off a line. By the uniqueness result in \cite{Chen-Zhu}, the manifold $(M^\infty,g^\infty(s))$ does not split off a line for any $s \leq 0$. By the strict maximum principle, the limit flow $(M^\infty,g^\infty(s))$ has positive sectional curvature for each $s \leq 0$. 

We claim that the limiting flow $(M^\infty,g^\infty(s))$ is rotationally symmetric. To prove this, we fix an arbitrary time $\bar{s} < 0$. Since $R_{\text{\rm max}}(\hat{t}_k) \leq r_k^{-2}$, it follows that $t_k + r_k^2 \, \bar{s} \in [\hat{t}_k-\hat{\varepsilon}_k^{-2} R_{\text{\rm max}}(\hat{t}_k)^{-1},t_k)$ if $k$ is sufficiently large. By Lemma \ref{varepsilon_k.symmetry.at.earlier.times}, the original flow is $\varepsilon_k$-symmetric at time $t_k + r_k^2 \, \bar{s}$, provided that $k$ is sufficiently large. By the Harnack inequality, $R_{\text{\rm max}}(t_k+r_k^2 \, \bar{s}) \leq r_k^{-2}$. On the other hand, since $(M^\infty,g^\infty(\bar{s}))$ has positive sectional curvature, we obtain $\liminf_{k \to \infty} r_k^2 \, \lambda_1(q_k,t_k+r_k^2 \, \bar{s}) > 0$. 
Therefore, the cap in $(M,g(t_k + r_k^2 \, \bar{s}))$ has diameter $\lesssim r_k$, the scalar curvature on the cap is $\sim r_k^{-2}$, and the cap has distance $\lesssim r_k$ from the point $q_k$. We now pass to the limit as $k \to \infty$. In the limit, we obtain a domain $D^\infty \subset M^\infty$ and vector fields $U^{(\infty,1)},U^{(\infty,2)},U^{(\infty,3)}$ on $D^\infty$ with the following properties: 
\begin{itemize} 
\item The boundary $\partial D^\infty$ is a leaf of the CMC foliation of $(M^\infty,g^\infty(\bar{s}))$.
\item The metric $g(\bar{s})$ is rotationally symmetric on $M^\infty \setminus D^\infty$. 
\item The vector fields $U^{(\infty,a)}$ satisfy $\mathscr{L}_{U^{(\infty,a)}}(g^\infty(\bar{s})) = 0$ in $D^\infty$.
\item The vector fields $U^{(\infty,a)}$ are tangential along $\partial D^\infty$.
\item $\text{\rm area}_{g^\infty(\bar{s})}(\partial D^\infty)^{-2} \int_{\partial D^\infty} \langle U^{(\infty,a)},U^{(\infty,b)} \rangle_{g^\infty(\bar{s})} \, d\mu_{g^\infty(\bar{s})} = \delta_{ab}$.
\end{itemize}
Thus, we conclude that the limiting manifold $(M^\infty,g^\infty(\bar{s}))$ is rotationally symmetric. 

To summarize, we have shown that $(M^\infty,g^\infty(s))$ is a noncompact ancient $\kappa$-solution which is rotationally symmetric and has positive sectional curvature. By Theorem \ref{main.thm.a}, the limiting flow $(M^\infty,g^\infty(s))$ must be isometric to the Bryant soliton, up to scaling. 

Finally, we claim that $R_{g^\infty(0)}(p_\infty) = 1$, where $p_\infty \in M^\infty$ denotes the tip of the limiting soliton $(M^\infty,g^\infty(s))$. To see this, consider a sequence of points $p_k \in M$ converging to $p_\infty$. Clearly, $R_{g^\infty(0)}(p_\infty) = \lim_{k \to \infty} r_k^2 \, R(p_k,t_k) \in (0,1]$. Using Proposition \ref{no.small.necks}, we can find a large constant $A$ such that 
\[\sup_{x \in M \setminus B_{g(t_k)}(p_k,A R(p_k,t_k)^{-\frac{1}{2}})} \frac{R(x,t_k)}{R(p_k,t_k)} \leq \frac{1}{2}\] 
if $k$ is sufficiently large. Moreover, since the scalar curvature of $(M^\infty,g^\infty(0))$ attains its maximum at the point $p_\infty$, we obtain 
\[\limsup_{k \to \infty} \sup_{x \in B_{g(t_k)}(p_k,A R(p_k,t_k)^{-\frac{1}{2}})} \frac{R(x,t_k)}{R(p_k,t_k)} \leq 1\] 
for every fixed $A$. Putting these facts together, we conclude that 
\[\limsup_{k \to \infty} \sup_{x \in M} \frac{R(x,t_k)}{R(p_k,t_k)} \leq 1.\] 
Thus, $R_{g^\infty(0)}(p_\infty) = \lim_{k \to \infty} r_k^2 \, R(p_k,t_k) \geq 1$. This completes the proof of Proposition \ref{solution.is.close.to.Bryant.soliton.at.time.t_k}. \\

\begin{corollary} 
\label{perturbation.of.Bryant.soliton}
There exists a sequence of positive real numbers $\delta_k \to 0$ such that $\delta_k \geq 2\varepsilon_k$ for each $k$ and the following statements hold when $k$ is sufficiently large: 
\begin{itemize}
\item For each $t \in [t_k-\delta_k^{-1} r_k^2,t_k]$, we have $\frac{1-\delta_k}{3} \, g \leq r_k^2 \, \text{\rm Ric} \leq \frac{1+\delta_k}{3} \, g$ at the point $(p_k,t)$.
\item The scalar curvature satisfies $\frac{1}{2K} \, (r_k^{-1} \, d_{g(t)}(p_k,x)+1)^{-1} \leq r_k^2 \, R(x,t) \leq 2K \, (r_k^{-1} \, d_{g(t)}(p_k,x)+1)^{-1}$ for all points $(x,t) \in B_{g(t_k)}(p_k,\delta_k^{-1} r_k) \times [t_k-\delta_k^{-1} r_k^2,t_k]$.
\item There exists a nonnegative function $f: B_{g(t_k)}(p_k,\delta_k^{-1} r_k) \times [t_k-\delta_k^{-1} r_k^2,t_k] \to \mathbb{R}$ such that $|\text{\rm Ric}-D^2 f| \leq \delta_k r_k^{-2}$, $|\Delta f + |\nabla f|^2 - r_k^{-2}| \leq \delta_k r_k^{-2}$, and $|\frac{\partial}{\partial t} f + |\nabla f|^2| \leq \delta_k r_k^{-2}$. 
\item The function $f$ satisfies $\frac{1}{2K} \, (r_k^{-1} \, d_{g(t)}(p_k,x)+1) \leq f(x,t)+1 \leq 2K \, (r_k^{-1} \, d_{g(t)}(p_k,x)+1)$ for all points $(x,t) \in B_{g(t_k)}(p_k,\delta_k^{-1} r_k) \times [t_k-\delta_k^{-1} r_k^2,t_k]$.
\end{itemize}
Here, $K \geq 10$ is a universal constant.
\end{corollary}

\textbf{Proof.} 
On the Bryant soliton, the eigenvalues of the Ricci tensor at the tip are equal to $\frac{1}{3}$. Moreover, on the Bryant soliton, the scalar curvature satisfies $\frac{1}{K} \, (d(p,x)+1)^{-1} \leq R \leq K \, (d(p,x)+1)^{-1}$, where $p$ denotes the tip of the Bryant soliton and $K$ is a universal constant. Furthermore, on the Bryant soliton, the potential function $f$ satisfies $\frac{1}{K} \, (d(p,x)+1) \leq f+1 \leq K \, (d(p,x)+1)$, where again $p$ denotes the tip of the Bryant soliton and $K$ is a universal constant. Finally, the potential function $f$ satisfies $\text{\rm Ric} = D^2 f$, $\Delta f + |\nabla f|^2 = 1$, and $\frac{\partial}{\partial t} f + |\nabla f|^2 = 0$. The assertion now follows from Proposition \ref{solution.is.close.to.Bryant.soliton.at.time.t_k}. \\

\begin{corollary}
\label{control.of.R_max}
For each $t \in [t_k-\delta_k^{-1} r_k^2,t_k]$, we have $(1-\delta_k) r_k^{-2} \leq R(p_k,t) \leq R_{\text{\rm max}}(t) \leq r_k^{-2}$.
\end{corollary}

\textbf{Proof.} 
The Harnack inequality (cf. \cite{Hamilton2}) implies that $R_{\text{\rm max}}(t) \leq r_k^{-2}$ for each $t \leq t_k$. Moreover, Corollary \ref{perturbation.of.Bryant.soliton} implies $R(p_k,t) \geq (1-\delta_k) r_k^{-2}$ for each $t \in [t_k-\delta_k^{-1} r_k^2,t_k]$. This completes the proof of Corollary \ref{control.of.R_max}. \\

\begin{lemma}
\label{derivative.of.distance.function}
The time derivative of the distance function satisfies $0 \leq -\frac{d}{dt} d_{g(t)}(p_k,x) \leq 80 r_k^{-1}$ for all $(x,t) \in M \times [t_k-\delta_k^{-1} r_k^2,t_k]$.
\end{lemma}

\textbf{Proof.} 
Using Lemma 8.3(b) in \cite{Perelman1}, we obtain $0 \leq -\frac{d}{dt} d_{g(t)}(p_k,x) \leq 80 \, R_{\text{\rm max}}(t)^{\frac{1}{2}} \leq 80r_k^{-1}$ for all $(x,t) \in M \times [t_k-\delta_k^{-1} r_k^2,t_k]$. This completes the proof of Lemma \ref{derivative.of.distance.function}. \\

In view of Theorem \ref{canonical.neighborhood.theorem} and Corollary \ref{control.of.R_max}, we can find a large constant $\Lambda$ with the following properties:
\begin{itemize} 
\item $L \sqrt{\frac{4K}{\Lambda}} \leq 10^{-6}$.
\item If $(\bar{x},\bar{t}) \in M \times [t_k-\delta_k^{-1} r_k^2,t_k]$ is a point in spacetime satisfying $d_{g(\bar{t})}(p_k,\bar{x}) \geq \Lambda r_k$, then $\lambda_1(x,t) < \frac{1}{2} \theta R(x,t)$ for all points $(x,t) \in B_{g(\bar{t})}(\bar{x},L \, R(\bar{x},\bar{t})^{-\frac{1}{2}}) \times [\bar{t} - L \, R(\bar{x},\bar{t})^{-1},\bar{t}]$.
\end{itemize}

\begin{lemma}
\label{improved.symmetry.for.faraway.points}
If $k$ is sufficiently large, then the following holds. If $(\bar{x},\bar{t}) \in M \times [t_k-\delta_k^{-1} r_k^2,t_k]$ satisfies $d_{g(\bar{t})}(p_k,\bar{x}) \geq \Lambda r_k$, then $(\bar{x},\bar{t})$ is $\frac{\varepsilon_k}{2}$-symmetric.
\end{lemma}

\textbf{Proof.} 
We distinguish two cases: 

\textit{Case 1:} Suppose first that $R(\bar{x},\bar{t}) \leq \hat{\varepsilon}_k R_{\text{\rm max}}(\bar{t})$. By our choice of $\varepsilon_k$, the point $(\bar{x},\bar{t})$ lies at the center of an evolving $\varepsilon_k^2$-neck, and this directly implies that $(\bar{x},\bar{t})$ is $\frac{\varepsilon_k}{2}$-symmetric. 

\textit{Case 2:} Suppose next that $R(\bar{x},\bar{t}) \geq \hat{\varepsilon}_k R_{\text{\rm max}}(\bar{t})$. Corollary \ref{control.of.R_max} implies $R_{\text{\rm max}}(\bar{t}) \geq \frac{1}{2} \, r_k^{-2}$, hence $R(\bar{x},\bar{t}) \geq \frac{1}{2} \, \hat{\varepsilon}_k r_k^{-2}$. On the other hand, $R_{\text{\rm max}}(\hat{t}_k) \leq r_k^{-2}$. Hence, if $k$ is sufficiently large, then we obtain 
\begin{align*} 
\bar{t} - L \, R(\bar{x},\bar{t})^{-1} 
&\geq t_k - \delta_k^{-1} r_k^2 - 2L \hat{\varepsilon}_k^{-1} r_k^2 \\ 
&\geq \hat{t}_k - \delta_k^{-1} R_{\text{\rm max}}(\hat{t}_k)^{-1} - 2L \hat{\varepsilon}_k^{-1} R_{\text{\rm max}}(\hat{t}_k)^{-1} \\ 
&\geq \hat{t}_k - \hat{\varepsilon}_k^{-2} R_{\text{\rm max}}(\hat{t}_k)^{-1}. 
\end{align*}
By definition of $\Lambda$, we have $\lambda_1(x,t) < \frac{1}{2} \theta R(x,t)$ for all points $(x,t) \in B_{g(\bar{t})}(\bar{x},L \, R(\bar{x},\bar{t})^{-\frac{1}{2}}) \times [\bar{t} - L \, R(\bar{x},\bar{t})^{-1},\bar{t}]$. Hence, by Proposition \ref{choice.of.theta}, every point in $B_{g(\bar{t})}(\bar{x},L \, R(\bar{x},\bar{t})^{-\frac{1}{2}}) \times [\bar{t} - L \, R(\bar{x},\bar{t})^{-1},\bar{t}]$ lies at the center of an evolving $\varepsilon_1$-neck. Moreover, by Lemma \ref{varepsilon_k.symmetry.at.earlier.times}, every point in $B_{g(\bar{t})}(\bar{x},L \, R(\bar{x},\bar{t})^{-\frac{1}{2}}) \times [\bar{t} - L \, R(\bar{x},\bar{t})^{-1},\bar{t})$ is $\varepsilon_k$-symmetric. Using the Neck Improvement Theorem, we conclude that the point $(\bar{x},\bar{t})$ is $\frac{\varepsilon_k}{2}$-symmetric. This completes the proof of Lemma \ref{improved.symmetry.for.faraway.points}. \\

\begin{proposition}
\label{iteration}
If $k$ is sufficiently large, then the following holds. If $(\bar{x},\bar{t}) \in M \times [t_k-2^{-j} \delta_k^{-1} r_k^2,t_k]$ satisfies $2^{\frac{j}{400}} \, \Lambda r_k \leq d_{g(\bar{t})}(p_k,\bar{x}) \leq (400KL)^{-j} \, \delta_k^{-1} \, r_k$, then $(\bar{x},\bar{t})$ is $2^{-j-1} \varepsilon_k$-symmetric. 
\end{proposition}

\textbf{Proof.} 
We argue by induction on $j$. For $j=0$, the assertion follows from Lemma \ref{improved.symmetry.for.faraway.points}. 

We now assume that $j \geq 1$ and the assertion holds for $j-1$. We will show that the assertion holds for $j$. To that end, we consider a point $(\bar{x},\bar{t}) \in M \times [t_k-2^{-j} \delta_k^{-1} r_k^2,t_k]$ such that $2^{\frac{j}{400}} \, \Lambda r_k \leq d_{g(\bar{t})}(p_k,\bar{x}) \leq (400KL)^{-j} \, \delta_k^{-1} \, r_k$. Clearly, $\lambda_1(\bar{x},\bar{t}) < \frac{1}{2} \theta R(\bar{x},\bar{t})$ by definition of $\Lambda$. By Proposition \ref{choice.of.theta}, $(\bar{x},\bar{t})$ lies at the center of an evolving $\varepsilon_1$-neck. Let $R(\bar{x},\bar{t}) = r^{-2}$. We will show that every point in $B_{g(\bar{t})}(\bar{x},Lr) \times [\bar{t}-Lr^2,\bar{t}]$ is $2^{-j} \varepsilon_k$-symmetric. By Corollary \ref{perturbation.of.Bryant.soliton}, $r^2 \leq 4K \, r_k \, d_{g(\bar{t})}(p_k,\bar{x})$. This implies 
\begin{align*} 
\bar{t}-Lr^2 
&\geq \bar{t} - 4KL \, r_k \, d_{g(\bar{t})}(p_k,\bar{x}) \\ 
&\geq \bar{t} - 4KL \, (400KL)^{-j} \, \delta_k^{-1} \, r_k^2 \\ 
&\geq \bar{t} - 2^{-j} \, \delta_k^{-1} \, r_k^2 \\ 
&\geq t_k - 2^{-j+1} \, \delta_k^{-1} \, r_k^2. 
\end{align*} 
In the next step, we observe that $r^2 \leq 4K \, r_k \, d_{g(\bar{t})}(p_k,\bar{x}) \leq \frac{4K}{\Lambda} \, d_{g(\bar{t})}(p_k,\bar{x})^2$. Since $L \sqrt{\frac{4K}{\Lambda}} \leq 10^{-6}$, we obtain $r \leq \sqrt{\frac{4K}{\Lambda}} \, d_{g(\bar{t})}(p_k,\bar{x}) \leq 10^{-6} \, L^{-1} \, d_{g(\bar{t})}(p_k,\bar{x})$. Consequently, 
\begin{align*} 
d_{g(\bar{t})}(p_k,x) 
&\geq d_{g(\bar{t})}(p_k,\bar{x}) - Lr \\ 
&\geq (1-10^{-6}) \, d_{g(\bar{t})}(p_k,\bar{x}) \\ 
&\geq (1-10^{-6}) \, 2^{\frac{j}{400}} \, \Lambda r_k \\ 
&\geq 2^{\frac{j-1}{400}} \, \Lambda r_k 
\end{align*}
for all $x \in B_{g(\bar{t})}(\bar{x},Lr)$. On the other hand, $r_k = R_{\text{\rm max}}(t_k)^{-\frac{1}{2}} \leq R_{\text{\rm max}}(\bar{t})^{-\frac{1}{2}} \leq r$ by the Harnack inequality. Using this inequality together with  the inequality $r^2 \leq 4K \, r_k \, d_{g(\bar{t})}(p_k,\bar{x})$, we obtain 
\begin{align*} 
d_{g(\bar{t})}(p_k,x) + 80Lr^2 r_k^{-1} 
&\leq d_{g(\bar{t})}(p_k,\bar{x}) + Lr + 80Lr^2 r_k^{-1} \\ 
&\leq d_{g(\bar{t})}(p_k,\bar{x}) + 81Lr^2 r_k^{-1} \\ 
&\leq 400KL \, d_{g(\bar{t})}(p_k,\bar{x}) \\ 
&\leq (400KL)^{-j+1} \, \delta_k^{-1} \, r_k 
\end{align*} 
for all $x \in B_{g(\bar{t})}(\bar{x},Lr)$. Lemma \ref{derivative.of.distance.function} gives 
\[d_{g(\bar{t})}(p_k,x) \leq d_{g(t)}(p_k,x) \leq d_{g(\bar{t})}(p_k,x) + 80Lr^2 r_k^{-1},\]
hence 
\[2^{\frac{j-1}{400}} \, \Lambda r_k \leq d_{g(t)}(p_k,x) \leq (400KL)^{-j+1} \, \delta_k^{-1} \, r_k\] 
for all $(x,t) \in B_{g(\bar{t})}(\bar{x},Lr) \times [\bar{t}-Lr^2,\bar{t}]$. Therefore, the induction hypothesis implies that every point in $B_{g(\bar{t})}(\bar{x},Lr) \times [\bar{t}-Lr^2,\bar{t}]$ is $2^{-j} \varepsilon_k$-symmetric. Consequently, the point $(\bar{x},\bar{t})$ is $2^{-j-1} \varepsilon_k$-symmetric by the Neck Improvement Theorem. This completes the proof of Proposition \ref{iteration}. \\

\begin{lemma}
\label{existence.of.almost.Killing.vector.fields}
If $j$ is sufficiently large and $k$ is sufficiently large depending on $j$, then the following holds. Given any $\bar{t} \in [t_k-2^{\frac{j}{100}} \, r_k^2,t_k]$, there exist time-independent vector fields $U^{(1)},U^{(2)},U^{(3)}$ on $B_{g(\bar{t})}(p_k,2^{\frac{j}{400}} \, \Lambda r_k)$ with the following properties: 
\begin{itemize} 
\item $|\mathscr{L}_{U^{(a)}}(g(t))|+r_k \, |D(\mathscr{L}_{U^{(a)}}(g(t)))| \leq C \, (r_k^{-1} \, d_{g(t)}(p_k,x)+1)^{-100} \, \varepsilon_k$ for all $(x,t) \in B_{g(\bar{t})}(p_k,2^{\frac{j}{400}} \, \Lambda r_k) \times [\bar{t}-r_k^2,\bar{t}]$. 
\item If $t \in [\bar{t}-r_k^2,\bar{t}]$ and $x \in B_{g(\bar{t})}(p_k,2^{\frac{j}{400}} \, \Lambda r_k) \setminus B_{g(\bar{t})}(p_k,2\Lambda r_k)$, then $r_k^{-1} \, |\langle U^{(a)},\nu \rangle| \leq C \, (r_k^{-1} \, d_{g(t)}(p_k,x)+1)^{-100} \, \varepsilon_k$, where $\nu$ denotes the unit normal to the CMC foliation of $(M,g(t))$.
\item If $t \in [\bar{t}-r_k^2,\bar{t}]$ and $x \in B_{g(\bar{t})}(p_k,2^{\frac{j}{400}} \, \Lambda r_k) \setminus B_{g(\bar{t})}(p_k,2\Lambda r_k)$, then 
\begin{align*} 
&\sum_{a,b=1}^3 \bigg | \delta_{ab} - \text{\rm area}_{g(t)}(\Sigma)^{-2} \int_\Sigma \langle U^{(a)},U^{(b)} \rangle_{g(t)} \, d\mu_{g(t)} \bigg | \\ 
&\leq C \, (r_k^{-1} \, d_{g(t)}(p_k,x)+1)^{-100} \, \varepsilon_k, 
\end{align*}
where $\Sigma$ denotes the leaf of the CMC foliation passing through $(x,t)$.
\end{itemize}
Moreover, on the ball $B_{g(\bar{t})}(p_k,2^{\frac{j}{400}} \, \Lambda r_k)$, the vector fields $U^{(1)},U^{(2)},U^{(3)}$ are close to the standard rotation vector fields on the Bryant soliton in the $C^2$-norm.
\end{lemma}

\textbf{Proof.} 
We proceed in two steps:

\textit{Step 1:} Suppose first that $\bar{t} \in [t_k-2^{\frac{j}{100}} \, r_k^2,t_k)$. By Lemma \ref{varepsilon_k.symmetry.at.earlier.times}, the flow is $\varepsilon_k$-symmetric at time $\bar{t}$. Moreover, if $\bar{x} \in B_{g(\bar{t})}(p_k,2^{\frac{j}{400}} \, \Lambda r_k) \setminus B_{g(\bar{t})}(p_k,\Lambda r_k)$, then the point $(\bar{x},\bar{t})$ is $C \, (r_k^{-1} \, d_{g(\bar{t})}(p_k,\bar{x}))^{-400} \, \varepsilon_k$-symmetric by Proposition \ref{iteration}. By a repeated application of Corollary \ref{gluing}, we can construct vector fields $U^{(1)},U^{(2)},U^{(3)}$ satisfying the conditions above. Moreover, in view of Definition \ref{symmetry.of.necks} and Definition \ref{symmetry.of.cap}, the Lie derivatives $\mathscr{L}_{U^{(1)}}(g),\mathscr{L}_{U^{(2)}}(g),\mathscr{L}_{U^{(3)}}(g)$ are small in the $C^2$-norm. Consequently, the vector fields $U^{(1)},U^{(2)},U^{(3)}$ are close to the standard rotation vector fields on the Bryant soliton in the $C^{2,\frac{1}{2}}$-norm. 

\textit{Step 2:} Suppose next that $\bar{t} = t_k$. In this case, the assertion follows from the result in Step 1 by passing to the limit. Since the vector fields constructed in Step 1 are bounded in $C^{2,\frac{1}{2}}$, we may take the limit in $C^2$. \\

\begin{lemma}
\label{vector.field.comparison.on.cap}
If $j$ is sufficiently large and $k$ is sufficiently large depending on $j$, then the following statement holds. Consider a time $\bar{t} \in [t_k-2^{\frac{j}{100}} \, r_k^2,t_k]$. Suppose that $U^{(1)},U^{(2)},U^{(3)}$ are vector fields on $B_{g(\bar{t})}(p_k,2^{\frac{j}{400}} \, \Lambda r_k)$ with the following properties: 
\begin{itemize} 
\item $|\mathscr{L}_{U^{(a)}}(g(\bar{t}))|+r_k \, |D(\mathscr{L}_{U^{(a)}}(g(\bar{t})))| \leq C \, (r_k^{-1} \, d_{g(\bar{t})}(p_k,x)+1)^{-100} \, \varepsilon_k$ for all $x \in B_{g(\bar{t})}(p_k,2^{\frac{j}{400}} \, \Lambda r_k)$. 
\item If $x \in B_{g(\bar{t})}(p_k,2^{\frac{j}{400}} \, \Lambda r_k) \setminus B_{g(\bar{t})}(p_k,4\Lambda r_k)$, then $r_k^{-1} \, |\langle U^{(a)},\nu \rangle| \leq C \, (r_k^{-1} \, d_{g(\bar{t})}(p_k,x)+1)^{-100} \, \varepsilon_k$, where $\nu$ denotes the unit normal to the CMC foliation of $(M,g(\bar{t}))$.
\item If $x \in B_{g(\bar{t})}(p_k,2^{\frac{j}{400}} \, \Lambda r_k) \setminus B_{g(\bar{t})}(p_k,4\Lambda r_k)$, then 
\begin{align*} 
&\sum_{a,b=1}^3 \bigg | \delta_{ab} - \text{\rm area}_{g(\bar{t})}(\Sigma)^{-2} \int_\Sigma \langle U^{(a)},U^{(b)} \rangle_{g(\bar{t})} \, d\mu_{g(\bar{t})} \bigg | \\ 
&\leq C \, (r_k^{-1} \, d_{g(\bar{t})}(p_k,x)+1)^{-100} \, \varepsilon_k, 
\end{align*} 
where $\Sigma$ denotes the leaf of the CMC foliation passing through $(x,\bar{t})$.
\end{itemize}
Moreover, suppose that $\tilde{U}^{(1)},\tilde{U}^{(2)},\tilde{U}^{(3)}$ are vector fields on $B_{g(\bar{t})}(p_k,2^{\frac{j}{400}} \, \Lambda r_k)$ with the following properties: 
\begin{itemize} 
\item $|\mathscr{L}_{\tilde{U}^{(a)}}(g(\bar{t}))|+r_k \, |D(\mathscr{L}_{\tilde{U}^{(a)}}(g(\bar{t})))| \leq C \, (r_k^{-1} \, d_{g(\bar{t})}(p_k,x)+1)^{-100} \, \varepsilon_k$ for all $x \in B_{g(\bar{t})}(p_k,2^{\frac{j}{400}} \, \Lambda r_k)$. 
\item If $x \in B_{g(\bar{t})}(p_k,2^{\frac{j}{400}} \, \Lambda r_k) \setminus B_{g(\bar{t})}(p_k,4\Lambda r_k)$, then $r_k^{-1} \, |\langle \tilde{U}^{(a)},\nu \rangle| \leq C \, (r_k^{-1} \, d_{g(\bar{t})}(p_k,x)+1)^{-100} \, \varepsilon_k$, where $\nu$ denotes the unit normal to the CMC foliation of $(M,g(\bar{t}))$.
\item If $x \in B_{g(\bar{t})}(p_k,2^{\frac{j}{400}} \, \Lambda r_k) \setminus B_{g(\bar{t})}(p_k,4\Lambda r_k)$, then 
\begin{align*} 
&\sum_{a,b=1}^3 \bigg | \delta_{ab} - \text{\rm area}_{g(\bar{t})}(\Sigma)^{-2} \int_\Sigma \langle \tilde{U}^{(a)},\tilde{U}^{(b)} \rangle_{g(\bar{t})} \, d\mu_{g(\bar{t})} \bigg | \\ 
&\leq C \, (r_k^{-1} \, d_{g(\bar{t})}(p_k,x)+1)^{-100} \, \varepsilon_k, 
\end{align*}
where $\Sigma$ denotes the leaf of the CMC foliation passing through $(x,\bar{t})$.
\end{itemize}
Then there exists a $3 \times 3$ matrix $\omega \in O(3)$ such that 
\[r_k^{-1} \sum_{a=1}^3 \Big | \sum_{b=1}^3 \omega_{ab} \, U^{(b)} - \tilde{U}^{(a)} \Big |_{g(\bar{t})} \leq C \, (r_k^{-1} \, d_{g(\bar{t})}(p_k,x)+1)^{-20} \, \varepsilon_k\] 
on $B_{g(\bar{t})}(p_k,2^{\frac{j-1}{400}} \, \Lambda r_k)$.
\end{lemma}

\textbf{Proof.} 
For each integer $m \in [8\Lambda,2^{\frac{j-1}{400}} \, \Lambda]$, Proposition \ref{rigidity.of.Killing.vector.fields.on.cylinder} implies that there exists a $3 \times 3$ matrix $\omega^{(m)} \in O(3)$ such that 
\[r_k^{-1} \sum_{a=1}^3 \Big | \sum_{b=1}^3 \omega_{ab}^{(m)} \, U^{(b)} - \tilde{U}^{(a)} \Big |_{g(\bar{t})} \leq C \, m^{-80} \, \varepsilon_k\] 
on $B_{g(\bar{t})}(p_k,(m+1)r_k) \setminus B_{g(\bar{t})}(p_k,(m-1)r_k)$. Note that $|\omega^{(m)}-\omega^{(m+1)}| \leq C \, m^{-60} \, \varepsilon_k$. Consequently, there exists a $3 \times 3$ matrix $\omega \in O(3)$ such that $|\omega^{(m)}-\omega| \leq C \, m^{-40} \, \varepsilon_k$. Hence, for every integer $m \in [8\Lambda,2^{\frac{j-1}{400}} \, \Lambda]$, we obtain 
\[r_k^{-1} \sum_{a=1}^3 \Big | \sum_{b=1}^3 \omega_{ab} \, U^{(b)} - \tilde{U}^{(a)} \Big |_{g(\bar{t})} \leq C \, m^{-20} \, \varepsilon_k\] 
on $B_{g(\bar{t})}(p_k,(m+1)r_k) \setminus B_{g(\bar{t})}(p_k,(m-1)r_k)$. Using Lemma \ref{ode.for.U}, we deduce that 
\[r_k^{-1} \sum_{a=1}^3 \Big | \sum_{b=1}^3 \omega_{ab} \, U^{(b)} - \tilde{U}^{(a)} \Big |_{g(\bar{t})} \leq C \varepsilon_k\] 
on $B_{g(\bar{t})}(p_k,16\Lambda r_k)$. This completes the proof of Lemma \ref{vector.field.comparison.on.cap}. \\

In the following, we define 
\[\Omega^{(j,k)} := \{(x,t) \in B_{g(t_k)}(p_k,\delta_k^{-1} r_k) \times [t_k-2^{\frac{j}{100}} \, r_k^2,t_k]: f(x,t) \leq 2^{\frac{j}{400}}\},\] 
where $f: B_{g(t_k)}(p_k,\delta_k^{-1} r_k) \times [t_k-\delta_k^{-1} r_k^2,t_k] \to \mathbb{R}$ is the function in Corollary \ref{perturbation.of.Bryant.soliton}. We now state the main result of this section:

\begin{proposition}
\label{cap.improvement.theorem}
Let $j$ be a large positive integer. If $k$ is sufficiently large (depending on $j$), then we can find time-independent vector fields $W^{(1)},W^{(2)},W^{(3)}$ such that 
\[\sum_{l=0}^{40} r_k^l \, |D^l(\mathscr{L}_{W^{(a)}}(g))| \leq C \, 2^{-\frac{j}{400}} \, \varepsilon_k\] 
for all points $(x,t) \in B_{g(t_k)}(p_k,4\Lambda r_k) \times [t_k-1000K\Lambda r_k^2,t_k]$. Here, $C$ is a constant which is independent of $j$ and $k$. Finally, on the set $B_{g(t_k)}(p_k,4\Lambda r_k) \times [t_k-1000K\Lambda r_k^2,t_k]$, the vector fields $W^{(1)},W^{(2)},W^{(3)}$ are close to the standard rotation vector fields on the Bryant soliton in the $C^{80}$-norm.
\end{proposition}

\textbf{Proof.} 
We will assume throughout that $j$ is large, and $k$ is sufficiently large depending on $j$. This ensures, after rescaling by the factor $r_k^{-1}$, the domain $\Omega^{(j,k)}$ is close to a piece of the Bryant soliton. By Corollary \ref{perturbation.of.Bryant.soliton}, the function $f: \Omega^{(j,k)} \to \mathbb{R}$ satisfies $R-\Delta f \leq 3\delta_k r_k^{-2}$ and $\Delta f + |\nabla f|^2 \leq (1+\delta_k) r_k^{-2}$, hence 
\[R+|\nabla f|^2 \leq (1+4\delta_k) r_k^{-2} \leq 2r_k^{-2}.\] 
Moreover, Corollary \ref{perturbation.of.Bryant.soliton} implies $\Delta f + |\nabla f|^2 \geq (1-\delta_k) r_k^{-2}$ and $\frac{\partial}{\partial t} f + |\nabla f|^2 \leq \delta_k r_k^{-2}$, hence 
\[\frac{\partial}{\partial t} f - \Delta f \leq -(1-2\delta_k) r_k^{-2} \leq -\frac{1}{2} \, r_k^{-2}.\] 
Note that $\frac{1}{2K} \, (r_k^{-1} \, d_{g(t)}(p_k,x) + 1) \leq f(x,t)+1 \leq 2K \,  (r_k^{-1} \, d_{g(t)}(p_k,x) + 1)$ by Corollary \ref{perturbation.of.Bryant.soliton}.

\textit{Step 1:} Using Lemma \ref{existence.of.almost.Killing.vector.fields} and Lemma \ref{vector.field.comparison.on.cap}, we can construct time-dependent vector fields $U^{(1)},U^{(2)},U^{(3)}$, defined on $\Omega^{(j,k)}$, with the following properties: 
\begin{itemize} 
\item $r_k \, |\frac{\partial}{\partial t} U^{(a)}| \leq C \, (f+100)^{-10} \, \varepsilon_k$ on $\Omega^{(j,k)}$. 
\item $|\mathscr{L}_{U^{(a)}}(g)|+r_k \, |D(\mathscr{L}_{U^{(a)}}(g))| \leq C \, (f+100)^{-100} \, \varepsilon_k$ on $\Omega^{(j,k)}$.
\end{itemize} 
Here, $C$ is a large constant that does not depend on $j$ or $k$. Moreover, we can arrange that, on the set $\Omega^{(j,k)}$, the vector fields $U^{(1)},U^{(2)},U^{(3)}$ are close to the standard rotation vector fields on the Bryant soliton in the $C^2$-norm. Note that $r_k \, |\Delta U^{(a)} + \text{\rm Ric}(U^{(a)})| \leq Cr_k \, |D(\mathscr{L}_{U^{(a)}}(g))| \leq C \, (f+100)^{-100} \, \varepsilon_k$ on $\Omega^{(j,k)}$.

\textit{Step 2:} Let $V^{(a)}$ denote the solution of the PDE $\frac{\partial}{\partial t} V^{(a)} = \Delta V^{(a)} + \text{\rm Ric}(V^{(a)})$ on $\Omega^{(j,k)}$ with Dirichlet boundary condition $V^{(a)}=U^{(a)}$ on the parabolic boundary of $\Omega^{(j,k)}$. Using the estimate 
\[r_k \, |\frac{\partial}{\partial t} U^{(a)} - \Delta U^{(a)} - \text{\rm Ric}(U^{(a)})| \leq C \, (f+100)^{-10} \, \varepsilon_k,\] 
we obtain
\begin{align*} 
&r_k \, |\frac{\partial}{\partial t} (V^{(a)}-U^{(a)}) - \Delta (V^{(a)}-U^{(a)}) - \text{\rm Ric}(V^{(a)}-U^{(a)})| \\ 
&\leq C \, (f+100)^{-10} \, \varepsilon_k 
\end{align*}
in $\Omega^{(j,k)}$, where $C$ is a large constant that does not depend on $j$ or $k$. Proposition \ref{estimate.for.norm.of.vector.field} gives 
\[r_k \, (\frac{\partial}{\partial t}-\Delta) |V^{(a)}-U^{(a)}| \leq C \, (f+100)^{-10} \, \varepsilon_k\] 
in $\Omega^{(j,k)}$, where $C$ is a large constant that does not depend on $j$ or $k$. 

Using the inequalities $(\frac{\partial}{\partial t}-\Delta)f \leq -\frac{1}{2} \, r_k^{-2}$ and $|\nabla f|^2 \leq 2r_k^{-2}$, we obtain 
\begin{align*} 
&(\frac{\partial}{\partial t}-\Delta) (f+100)^{-8} \\ 
&= -8 \, (f+100)^{-9} \, (\frac{\partial}{\partial t} - \Delta) f - 72 \, (f+100)^{-10} \, |\nabla f|^2 \\ 
&\geq 4 \, (f+100)^{-9} \, r_k^{-2} - 144 \, (f+100)^{-10} \, r_k^{-2} \\ 
&\geq (f+100)^{-9} \, r_k^{-2}. 
\end{align*}
in $\Omega^{(j,k)}$. Using the maximum principle, we conclude that 
\[r_k^{-1} \, |V^{(a)}-U^{(a)}| \leq C \, (f+100)^{-8} \, \varepsilon_k\] 
in $\Omega^{(j,k)}$, where $C$ is a large constant that does not depend on $j$ or $k$. Using standard interior estimates for linear parabolic equations, we obtain 
\[|D(V^{(a)}-U^{(a)})| \leq C \, (f+100)^{-8} \, \varepsilon_k\] 
in $\Omega^{(j-1,k)}$, where $C$ is a large constant that does not depend on $j$ or $k$. In particular, on the set $\Omega^{(j-1,k)}$, the vector fields $V^{(1)},V^{(2)},V^{(3)}$ are close to the standard rotation vector fields on the Bryant soliton in the $C^1$-norm. Consequently, on the set $B_{g(t_k)}(p_k,8\Lambda r_k) \times [t_k-2000K\Lambda r_k^2,t_k]$, the vector fields $V^{(1)},V^{(2)},V^{(3)}$ are close to the standard rotation vector fields on the Bryant soliton in the $C^{100}$-norm. 

\textit{Step 3:} We now define $h^{(a)}(t) := \mathscr{L}_{V^{(a)}(t)}(g(t))$. Since $\frac{\partial}{\partial t} V^{(a)} = \Delta V^{(a)} + \text{\rm Ric}(V^{(a)})$, we obtain 
\[\frac{\partial}{\partial t} h^{(a)}(t) = \Delta_{L,g(t)} h^{(a)}(t)\] 
by Corollary \ref{pde.for.lie.derivative}. The estimate for $V^{(a)}-U^{(a)}$ in Step 2 implies 
\[|h^{(a)}| \leq |\mathscr{L}_{U^{(a)}}(g)| + C \, |D(V^{(a)}-U^{(a)})| \leq C \, (f+100)^{-8} \, \varepsilon_k\] 
in $\Omega^{(j-1,k)}$, where $C$ is a large constant that does not depend on $j$ or $k$.

Let $C_\#$ and $c_\#$ denote the constants in Theorem \ref{anderson.chow.estimate}. If $j$ is sufficiently large and $k$ is sufficiently large depending on $j$, then $C_\# \, 2^{-\frac{j}{200}} \, r_k^{-2} \leq R(x,t) \leq r_k^{-2}$ for all points $(x,t) \in \Omega^{(j,k)}$. Therefore, we may apply Theorem \ref{anderson.chow.estimate} with $\rho := 2^{-\frac{j}{200}} \, r_k^{-2}$. Consequently, the function 
\[\psi^{(a)} := \exp(-c_\# \, 2^{-\frac{j}{200}} \, r_k^{-2} \, (t_k-t)) \, \frac{|h^{(a)}|}{r_k^2 R-2^{-\frac{j}{200}}}\] 
satisfies 
\[\Big ( \frac{\partial}{\partial t}-\Delta -\frac{2r_k^2}{r_k^2 R-2^{-\frac{j}{200}}} \, D^i R \, D_i \Big) (\psi^{(a)})^2 \leq 0\]
on the set $\Omega^{(j-1,k)}$. Applying the maximum principle to $\psi^{(a)}$ on the set $\Omega^{(j-1,k)}$, we obtain
\[|h^{(a)}| \leq C \, 2^{-\frac{j}{400}} \, \varepsilon_k\] 
on the set $B_{g(t_k)}(p_k,16\Lambda r_k) \times [t_k-2000K\Lambda r_k^2,t_k]$, where $C$ is a large constant that does not depend on $j$ or $k$. Using standard interior estimates for linear parabolic equations, we obtain 
\[\sum_{l=0}^{100} r_k^l \, |D^l h^{(a)}| \leq C \, 2^{-\frac{j}{400}} \, \varepsilon_k\] 
on the set $B_{g(t_k)}(p_k,8\Lambda r_k) \times [t_k-1000K\Lambda r_k^2,t_k]$. To summarize, we have shown that 
\[\sum_{l=0}^{100} r_k^l \, |D^l (\mathscr{L}_{V^{(a)}}(g))| \leq C \, 2^{-\frac{j}{400}} \, \varepsilon_k\] 
on the set $B_{g(t_k)}(p_k,8\Lambda r_k) \times [t_k-1000K\Lambda r_k^2,t_k]$. Moreover, using the identity 
\[\frac{\partial}{\partial t} V^{(a)} = \Delta V^{(a)} + \text{\rm Ric}(V^{(a)}) = \text{\rm div} \, h^{(a)} - \frac{1}{2} \, \nabla(\text{\rm tr} \, h^{(a)}),\] 
we obtain 
\[\sum_{l=0}^{80} r_k^{l+1} \, \Big | D^l \Big ( \frac{\partial}{\partial t} V^{(a)} \Big ) \Big | \leq C \, 2^{-\frac{j}{400}} \, \varepsilon_k,\] 
on the set $B_{g(t_k)}(p_k,8\Lambda r_k) \times [t_k-1000K\Lambda r_k^2,t_k]$. 

\textit{Step 4:} Let $W^{(1)},W^{(2)},W^{(3)}$ be time-independent vector fields such that $W^{(a)}=V^{(a)}$ at time $t_k$. On the set $B_{g(t_k)}(p_k,8\Lambda r_k) \times [t_k-1000K\Lambda r_k^2,t_k]$, the vector fields $W^{(1)},W^{(2)},W^{(3)}$ are close to the standard rotation vector fields on the Bryant soliton. Using the estimate for $\frac{\partial}{\partial t} V^{(a)}$ in Step 3, we obtain  
\[\sum_{l=0}^{60} r_k^{l-1} \, |D^l(W^{(a)}-V^{(a)})| \leq C \, 2^{-\frac{j}{400}} \, \varepsilon_k,\] 
hence 
\[\sum_{l=0}^{40} r_k^l \, |D^l(\mathscr{L}_{W^{(a)}}(g))| \leq C \, 2^{-\frac{j}{400}} \, \varepsilon_k\] 
on the set $B_{g(t_k)}(p_k,8\Lambda r_k) \times [t_k-1000K\Lambda r_k^2,t_k]$. This completes the proof of Proposition \ref{cap.improvement.theorem}. \\

For each $k$ large, we choose a compact domain $D_k \subset M$ with the following properties: 
\begin{itemize} 
\item There exists a point $x \in \partial D_k$ such that $\lambda_1(x,t_k) = \frac{2}{3} \theta R(x,t_k)$. 
\item For each $x \in D_k$, we have $\lambda_1(x,t_k) \geq \frac{2}{3} \theta R(x,t_k)$. 
\item $\partial D_k$ is a leaf of the CMC foliation of $(M,g(t_k))$.
\end{itemize} 
Note that $D_k \subset \{x \in M: \lambda_1(x,t_k) > \frac{1}{2} \theta R(x,t_k)\} \subset B_{g(t_k)}(p_k,\Lambda r_k)$ in view of our choice of $\Lambda$. Moreover, if $\bar{x} \in M \setminus D_k$, then the point $(\bar{x},t_k)$ lies at the center of an evolving $\varepsilon_1$-neck by Proposition \ref{choice.of.theta}.

\begin{proposition}
\label{W.tangential.to.CMC.foliation}
Let $j$ be a large positive integer. If $k$ is sufficiently large (depending on $j$), then the vector fields $W^{(1)},W^{(2)},W^{(3)}$ constructed in Proposition \ref{cap.improvement.theorem} have the following property. For each point $\bar{x} \in B_{g(t_k)}(p_k,\Lambda r_k) \setminus D_k$, we have 
\[\sum_{l=0}^{10} r_k^{l-1} \, |D^l(\langle W^{(a)},\nu \rangle)| \leq C \, 2^{-\frac{j}{400}} \, \varepsilon_k\] 
on the parabolic neighborhood $B_{g(t_k)}(\bar{x},600 R(\bar{x},t_k)^{-\frac{1}{2}}) \times [t_k-200 R(\bar{x},t_k)^{-1},t_k]$. Here, $\nu$ denotes the unit normal to the CMC foliation and $C$ is a constant which is independent of $j$ and $k$. 
\end{proposition}

\textbf{Proof.} 
Let us consider a point $\bar{x} \in B_{g(t_k)}(p_k,\Lambda r_k) \setminus D_k$. Recall that the point $(\bar{x},t_k)$ lies at the center of an evolving $\varepsilon_1$-neck. By Corollary \ref{perturbation.of.Bryant.soliton}, $R(\bar{x},t_k)^{-1} \leq 4K\Lambda r_k^2$. Since $\sqrt{\frac{4K}{\Lambda}} \leq 10^{-6}$, we obtain $R(\bar{x},t_k)^{-\frac{1}{2}} \leq 10^{-6} \Lambda r_k$. Hence, the parabolic neighborhood $B_{g(t_k)}(\bar{x},1000 R(\bar{x},t_k)^{-\frac{1}{2}}) \times [t_k-200 R(\bar{x},t_k)^{-1},t_k]$ is contained in $B_{g(t_k)}(p_k,4\Lambda r_k) \times [t_k-1000K\Lambda r_k^2,t_k]$. In particular, the estimates in Proposition \ref{cap.improvement.theorem} hold on the parabolic neighborhood $B_{g(t_k)}(\bar{x},1000 R(\bar{x},t_k)^{-\frac{1}{2}}) \times [t_k-200 R(\bar{x},t_k)^{-1},t_k]$. 

Let us fix a time $t \in [t_k-200 R(\bar{x},t_k)^{-1},t_k]$, and let $\Sigma_s$ denote the CMC foliation of $(M,g(t))$. Note that the foliation depends on $t$, but we suppress this dependence in our notation. In the following, we only consider those leaves of the foliation that are contained in $B_{g(t_k)}(\bar{x},800 R(\bar{x},t_k)^{-\frac{1}{2}})$. We define a function $F^{(a)}: \Sigma_s \to \mathbb{R}$ by $F^{(a)} := \langle W^{(a)},\nu \rangle$. The quantity 
\[\Delta_{\Sigma_s} F^{(a)} +  (|A|^2+\text{\rm Ric}(\nu,\nu)) \, F^{(a)} =: H^{(a)}\] 
can be expressed in terms of $\mathscr{L}_{W^{(a)}}(g)$ and the first derivatives of $\mathscr{L}_{W^{(a)}}(g)$. Using the estimate for $\mathscr{L}_{W^{(a)}}(g)$ in Step 4, we obtain 
\[\sum_{l=0}^{20} r_k^{l+1} \, |D^l H^{(a)}| \leq C \, 2^{-\frac{j}{400}} \, \varepsilon_k.\] 
We now consider the spectrum of the Jacobi operator $\Delta_{\Sigma_s} + (|A|^2+\text{\rm Ric}(\nu,\nu))$ on $\Sigma_s$. Since $\Sigma_s \subset B_{g(t_k)}(\bar{x},800 R(\bar{x},t_k)^{-\frac{1}{2}}) \subset B_{g(t_k)}(p_k,4\Lambda r_k)$, the eigenvalues of the Jacobi operator lie outside the interval $[-cr_k^{-2},cr_k^{-2}]$ for some small positive constant $c$ which is independent of $j$ and $k$. (This can be easily verified on the Bryant soliton. For the general case, we observe that the actual solution is a small perturbation of the Bryant soliton in the relevant region.) Consequently, we can invert the Jacobi operator. Using the estimate $\sum_{l=0}^{20} r_k^{l+1} \, |D^l H^{(a)}| \leq C \, 2^{-\frac{j}{400}} \, \varepsilon_k$, we obtain 
\[\sum_{l=0}^{10} r_k^{l-1} \, |D^l F^{(a)}| \leq C \, 2^{-\frac{j}{400}} \, \varepsilon_k.\] 
Since $t \in [t_k-200 R(\bar{x},t_k)^{-1},t_k]$ is arbitrary, we conclude that 
\[\sum_{l=0}^{10} r_k^{l-1} \, |D^l(\langle W^{(a)},\nu \rangle)| \leq C \, 2^{-\frac{j}{400}} \, \varepsilon_k\] 
on $B_{g(t_k)}(\bar{x},600 R(\bar{x},t_k)^{-\frac{1}{2}}) \times [t_k-200 R(\bar{x},t_k)^{-1},t_k]$. This completes the proof of Proposition \ref{W.tangential.to.CMC.foliation}. \\

\begin{proposition}
\label{L2.orthonormality}
Let $j$ be a large integer. If $k$ is sufficiently large (depending on $j$), then the vector fields $W^{(1)},W^{(2)},W^{(3)}$ constructed in Proposition \ref{cap.improvement.theorem} have the following property. For each point $\bar{x} \in B_{g(t_k)}(p_k,\Lambda r_k) \setminus D_k$, we can find a symmetric $3 \times 3$ matrix $Q_{ab}$ such that 
\[\bigg | Q_{ab} - \text{\rm area}_{g(t)}(\Sigma)^{-2} \int_\Sigma \langle W^{(a)},W^{(b)} \rangle_{g(t)} \, d\mu_{g(t)} \bigg | \leq C \, 2^{-\frac{j}{400}} \, \varepsilon_k\] 
whenever $t \in [t_k - 200 R(\bar{x},t_k)^{-1},t_k]$ and $\Sigma \subset B_{g(t_k)}(\bar{x},200 R(\bar{x},t_k)^{-\frac{1}{2}})$ is a leaf of the CMC foliation of $(M,g(t))$. Note that $Q_{ab}$ is independent of $t$ and $\Sigma$. Moreover, the eigenvalues of the matrix $Q_{ab}$ lie in the interval $[\frac{1}{C},C]$, where $C$ is independent of $j$ and $k$.
\end{proposition} 

\textbf{Proof.} 
Let us consider a point $\bar{x} \in B_{g(t_k)}(p_k,\Lambda r_k) \setminus D_k$. Recall that the point $(\bar{x},t_k)$ lies at the center of an evolving $\varepsilon_1$-neck. Moreover, since $R(\bar{x},t_k)^{-1} \leq 4K\Lambda r_k^2$ and $R(\bar{x},t_k)^{-\frac{1}{2}} \leq 10^{-6} \Lambda r_k$, the parabolic neighborhood $B_{g(t_k)}(\bar{x},600 R(\bar{x},t_k)^{-\frac{1}{2}}) \times [t_k-200 R(\bar{x},t_k)^{-1},t_k]$ is contained in $B_{g(t_k)}(p_k,4\Lambda r_k) \times [t_k-1000K\Lambda r_k^2,t_k]$. Therefore, the estimates in Proposition \ref{cap.improvement.theorem} and Proposition \ref{W.tangential.to.CMC.foliation} hold on the parabolic neighborhood $B_{g(t_k)}(\bar{x},600 R(\bar{x},t_k)^{-\frac{1}{2}}) \times [t_k-200 R(\bar{x},t_k)^{-1},t_k]$. 

We now argue as in Steps 11 -- 14 in the proof of the Neck Improvement Theorem. This implies that there exists a symmetric $3 \times 3$ matrix $Q_{ab}$ such that  
\[\bigg | Q_{ab} - \text{\rm area}_{g(t)}(\Sigma)^{-2} \int_\Sigma \langle W^{(a)},W^{(b)} \rangle_{g(t)} \, d\mu_{g(t)} \bigg | \leq C \, 2^{-\frac{j}{400}} \, \varepsilon_k\] 
whenever $t \in [t_k - 200 R(\bar{x},t_k)^{-1},t_k]$ and $\Sigma \subset B_{g(t_k)}(\bar{x},200 R(\bar{x},t_k)^{-\frac{1}{2}})$ is a leaf of the CMC foliation of $(M,g(t))$. Finally, since the vector fields $W^{(1)},W^{(2)},W^{(3)}$ are close to the standard rotation vector fields on the Bryant soliton, the eigenvalues of the matrix $Q_{ab}$ are uniformly bounded from above and below. This completes the proof of Proposition \ref{L2.orthonormality}. \\

\begin{corollary}
\label{improved.symmetry.outside.D_k}
If $k$ is sufficiently large, then $(\bar{x},t_k)$ is $\frac{\varepsilon_k}{2}$-symmetric for all $\bar{x} \in B_{g(t_k)}(p_k,\Lambda r_k) \setminus D_k$. 
\end{corollary}

\textbf{Proof.} 
This follows by combining Proposition \ref{cap.improvement.theorem}, Proposition \ref{W.tangential.to.CMC.foliation}, and Proposition \ref{L2.orthonormality}. \\

We can now complete the proof of Theorem \ref{main.thm.b}. Combining Lemma \ref{improved.symmetry.for.faraway.points} and Corollary \ref{improved.symmetry.outside.D_k}, we conclude that $(\bar{x},t_k)$ is $\frac{\varepsilon_k}{2}$-symmetric for all $\bar{x} \in M \setminus D_k$. Moreover, if $k$ is sufficiently large, it follows from Proposition \ref{cap.improvement.theorem}, Proposition \ref{W.tangential.to.CMC.foliation}, and Proposition \ref{L2.orthonormality} that there exist vector fields on the cap $D_k$ which satisfy the requirements of Definition \ref{symmetry.of.cap} with $\varepsilon = \frac{\varepsilon_k}{2}$. Therefore, the flow is $\frac{\varepsilon_k}{2}$-symmetric at time $t_k$ if $k$ is sufficiently large. By Lemma \ref{openness.2}, we can find a time $\tilde{t}_k>t_k$ with the property that the flow is $\varepsilon_k$-symmetric at time $t$ for all $t \in [t_k,\tilde{t}_k]$. This contradicts the definition of $t_k$. This completes the proof of Theorem \ref{main.thm.b}.

\appendix 

\section{Summary of known results about ancient $\kappa$-solutions}

In this appendix, we collect some of the main known results on ancient $\kappa$-solutions, which we use in this paper. We first recall a basic Riemannian geometry fact:

\begin{proposition}
\label{no.small.necks}
Let $(M,g)$ be a complete, noncompact manifold with positive sectional curvature, and let $N$ be a neck in $M$. Let $U$ denote the unbounded connected component of $M \setminus N$. If every point in $N \cup U$ lies at the center of a neck, then $\sup_U R \leq C \sup_N R$.
\end{proposition}

\textbf{Proof.} 
The assertion is a consequence of Corollary 2.21 in \cite{Morgan-Tian}. (Note that the soul cannot lie at the center of a neck, and therefore must be contained in $M \setminus (N \cup U)$.) 

In the following, we give an alternative argument for the convenience of the reader. By assumption, every point in $N \cup U$ lies at the center of a neck. Hence, by work of Hamilton, there is a canonical CMC foliation $\Sigma_s$, $s \in [0,\infty)$, such that $\Sigma_0 \subset N$, and $U \subset \bigcup_{s \in [0,\infty)} \Sigma_s \subset N \cup U$. Let $v$ denote the lapse function of this CMC foliation. We assume that the lapse function $v$ has mean value $1$, so that $\int_{\Sigma_s} v = \text{\rm area}(\Sigma_s)$ for each $s$. Note that $\sup_{\Sigma_s} |v-1|$ is very small; in particular, $v$ is positive. We next compute 
\[-\frac{d}{ds} H = \Delta_{\Sigma_s} v + (|A|^2+\text{\rm Ric}(\nu,\nu)) \, v \geq \Delta_{\Sigma_s} v + \frac{1}{2} \, H^2 \, v\] 
at each point on $\Sigma_s$. We now take the mean value over $\Sigma_s$. Clearly, $\Delta_{\Sigma_s} v$ has mean value $0$ by the divergence theorem. Moreover, since $H$ is constant on $\Sigma_s$ and $v$ has mean value $1$, it follows that the function $H^2 \, v$ has mean value $H^2$. This gives $-\frac{d}{ds} H \geq \frac{1}{2} \, H^2$. Hence, if $H(s)<0$ for some $s$, then $H(s)$ converges to $-\infty$ at a finite value of $s$, which is impossible. Therefore, $H(s) \geq 0$ for all $s$. Consequently, $\text{\rm area}(\Sigma_s)$ is an increasing function of $s$. This implies $\frac{1}{C} \sup_{\Sigma_s} R \leq \text{\rm area}(\Sigma_s)^{-1} \leq \text{\rm area}(\Sigma_0)^{-1} \leq C \sup_{\Sigma_0} R$ for all $s \in [0,\infty)$. From this, the assertion follows. \\

We now recall the following fundamental theorem due to Perelman:

\begin{theorem}[G.~Perelman \cite{Perelman1}, Section 11.8]
\label{canonical.neighborhood.theorem}
Let $(M,g(t))$ be a three-dimensional ancient $\kappa$-solution which is noncompact and has positive sectional curvature. Given any $\varepsilon>0$, we can find a compact domain $\Omega_t \subset M$ with the following properties: 
\begin{itemize}
\item For each $x \in M \setminus \Omega_t$, the point $(x,t)$ lies at the center of an evolving $\varepsilon$-neck.
\item The boundary $\partial \Omega_t$ is a leaf of the CMC foliation at time $t$.
\item $\sup_{x \in \Omega_t} R(x,t) \leq C(\varepsilon) \inf_{x \in \Omega_t} R(x,t)$. 
\item $\text{\rm diam}_{g(t)}(\Omega_t)^2 \, \sup_{x \in \Omega_t} R(x,t) \leq C(\varepsilon)$. 
\end{itemize}
\end{theorem}

Combining Theorem \ref{canonical.neighborhood.theorem} with Proposition \ref{no.small.necks} gives:

\begin{corollary} 
\label{bound.for.R_max}
Let $(M,g(t))$ be a three-dimensional ancient $\kappa$-solution which is noncompact and has positive sectional curvature. Let $\varepsilon$ be a small positive real number, and let $\Omega_t$ be as in Theorem \ref{canonical.neighborhood.theorem}. Then $\sup_{x \in M} R(x,t) \leq C(\varepsilon) \inf_{x \in \Omega_t} R(x,t)$ and $\text{\rm diam}_{g(t)}(\Omega_t)^2 \, \sup_{x \in M} R(x,t) \leq C(\varepsilon)$. 
\end{corollary}

\textbf{Proof.} 
Proposition \ref{no.small.necks} implies that $\sup_{x \in M \setminus \Omega_t} R(x,t) \leq C \sup_{x \in \partial \Omega_t} R(x,t)$. This gives $\sup_{x \in M} R(x,t) \leq C \sup_{x \in \Omega_t} R(x,t)$. Hence, the assertion follows from Theorem \ref{canonical.neighborhood.theorem}. \\



The next result is a consequence of the Neck Stability Theorem of Kleiner and Lott: 

\begin{theorem}[cf. Kleiner-Lott \cite{Kleiner-Lott}, Theorem 6.1]
\label{consequence.of.kleiner.lott}
Let $(M,g(t))$ be a three-dimensional ancient $\kappa$-solution which is noncompact and has positive sectional curvature. Then there exists a point $q \in M$ such that $\sup_{t \leq 0} (-t) \, R(q,t) \leq 100$.
\end{theorem}

\textbf{Proof.} 
In the following, we give a proof for the convenience of the reader. Suppose that the assertion is false, so that $\sup_{t \leq 0} (-t) \, R(q,t) > 100$ for each point $q \in M$. Let $q_k$ be a sequence of points going to infinity. For each $k$, we denote by $\ell_k(x,t)$ the reduced distance of $(x,t)$ from $(q_k,0)$. Moreover, we denote by 
\[V_k(t) = (-t)^{-\frac{3}{2}} \int_M e^{-\ell_k(x,t)} \, d\text{\rm vol}_{g(t)}\] 
the reduced volume at time $t$. 

By work of Perelman \cite{Perelman1}, we can find a sequence $\varepsilon_k \to 0$ such that the point $(q_k,0)$ lies at the center of an evolving $\varepsilon_k$-neck (cf. Theorem \ref{canonical.neighborhood.theorem}). This implies $(-t) \, R(x,t) \leq 10$ for all $t \in [-\varepsilon_k^{-1} R(q_k,0)^{-1},0]$. Therefore, $\ell_k(q_k,t) \leq 100$ for all $t \in [-\varepsilon_k^{-1} R(q_k,0)^{-1},0)$. By a result of Ye, there exists a universal constant $C$ such that 
\[\frac{d_{g(t)}(x,y)^2}{(-t)} \leq C \, (\ell_k(x,t)+\ell_k(y,t)+1)\] 
for all $t < 0$ and all $x,y \in M$ (see \cite{Ye}, Lemma 3.2). Putting $y=q_k$ gives 
\[\frac{d_{g(t)}(x,q_k)^2}{(-t)} \leq C \, (\ell_k(x,t)+1)\] 
for all $t \in [-\varepsilon_k^{-1} R(q_k,0)^{-1},0)$ and all $x \in M$.

Recall that the point $(q_k,0)$ lies at the center of an evolving $\varepsilon_k$-neck. Using this fact together with Ye's estimate, we obtain 
\[\limsup_{k \to \infty} V_k(\tau \, R(q_k,0)^{-1}) \leq V_{\text{\rm cyl}}(\tau)\] 
for each $\tau < 0$, where $V_{\text{\rm cyl}}(\tau)$ denotes the reduced volume for a family of shrinking cylinders. Using the monotonicity of the reduced volume, we deduce that 
\[\limsup_{k \to \infty} V_k(-\varepsilon_k^{-1} R(q_k,0)^{-1}) \leq V_{\text{\rm cyl}}(\tau)\] 
for each $\tau < 0$. Taking the limit as $\tau \to -\infty$ gives 
\[\limsup_{k \to \infty} V_k(-\varepsilon_k^{-1} R(q_k,0)^{-1}) \leq V_{\text{\rm cyl}}(-\infty),\] 
where $V_{\text{\rm cyl}}(-\infty) := \lim_{\tau \to -\infty} V_{\text{\rm cyl}}(\tau)$. On the other hand, since the asymptotic shrinking soliton is a family of shrinking cylinders, we have 
\[\lim_{t \to -\infty} V_k(t) \geq V_{\text{\rm cyl}}(-\infty)\] 
for each $k$. Since $V_k(t)$ is monotone increasing in $t$, it follows that  
\[V_k(t) \geq V_{\text{\rm cyl}}(-\infty)\] 
for all $k$ and all $t$.

For each $k$, we define $t_k := \sup \{t \leq 0: (-t) \, R(q_k,t) > 10\}$. Clearly, $t_k \leq -\varepsilon_k^{-1} R(q_k,0)^{-1}$, $(-t_k) \, R(q_k,t_k) = 10$, and $(-t) \, R(q_k,t) \leq 10$ for all $t \in [t_k,0]$. This implies $\ell_k(q_k,t_k) \leq \frac{1}{2\sqrt{-t_k}} \int_{t_k}^0 \sqrt{-t} \, R(q_k,t) \, dt \leq 100$. The discussion above gives 
\[\inf_{t \in (-\infty,t_k]} V_k(t) \geq V_{\text{\rm cyl}}(-\infty)\] 
and 
\[\sup_{t \in (-\infty,t_k]} V_k(t) \leq V_k(-\varepsilon_k^{-1} R(q_k,0)^{-1}) \to V_{\text{\rm cyl}}(-\infty)\]
as $k \to \infty$. Hence, if we dilate the flow $(M,g(t))$ around $(q_k,t_k)$ by the factor $(-t_k)^{-\frac{1}{2}}$, then the rescaled flows converge in the Cheeger-Gromov sense to a shrinking gradient Ricci soliton (see \cite{Perelman1}, Section 11). Since $(-t_k) \, R(q_k,t_k) = 10$ for each $k$, this limiting soliton is non-flat, and consequently must be a cylinder with scalar curvature $1$ (cf. \cite{Perelman2}, Section 1). In particular, $(-t_k) \, R(q_k,t_k) \to 1$ as $k \to \infty$. This contradicts the fact that $(-t_k) \, R(q_k,t_k) = 10$ for each $k$. This completes the proof of Theorem \ref{consequence.of.kleiner.lott}. \\

\section{A variant of the Anderson-Chow estimate}

\label{anderson.chow.est}

In \cite{Anderson-Chow}, Anderson and Chow proved an important estimate for solutions of the parabolic Lichnerowicz equation. In this appendix, we state a variant of that estimate which is due to Kyeongsu Choi:

\begin{theorem}[K.~Choi]
\label{anderson.chow.estimate}
There exists a large constant $C_\# \geq 10$ and a small positive constant $c_\#$ such that the following holds. Let $(M,g(t))$ be a solution to the Ricci flow in dimension $3$ with nonnegative Ricci curvature, let $h(t)$ be a one-parameter family of symmetric $(0,2)$-tensors satisfying the parabolic Lichnerowicz equation $\frac{\partial}{\partial t} h(t) = \Delta_{L,g(t)} h(t)$, and let $\rho$ denote a positive real number. Then 
\[\Big (\frac{\partial}{\partial t}-\Delta -\frac{2}{R-\rho} \, D^i R \, D_i \Big ) \, \Big ( \exp(2c_\# \rho t) \, \frac{|h|^2}{(R-\rho)^2} \Big ) \leq 0\] 
whenever $R \geq C_\# \rho$.
\end{theorem}

In the following, we sketch the proof of Theorem \ref{anderson.chow.estimate}. We assume throughout that $R>\rho>0$. The computation of Anderson-Chow yields
\begin{align*}
&\Big ( \frac{\partial}{\partial t}-\Delta -\frac{2}{R-\rho} \, D^i R \, D_i \Big) \frac{|h|^2}{(R-\rho)^2} \\ 
&= -\frac{2}{(R-\rho)^4} \, |(R-\rho) \, D_i h_{jk} - D_i R \, h_{jk}|^2 - \frac{4S}{(R-\rho)^3}, 
\end{align*} 
where 
\begin{align*} 
S &:= -2(R-\rho) \, \langle \text{\rm Ric},h \rangle \, \text{\rm tr}(h) + 2(R-\rho) \, \langle \text{\rm Ric},h^2 \rangle \\ 
&- \frac{1}{2} \, R(R-\rho) \, (|h|^2 - \text{\rm tr}(h)^2) + |h|^2 \, |\text{\rm Ric}|^2 
\end{align*} 
(cf. \cite{Anderson-Chow}, p.~8). Let us fix a point $p \in M$, and consider an orthonormal basis of $T_p M$ with the property that $h$ is diagonal. We denote by $h_1,h_2,h_3$ the diagonal entries of $h$. Moreover, we denote by $r_1,r_2,r_3$ the diagonal entries of $\text{\rm Ric}$. We may assume that $r_1 \leq r_2 \leq r_3$. Clearly, $R=r_1+r_2+r_3$ and $|\text{\rm Ric}|^2 \geq r_1^2+r_2^2+r_3^2$. This implies 
\[2S \geq \begin{bmatrix} h_1 & h_2 & h_3 \end{bmatrix} A_\rho  \begin{bmatrix} h_1 \\ h_2 \\ h_3 \end{bmatrix},\] 
where $A_\rho$ is defined by
\[A_\rho = 
\begin{bmatrix} 
2(r_1^2+r_2^2+r_3^2) & (\rho-R)(r_1+r_2-r_3) & (\rho-R)(r_3+r_1-r_2) \\ 
(\rho-R)(r_1+r_2-r_3) & 2(r_1^2+r_2^2+r_3^2) & (\rho-R)(r_2+r_3-r_1) \\ 
(\rho-R)(r_3+r_1-r_2) & (\rho-R)(r_2+r_3-r_1) & 2(r_1^2+r_2^2+r_3^2) 
\end{bmatrix}.\]
We claim that the matrix $A_\rho$ is positive definite. To prove this, we use Sylvester's criterion. The first minor is clearly positive. The second minor satisfies
\begin{align*}
&4(r_1^2+r_2^2+r_3^2)^2-(\rho-R)^2(r_1+r_2-r_3)^2 \\ 
&\geq 4(r_1^2+r_2^2+r_3^2)^2-R^2(r_1+r_2-r_3)^2 \\ 
&= 4(r_1^2+r_2^2+r_3^2)^2-((r_1+r_2)^2-r_3^2)^2 \\ 
&> 0, 
\end{align*}
where in the last step we have used the inequality $-2(r_1^2+r_2^2+r_3^2) < (r_1+r_2)^2 - r_3^2 < 2(r_1^2+r_2^2+r_3^2)$. 

Finally, we consider the third minor of $A_\rho$. Expanding $\det A_\rho$ in powers of $\rho$ gives 
\begin{align*}
\det A_\rho 
&\geq \det A_0 \\ 
&+ 4\rho R (r_1^2+r_2^2+r_3^2) \, [(r_1+r_2-r_3)^2+(r_2+r_3-r_1)^2+(r_3+r_1-r_2)^2] \\ 
&+ 6\rho R^2  (r_1+r_2-r_3)(r_2+r_3-r_1)(r_3+r_1-r_2) \\ 
&- C \rho^2 R^4 - C \rho^3  R^3. 
\end{align*} 
By work of Anderson-Chow, $\det A_0 \geq 0$ (see \cite{Anderson-Chow}, pp.~10--11). Moreover, 
\[R(r_1+r_2-r_3) = (r_1+r_2)^2-r_3^2 \geq -r_3^2\] 
and 
\[0 \leq (r_2+r_3-r_1)(r_3+r_1-r_2) \leq \frac{1}{2} \, [(r_2+r_3-r_1)^2+(r_3+r_1-r_2)^2].\] 
This implies 
\[R(r_1+r_2-r_3) (r_2+r_3-r_1)(r_3+r_1-r_2) \geq -\frac{1}{2} \, r_3^2 \, [(r_2+r_3-r_1)^2+(r_3+r_1-r_2)^2].\] 
Putting these facts together, we obtain 
\begin{align*}
\det A_\rho 
&\geq 4\rho R (r_1^2+r_2^2+r_3^2) \, [(r_1+r_2-r_3)^2+(r_2+r_3-r_1)^2+(r_3+r_1-r_2)^2] \\ 
&- 3\rho R r_3^2 \, [(r_2+r_3-r_1)^2+(r_3+r_1-r_2)^2] \\ 
&- C \rho^2 R^4 - C \rho^3  R^3. 
\end{align*} 
Hence, we can find a large constant $C_\# \geq 10$ and a small positive constant $c$ with the property that $\det A_\rho \geq c \rho R^5 > 0$ whenever $R \geq C_\# \rho$. By Sylvester's criterion, the matrix $A_\rho$ is positive definite whenever $R \geq C_\# \rho$. Moreover, the largest eigenvalue of $A_\rho$ is bounded by $C R^2$ whenever $R > \rho$. Since $\det A_\rho \geq c \rho R^5$, it follows that the smallest eigenvalue of $A_\rho$ is greater than $c \rho R$ whenever $R \geq C_\# \rho$.

To summarize, we have shown that there exists a small positive constant $c_\#$ such that $2S \geq c_\# \rho R \, |h|^2$ whenever $R \geq C_\# \rho$. Putting these facts together, we conclude that 
\[\Big (\frac{\partial}{\partial t}-\Delta -\frac{2}{R-\rho} \, D^i R \, D_i \Big ) \, \Big ( \frac{|h|^2}{(R-\rho)^2} \Big ) \leq -\frac{4S}{(R-\rho)^3} \leq -2c_\# \rho \, \frac{|h|^2}{(R-\rho)^2}\] 
whenever $R \geq C_\# \rho$. From this, the assertion follows. \\


\end{document}